\documentclass[11pt,a4paper]{article}
\usepackage{pdfsync} 

\usepackage{amsfonts,amsmath,amsthm,enumitem,tikz,pgfplots,subfig}
\usetikzlibrary{positioning}
\newtheorem{theorem}{Theorem}
\numberwithin{theorem}{section}
\newtheorem{lemma}[theorem]{Lemma}
\theoremstyle{definition}

\numberwithin{equation}{section}
\usepackage{amssymb}
\usepackage{hyperref}
\usepackage{graphicx}
\usepackage{wrapfig}
\usepackage{subfig}
\usepackage{epsfig}
\usepackage{float,mathrsfs}
\usepackage{bbm}

\usepackage{algpseudocode}
\usepackage[font=scriptsize]{caption}
\counterwithout{figure}{section}
\usepackage{color,xcolor}
\definecolor{refkey}{rgb}{0.9451,0.2706,0.4941}
\definecolor{labelkey}{rgb}{0.9451,0.2706,0.4941}

\definecolor{darkred}{RGB}{139,0,0}
\definecolor{darkgreen}{RGB}{0,100,0}
\definecolor{darkmagenta}{RGB}{139,0,139}

\newcommand{\setu}{{\mathrm{\mathfrak{u}}}}
\newcommand{\setv}{{\mathrm{\mathfrak{v}}}}

\newcommand{\bsb}{{\boldsymbol{b}}}
\newcommand{\bse}{{\boldsymbol{e}}}

\newcommand{\bsm}{{\boldsymbol{m}}}

\newcommand{\bsw}{{\boldsymbol{w}}}
\newcommand{\bsx}{{\boldsymbol{x}}}
\newcommand{\bsy}{{\boldsymbol{y}}}
\newcommand{\bsz}{{\boldsymbol{z}}}

\newcommand{\bsV}{{\boldsymbol{V}}}

\newcommand{\bsell}{{\boldsymbol{\ell}}}
\newcommand{\bsalpha}{{\boldsymbol{\alpha}}}
\newcommand{\bsbeta}{{\boldsymbol{\beta}}}
\newcommand{\bsgamma}{{\boldsymbol{\gamma}}}
\newcommand{\bslambda}{{\boldsymbol{\lambda}}}

\newcommand{\bsmu}{{\boldsymbol{\mu}}}

\newcommand{\bspsi}{{\boldsymbol{\psi}}}
\newcommand{\bsnu}{{\boldsymbol{\nu}}}

\newcommand{\bszero}{{\boldsymbol{0}}}

\newcommand{\rd}{{\mathrm{d}}}
\newcommand{\ri}{{\mathrm{i}}}
\newcommand{\re}{{\mathrm{e}}}

\newcommand{\bbA}{{\mathbb{A}}}
\newcommand{\bbB}{{\mathbb{B}}}

\newcommand{\bbP}{{\mathbb{P}}}

\newcommand{\bbZ}{{\mathbb{Z}}}

\newcommand{\bbR}{{\mathbb{R}}}

\newcommand{\calO}{\mathcal{O}}

\newcommand{\calC}{\mathcal{C}}

\newcommand{\calF}{\mathcal{F}}

\newcommand{\tr}{{\rm T}}
\newcommand{\satop}[2]{\stackrel{\scriptstyle{#1}}{\scriptstyle{#2}}}
\newcommand{\esup}{\operatornamewithlimits{ess\,sup}}

\newcommand{\mask}[1]{{}}
\newcommand{\Dref}{{D_{\rm ref}}}
\renewcommand{\top}{{\rm T}}
\newcommand{\supp}{{\rm supp}}

\usepackage{subfig}
\graphicspath{{Figures/}{PDF/}}

\title{Uncertainty quantification for random domains\\using periodic random variables}

\author{H. Hakula\footnotemark[1]\and H. Harbrecht\footnotemark[2]\and V. Kaarnioja\footnotemark[3]\and  F.~Y. Kuo\footnotemark[4]\and I.~H. Sloan\footnotemark[4]}

\begin{document}
\footnotetext[1]{Department of Mathematics and Systems Analysis, Aalto University, P.O. Box 11100, 00076 Aalto, Finland
 ({\tt harri.hakula@aalto.fi}).}
\footnotetext[2]{Departement Mathematik und Informatik, Universit\"{a}t Basel, Spiegelgasse 1, 4051 Basel, Switzerland
 ({\tt helmut.harbrecht@unibas.ch}).}
\footnotetext[3]{Fachbereich Mathematik und Informatik, Freie Universit\"at Berlin, Arnimallee 6, 14195 Berlin, Germany
 ({\tt vesa.kaarnioja@fu-berlin.de}).}
 \footnotetext[4]{School of Mathematics and Statistics, University of New South Wales, Sydney NSW 2052, Australia
 ({\tt f.kuo@unsw.edu.au}, {\tt i.sloan@unsw.edu.au}).}
\maketitle

\begin{abstract}
We consider uncertainty quantification for the Poisson problem subject to domain uncertainty. For the stochastic parameterization of the random domain, we use the model recently introduced by Kaarnioja, Kuo, and Sloan (SIAM J.~Numer.~Anal., 2020) in which a countably infinite number of independent random variables enter the random field as periodic functions. We develop lattice quasi-Monte Carlo (QMC) cubature rules for computing the expected value of the solution to the Poisson problem subject to domain uncertainty. These QMC rules can be shown to exhibit higher order cubature convergence rates permitted by the periodic setting independently of the stochastic dimension of the problem. In addition, we present a complete error analysis for the problem by taking into account the approximation errors incurred by truncating the input random field to a finite number of terms and discretizing the spatial domain using finite elements. The paper concludes with numerical experiments demonstrating the theoretical error estimates.
\end{abstract}

\section{Introduction}

Modeling domain uncertainty is pertinent in many engineering
applications, where the shape of the object may not be perfectly
known. For example, one might think of manufacturing imperfections in the
shape of products fabricated by line production, or shapes which stem from
inverse problems such as tomography.

If the domain perturbations are small, one can apply the {\em perturbation
technique} as in, e.g.,~\cite{BC02,H10,HSS08}. In that approach,
{which is based on Eulerian coordinates,} the shape derivative
is used to linearize the problem under consideration around a
nominal reference geometry. By Hadamard's theorem, the resulting
linearized equations for the first order shape sensitivities are
homogeneous equations posed on the nominal geometry, with inhomogeneous
boundary data only. By using a {first order} shape Taylor expansion, one can
derive tensor deterministic partial differential equations {(PDEs)} for the
statistical moments of the problem.

The approach we consider in this work is the {\em domain mapping}
approach. It transfers the shape uncertainty onto a fixed reference
domain {and hence reflects the Lagrangian setting.}
Especially, it allows us to deal with large deformations, see~\cite{mohan,xiu} for example. Analytic dependency of the solution on the
random domain mapping in the case of the Poisson equation was shown
in~\cite{CNT16,nummat2016} and in the case of linear elasticity in~\cite{HKS23}.
Related shape holomorphy for acoustic and electromagnetic scattering
problems has been shown in~\cite{HSSS18,JSZ17} and for Stokes and
Navier-Stokes equations in~\cite{Navier}. Mixed spatial and stochastic
regularity {of the solutions to the Poisson equation on random domains},
required for multilevel cubature methods, has been recently proved
in~\cite{HS21}.

In order to introduce the domain mapping approach,
let $(\Omega,\mathcal{F},\mathbb{P})$ be a probability space,
where $\Omega$ is the set of possible outcomes, the $\sigma$-algebra
$\mathcal{F}\subset 2^{\Omega}$ is the set of all events, and $\mathbb{P}$
is a probability measure. We are interested in uncertainty quantification
for the Poisson problem
\begin{align}
\begin{cases}
-\Delta u(\bsx,\omega)=f(\bsx),&\bsx\in D(\omega),\\
u(\bsx,\omega)=0,&\bsx\in \partial D(\omega),
\end{cases}\quad\text{for $\mathbb P$-almost every}~\omega\in\Omega,\label{eq:model}
\end{align}
where the domain $D(\omega)\subset\bbR^d$, $d\in\{2,3\}$, is assumed to be
uncertain. In the domain mapping framework, one starts with a {\em
reference domain} $\Dref\subset\mathbb R^d$. For example, in a
manufacturing application, $\Dref$ might be the planned domain provided by
computer-aided design with no imperfections. It is assumed that a
random perturbation field
$$
\bsV(\omega)\!:\overline{\Dref}\to\overline{D(\omega)}
$$
is prescribed. The Poisson problem~\eqref{eq:model} on the
perturbed domain is then transformed onto the fixed reference domain
$\Dref$ by a change of variable. This results in a {PDE} on the reference domain, equipped with a random coefficient
and a random source term, which are correlated by the random
deformation, and complicated by the occurrence of the Jacobian of the
transformation. Note that the dependence on the random perturbation field
is nonlinear. {Nonetheless, in the domain mapping approach, the
reference domain $\Dref$ needs only to be Lipschitz, while the
random perturbations $\bsV(\omega)$ have to be {$\calC^2$}-smooth.
This is in contrast to the perturbation approach mentioned above,
where both the reference domain and the random perturbations
need to be {$\calC^2$}-smooth.}

The approach adopted in this paper is first to truncate the initially
infinite number of scalar random variables to a finite (but possibly
large) number; then to approximate the transformed and truncated
{PDE} {on the reference domain} by a
finite element method; and finally to approximate the expected
value of the solution (which is an integral, possibly high dimensional)
by a carefully designed quasi-Monte Carlo ({QMC}) method{---i.e., an equal weight cubature rule}.

In \cite{CNT16,nummat2016}, the perturbation field was represented by an
affine, vector-valued Karhunen--Lo\`{e}ve expansion, under the
assumption of uniformly independent and identically distributed random
variables.  In this work, we instead expand the perturbation field
$\bsV(\omega)$ using \emph{periodic random variables}, following the
recent work~\cite{periodicpaper}. The advantage of the periodic setting is
that it permits higher order convergence of the {QMC} approximation to the
expectation of the solution to~\eqref{eq:model}. The periodic
representation is \emph{equivalent in law} to an affine
representation of the input random field, where the random variables
entering the series are i.i.d.~with the Chebyshev probability density
$\rho(z)=\frac{1}{\pi\sqrt{1-z^2}}$, $z\in(-1,1)$. The density $\rho$
is associated with Chebyshev polynomials of the first kind, which is
a popular family of basis functions in the method of generalized
polynomial chaos (gPC). The choice of density is a modeling assumption,
and it might be argued that in many applications choosing the Chebyshev
density over the uniform density is a matter of taste rather than
conviction, see~\cite[Sections~1--2]{kks24} for more discussion.

The paper is organized as follows. The problem is formulated in
Section~\ref{sec:2}, after which Section~\ref{sec:3} establishes the regularity of the
solution with respect to the stochastic variables; a more difficult task
than in most such applications because of the nonlinear nature of the
transformation, but a prerequisite for the design of good  {QMC} rules.
Section~\ref{sec:4} considers error analysis, estimating separately the errors from
dimension truncation, finite element approximation and QMC cubature.  In
the latter case, the error analysis leads to the design of appropriate
weight parameters  for the function space, ones that allow rigorous error
bounds for the QMC contribution to the error. Numerical experiments that test the theoretical analysis are presented in Section~\ref{sec:5},
and some conclusions are drawn in Section~\ref{sec:6}. The appendix~\ref{sec:tech} contains several technical results required for the parametric regularity analysis.%

\section{Problem formulation}\label{sec:2}

\subsection{Notations}

Let $U:={[0,1]^\mathbb{N}}$ denote a set of parameters.
Let us fix a bounded domain $\Dref\subset\bbR^d$ with Lipschitz boundary
and $d\in\{2,3\}$ as the {reference domain}.

We use multi-index notation throughout this paper. The set of all finitely
supported multi-indices is denoted by
$\mathscr{F}:=\{\bsnu\in\mathbb{N}_0^{\mathbb{N}}:|{\rm
supp}(\bsnu)|<\infty\}$, where ${\rm
supp}(\bsnu):=\{j\in\mathbb{N}:\nu_j\neq 0\}$ is the \emph{support} of a
multi-index $\bsnu:=(\nu_1,\nu_2,\ldots)$. For any sequence
$\bsb=(b_j)_{j=1}^\infty$ of real numbers and any $\bsnu\in\mathscr{F}$,
we define
$$
\bsb^{\bsnu}:=\prod_{j\in{\rm supp}(\bsnu)}b_j^{\nu_j},
$$
where the product over the empty set is defined as $1$ by convention. {We also define the multi-index notation
$$
\partial_{\bsy}^{\bsnu}:=\prod_{j\in{\rm supp}(\bsnu)}\frac{\partial^{\nu_j}}{\partial y_j^{\nu_j}}
$$
for higher order partial derivatives with respect to variable $\bsy$.}

For any matrix $M$, let $\sigma(M)$ denote the set of all singular values
of $M$ and {let} $\|M\|_2:={\max\sigma(M)}$ denote the matrix spectral norm. For
a function $v$ on $\Dref$ we define
\[
  \|v\|_{L^\infty(\Dref)}
  \,:=\,
  \begin{cases}
  \esup_{\bsx\in \Dref} |v(\bsx)| & \mbox{if } v: \Dref\to\bbR, \vspace{0.1cm} \\
  \esup_{\bsx\in \Dref} \|v(\bsx)\|_2 & \mbox{if } v: \Dref\to\bbR^d, \vspace{0.1cm} \\
  \esup_{\bsx\in \Dref} \|v(\bsx)\|_2 & \mbox{if } v: \Dref\to\bbR^{d\times d},
  \end{cases}
\]
where we apply the vector 2-norm (Euclidean norm) or matrix 2-norm
(spectral norm) depending on whether $v(\bsx)$ is a vector or a matrix.
Similarly, we define
\[
  \|v\|_{W^{1,\infty}(\Dref)}
  :=
  \begin{cases}
  \max\!\Big(\displaystyle\esup_{\bsx\in \Dref} |v(\bsx)|,\; \esup_{\bsx\in \Dref} \|\nabla v(\bsx)\|_2\Big)
  & \mbox{if } v\!: \Dref\to\bbR, \vspace{0.1cm} \\
  \max\!\Big(\displaystyle\esup_{\bsx\in \Dref} \|v(\bsx)\|_2,\; \esup_{\bsx\in \Dref} \|v'(\bsx)\|_2\Big)
  & \mbox{if } v\!: \Dref\to\bbR^d,
  \end{cases}
\]
where $\nabla v(\bsx)$ is the gradient vector if $v(\bsx)$ is a scalar,
and $v'(\bsx)$ is the Jacobian matrix if $v(\bsx)$ is a vector. We will
also need the standard Sobolev norm
\[
  \|v\|_{H_0^1(\Dref)} \,:=\,
  \|\nabla v\|_{L^2(\Dref)} \,=\, \left(\int_{\Dref} \|\nabla v(\bsx)\|_2^2\,\rd\bsx\right)^{1/2}
\]
for a scalar function $v:\Dref\to\bbR$ in $H_0^1(\Dref)$. {We also define the norm
$$
\|v\|_{\calC^k(\overline{D})}:=\max_{|\bsnu|\leq k}\sup_{\bsx\in \overline{D}}|\partial_\bsx^{\bsnu}v(\bsx)|\quad \text{for $v\in \calC^k(\overline{D})$},
$$
where $D\subset \mathbb R^d$ is a nonempty domain.}

In this paper we will make use of \emph{Stirling numbers of the second
kind} defined by
\[
  S(n,m) := \frac{1}{m!} \sum_{j=0}^m (-1)^{m-j} \binom{m}{j} j^n
\]
for integers $n\geq m\geq 0$, except for $S(0,0):=1$.

\subsection{Parameterization of domain uncertainty}

Let $\bsV\!:\overline{\Dref}\times U\to \bbR^d$ be a vector field such
that
\begin{align}
  \bsV(\bsx,\bsy)\,:=\, \bsx + \frac{1}{\sqrt 6}\sum_{i=1}^\infty \sin(2\pi y_i)\,\bspsi_i(\bsx),
  \quad \bsx\in {\Dref},~\bsy\in U,\label{eq:V}
\end{align}
with {\em stochastic fluctuations} $\bspsi_i : \Dref\to \bbR^d$. Denoting
the Jacobian matrix of $\bspsi_i$ by $\bspsi_i'$, the Jacobian matrix $J:
\Dref \to \bbR^{d\times d}$ of vector field $\bsV$ is
\begin{align*}
  J(\bsx,\bsy) \,:=\, I+\frac{1}{\sqrt{6}}\sum_{i=1}^\infty\sin(2\pi y_i)\,\bspsi_i'(\bsx),\quad \bsx\in \Dref,~\bsy\in U.
\end{align*}
We assume that the family of {\em admissible domains} $\{D(\bsy)\}_{\bsy\in U}$ is parameterized by
$$
D(\bsy):=\bsV(\Dref,\bsy),\quad \bsy\in U,
$$
and define the \emph{hold-all domain} by setting
$$
\mathcal{D}:=\bigcup_{\bsy\in U}D(\bsy).
$$

{The convention of explicitly writing down the factor $1/\sqrt 6$ in~\eqref{eq:V} is not a crucial part of the analysis, but it ensures that the mean and covariance of the vector field $\bsV$ match those of an affine and uniform parameterization for representing domain uncertainty using the same sequence of stochastic fluctuations $(\bspsi_i)_{i=1}^\infty$, with the uniform random variables supported on $[-1/2,1/2]$. See also the discussion in~\cite{periodicpaper}.

{In order to ensure that the aforementioned domain parameterization is well-posed and to enable our subsequent analysis of dimension truncation and finite element errors, we state the following assumptions for later use:}%
\begin{enumerate}[label=(A\arabic*),align=right,leftmargin=1.7\parindent]%
\item For each $\bsy\in U$,
    $\bsV(\cdot,\bsy)\!:\overline{\Dref}\to\bbR^d$ is an invertible,
    twice continuously differentiable vector field.\label{assumption1}
\item For some $C>0$, {there} holds
    $$\|\bsV(\cdot,\bsy)\|_{\calC^2(\overline{\Dref})}\leq C\quad\text{and}\quad
    \|\bsV^{-1}(\cdot,\bsy)\|_{\calC^2(\overline{D(\bsy)})}\leq C$$ for all $\bsy\in U$.\label{assumption2}
\item There exist constants $0 < \sigma_{\min} \leq 1 \leq
    \sigma_{\max}<\infty$ such that
$$
  \sigma_{\min}\leq \min \sigma(J(\bsx,\bsy))
  \leq \max\sigma(J(\bsx,\bsy))
  \leq \sigma_{\max}~\text{for all}~\bsx\in \Dref,~\bsy\in U,
$$
where $\sigma(J(\bsx,\bsy))$ denotes the set of all singular
values of the Jacobian matrix.\label{assumption3}

\item {There} holds $\|\bspsi_i\|_{W^{1,\infty}(\Dref)}<\infty$
    for all $i\in\mathbb{N}$, and $$\sum_{i=1}^\infty
    \|\bspsi_i\|_{W^{1,\infty}(\Dref)}<\infty.$$\label{assumption4}
\item For some $p\in(0,1)$, {there} holds\label{assumption5}
$$
 \sum_{i=1}^\infty \|\bspsi_i\|_{W^{1,\infty}(D_{\rm ref})}^p<\infty.
$$
\item $\|\bspsi_1\|_{W^{1,\infty}(\Dref)}\geq
    \|\bspsi_2\|_{W^{1,\infty}(\Dref)}\geq\cdots$.\label{assumption6}
\item The reference domain $D_{\rm ref}\subset\mathbb R^d$
is a convex, bounded polyhedron.\label{assumption7}
\end{enumerate}

For later convenience
    we define the sequence $\bsb=(b_i)_{i=1}^\infty$ and constant
    $\xi_{\bsb} \geq 0$ by setting
\begin{align} \label{eq:b-def}
 b_i:=\frac{1}{\sqrt{6}}\|\bspsi_i\|_{W^{1,\infty}(\Dref)}\quad\text{and}\quad \xi_{\bsb} :=\sum_{i=1}^\infty b_i<\infty.
\end{align}
It follows that we can take $\sigma_{\max} := 1 + \xi_{\bsb}$.

{\emph{Remark.} Under assumption~\ref{assumption3}, there holds
\begin{align}
\det J(\bsx,\bsy)>0\quad\text{for all}~\bsx\in\overline{\Dref},~\bsy\in U.\label{eq:detJpositive}
\end{align}
This follows from the continuity of the determinant and $\det J(\bsx,\mathbf 0)=1$.}

\subsection{The variational formulation on the reference domain}

The variational formulation of the model problem~\eqref{eq:model} can be stated as follows: for $\bsy\in U$, find $u(\cdot,\bsy)\in H_0^1(D(\bsy))$ such that
\begin{align}
 \int_{D(\bsy)}\nabla u(\bsx,\bsy)\cdot\nabla v(\bsx)\,{\rm d}\bsx
 \,=\, \int_{D(\bsy)}f(\bsx)\,v(\bsx)\,{\rm d}\bsx\quad\text{for all}~v\in H_0^1(D(\bsy)),\label{eq:weak1}
\end{align}
where $f\in \calC^\infty(\mathcal{D})$ is assumed to be an analytic
function.

We can transport the variational formulation~\eqref{eq:weak1} to the
reference domain by a change of variable. Let us first define the
matrix-valued function
\begin{align}
 B(\bsx,\bsy) \,:=\, J(\bsx,\bsy)^\top J(\bsx,\bsy),\quad \bsx\in D_{\rm ref},~\bsy\in U,\label{eq:B}
\end{align}
the transport coefficient
\begin{align}
  A(\bsx,\bsy) \,:=\, B^{-1}(\bsx,\bsy)\,\det J(\bsx,\bsy),\quad \bsx\in \Dref,~\bsy\in U,\label{eq:A}
\end{align}
and the transported source term
\begin{align*}%
  f_{\rm ref}(\bsx,\bsy) := \widehat{f}(\bsx,\bsy)\,\det J(\bsx,\bsy),\quad
  \widehat{f}(\bsx,\bsy) := f\big(\bsV(\bsx,\bsy)\big),
  \quad\bsx\in \Dref,~\bsy\in U.
\end{align*}
Then we can recast the problem~\eqref{eq:weak1} on the reference domain as
follows: for  $\bsy\in U$, find $\widehat{u}(\cdot,\bsy)\in H_0^1(\Dref)$
such that
\begin{align}
 \int_{\Dref} \big(A(\bsx,\bsy)\nabla \widehat{u}(\bsx,\bsy)\big) \cdot \nabla \widehat{v}(\bsx)\,{\rm d}\bsx
 \,=\, \int_{\Dref}f_{\rm ref}(\bsx,\bsy)\,\widehat{v}(\bsx)\,{\rm d}\bsx
\label{eq:weak2}
\end{align}
for all $\widehat{v}\in H_0^1(\Dref)$.

The solutions to problems~\eqref{eq:weak1} and~\eqref{eq:weak2} are
connected to one another by
$$
  u(\cdot,\bsy) \,=\, \widehat{u}\big(\bsV^{-1}(\cdot,\bsy),\bsy\big)
  \quad\iff \quad \widehat{u}(\cdot,\bsy) \,=\, u\big(\bsV(\cdot,\bsy),\bsy\big),\quad \bsy\in U.
$$
In the sequel, we focus on analyzing the problem~\eqref{eq:weak2}.

\section{Parametric regularity of the solution}\label{sec:3}

In order to develop higher order rank-1 lattice cubature rules for the
purpose of integrating the solution $\bsy\mapsto \widehat{u}(\cdot,\bsy)$
of~\eqref{eq:weak2} with respect to the parametric variable, we need to
derive bounds on the partial derivatives
$\partial_\bsy^{\bsnu}\widehat{u}(\cdot,\bsy)$ in the Sobolev norm
$\|\cdot\|_{H_0^1(\Dref)}$.

The analysis presented in this section is based on~\cite{nummat2016} and proceeds as follows:
\begin{itemize}
\item We derive bounds for $\partial_\bsy^\bsnu B^{-1}(\bsx,\bsy)$ and
    $\partial_\bsy^\bsnu\det J(\bsx,\bsy)$ for all $\bsnu\in
    \mathscr{F}$ in
    Lemmata~\ref{lemma:JTJinvbound} {and} \ref{lemma:detJbound}. These
    results are then used to derive  a bound on the partial
    derivatives of the coefficient $\partial_\bsy^\bsnu A(\bsx,\bsy)$
    in Lemma~\ref{lemma:Ader}.
\item Since the right-hand side of~\eqref{eq:weak2} depends on the
    parametric variable $\bsy\in U$, we develop the derivative bounds
    on $\partial_\bsy^{\bsnu}f_{\rm ref}(\bsx,\bsy)$ in
    Lemma~\ref{lemma:fref}, with the aid of the derivative bound
    developed for $\partial_\bsy^\bsnu\widehat{f}(\bsx,\bsy)$ in
    Lemma~\ref{lemma:fhat}.
\item An {\em a priori} bound is developed for the the solution
    $\widehat{u}(\bsx,\bsy)$ of~\eqref{eq:weak2} in
    Lemma~\ref{lemma:apriori}.
\item Finally, the main result of this section is stated in
    Theorem~\ref{thm:uhatregularity} which contains the partial
    derivative bound for $\partial_\bsy^\bsnu \widehat{u}(\bsx,\bsy)$,
    $\bsnu\in \mathscr{F}$, making use of the aforementioned results.
\end{itemize}

\subsection{Parametric regularity of the transport coefficient}

\begin{lemma}\label{lemma:JTJinvbound}
Let $\bsy\in U$ and $\bsnu\in \mathscr{F}$. {Under assumptions~{\rm \ref{assumption1}}--\,{\rm \ref{assumption4}}}, the matrix-valued
function $B$ defined by~\eqref{eq:B} satisfies
\begin{align*}
 &\big\|\partial_{\bsy}^{\bsnu} B^{-1}(\cdot,\bsy) \big\|_{L^\infty(\Dref)}\\
 &\le\, \frac{1}{\sigma_{\min}^2}\,
 \bigg(\frac{4\pi(1+\xi_{\bsb})}{\sigma_{\min}^2}\bigg)^{|\bsnu|}
 \sum_{\bsm\leq \bsnu} |\bsm|!\,a_{|\bsm|}\, \bsm!\,\bsb^{\bsm}\prod_{i\geq 1}S(\nu_i,m_i),
\end{align*}%
with the sequence
\begin{align} \label{eq:ak}
  a_k \,:=\, \frac{(1+\sqrt 3)^{k+1}-(1-\sqrt 3)^{k+1}}{2^{k+1}\sqrt 3}, \qquad k\ge 0.
\end{align}
\end{lemma}

\proof We will first derive a regularity bound for $\|\partial_{\bsy}^{\bsnu}B(\cdot,\bsy)\|_{L^\infty(D_{\rm ref})}$, $\bsnu\in\mathscr F$. The bound for $\|\partial_{\bsy}^{\bsnu}B(\cdot,\bsy)^{-1}\|_{L^\infty(D_{\rm ref})}$ follows by applying implicit differentiation to the identity $B^{-1}B=I$.

From the definition of $J$ it is easy to see that for any
multi-index $\bsm\in\mathscr{F}$ and any $\bsx\in\Dref$ and $\bsy\in U$ we
have
\begin{align} \label{eq:derJ}
  \partial^\bsm_\bsy J(\bsx,\bsy) \,=\,
  \begin{cases}
  I + \frac{1}{\sqrt{6}}\sum_{i\ge1}\sin(2\pi y_i)\,\bspsi_i'(\bsx) & \mbox{if } \bsm = \bszero, \\
  \frac{1}{\sqrt{6}}\,(2\pi)^k\,\sin(2\pi y_j + k\, \frac{\pi}{2})\,\bspsi_j'(\bsx) & \mbox{if } \bsm = k\,\bse_j,\; k\ge 1,\\
  0 & \mbox{otherwise.}
  \end{cases}
\end{align}
Taking the matrix $2$-norm, we obtain
\begin{align*}%
  \|\partial^\bsm_\bsy J(\bsx,\bsy)\|_2 \,\le\,
  \begin{cases}
  1 + \frac{1}{\sqrt{6}} \sum_{i\ge 1} \|\bspsi_i'(\bsx)\|_2& \mbox{if } \bsm = \bszero, \\
  \frac{1}{\sqrt{6}}\,(2\pi)^k\,\|\bspsi_j'(\bsx)\|_2 & \mbox{if } \bsm = k\,\bse_j,\; k\ge 1,\\
  0 & \mbox{otherwise,}
  \end{cases}
\end{align*}
and hence with \eqref{eq:b-def} we obtain
\begin{align} \label{eq:derJinf}
  \|\partial^\bsm_\bsy J(\cdot,\bsy)\|_{L^\infty(\Dref)} \,\le\,
  \begin{cases}
  1 + \xi_{\bsb}  & \mbox{if } \bsm = \bszero, \\
  (2\pi)^k\,b_j & \mbox{if } \bsm = k\,\bse_j,\; k\ge 1,\\
  0 & \mbox{otherwise.}
  \end{cases}
\end{align}

Now the Leibniz product rule yields for $\bsnu\in\mathscr{F}$,
\[
  \partial^\bsnu_\bsy B(\cdot,\bsy)
  \,=\, \sum_{\bsm\le\bsnu} \binom{\bsnu}{\bsm}
  \big(\partial^\bsm_\bsy J(\bsx,\bsy)^\top\big) \big(\partial^{\bsnu-\bsm}_\bsy J(\bsx,\bsy)\big),
\]
and after taking the matrix 2-norm on both sides and then the
$L^\infty$-norm over $\bsx$ we obtain
\[
  \big\|\partial^\bsnu_\bsy B(\cdot,\bsy)\big\|_{L^\infty(\Dref)}
  \le \sum_{\bsm\le\bsnu} \binom{\bsnu}{\bsm}
  \big\|\partial^\bsm_\bsy J(\cdot,\bsy)\big\|_{L^\infty(\Dref)}
  \big\|\partial^{\bsnu-\bsm}_\bsy J(\cdot,\bsy)\big\|_{L^\infty(\Dref)}.
\]
The bounds in \eqref{eq:derJinf} indicate that only the derivatives with
$|{\rm supp}(\bsnu)|\le 2$ will survive. Indeed, for $\bsnu = k\,\bse_j$
with $k\ge 1$, only $\bsm = \ell\,\bse_j$ for $\ell=0,\ldots, k$ remain,
while for $\bsnu = k\,\bse_j + k'\bse_{j'}$ with $k,k'\ge 1$ and $j\ne
j'$, only $\bsm = k\,\bse_j$ and $\bsm = k'\bse_{j'}$ remain. On the other
hand, for $\bsnu$ with $|{\rm supp}(\bsnu)|\geq 3$, in each term at least
one of the factors must vanish. We conclude that
\begin{align} %
  &\big\|\partial^\bsnu_\bsy B(\cdot,\bsy)\big\|_{L^\infty(\Dref)} \nonumber\\
  &\,\le\,
  \begin{cases}
  (1 + \xi_{\bsb} )^2 & \mbox{if } \bsnu = \bszero, \\
  2\,(1 + \xi_{\bsb} )\, (2\pi)^k\,b_j + (2^k-2)\, (2\pi)^k\,b_j^2  & \mbox{if } \bsnu = k\,\bse_j,\; k\ge 1,\\
  2\,(2\pi)^{k+k'}\,b_j\,b_{j'} & \mbox{if } \bsnu=k\,\bse_j\!+\!k'\bse_{j'},\,\! k,k'\!\ge\!1,\, j\!\ne\!j'\\
  0 & \mbox{if $|{\rm supp}(\bsnu)|\geq 3$,}
  \end{cases}\notag
  \\
  &\,\le\,
  \begin{cases}
  (1+\xi_{\bsb} )^2&\text{if}~\bsnu=\bszero,\\
  (4\pi)^{|\bsnu|}(1+\xi_{\bsb} )\prod_{j\in{\rm supp}(\bsnu)}b_j&\text{if}~|{\rm supp}(\bsnu)|\in\{1,2\},\\
  0 &\text{if}~|{\rm supp}(\bsnu)|\geq 3,
  \end{cases}\label{eq:derBinf}
\end{align}
where we used $b_j\leq \xi_{\bsb} $ for all $j\in\mathbb N$.

Our assumption~\ref{assumption3} immediately yields for all $\bsy\in U$ that
\[
  \|B^{-1}(\cdot,\bsy)\|_{L^\infty(\Dref)}
  \,=\,\esup_{\bsx\in \Dref} \Big(\max\,\sigma\big(J^{-1}(\bsx,\bsy)\big) \Big)^2
  \,\leq\, \frac{1}{\sigma_{\min}^2},
\]
which proves the lemma for the case $\bsnu = \bszero$. Now we consider
$\bsnu\in\mathscr F\setminus\{\bszero\}$. To obtain the derivative bounds
on $B^{-1}$ we start with $B^{-1}B = I$ and apply the Leibniz product rule
to obtain
\[
  \partial^\bsnu_\bsy \big(B^{-1}(\bsx,\bsy) B(\bsx,\bsy)\big)
  \,=\, \sum_{\bsm\le\bsnu} \binom{\bsnu}{\bsm}
  \big(\partial^{\bsnu-\bsm}_\bsy B^{-1}(\bsx,\bsy)\big) \big(\partial^\bsm_\bsy B(\bsx,\bsy)\big) \,=\, {0}.
\]
Separating out the $\bsm = \bszero$ term, post-multiplying both sides by
$B^{-1}(\bsx,\bsy)$, taking the matrix 2-norm followed by the
$L^\infty$-norm over $\bsx$, and applying \eqref{eq:derBinf}, we obtain
\begin{align*}
&\big\|\partial^\bsnu_\bsy B^{-1}(\cdot,\bsy)\big\|_{L^\infty(\Dref)}\\
&\leq\! \|B^{-1}(\cdot,\bsy)\|_{L^\infty\!(D_{\rm ref})}\!\!\!\!
  \sum_{\bszero\ne\bsm\le\bsnu} \!\!\!\binom{\!\bsnu\!}{\!\bsm\!}\!
  \big\|\partial^{\bsnu-\bsm}_\bsy B^{-1}(\cdot,\bsy)\big\|_{L^\infty\!(\Dref)}
  \big\|\partial^\bsm_\bsy B(\cdot,\bsy)\big\|_{L^\infty\!(\Dref)} \\
&\leq \frac{1+\xi_{\bsb} }{\sigma_{\min}^2}
 \sum_{\substack{\bszero\neq \bsm\leq \bsnu\\ |{\rm supp}(\bsm)|\leq 2}}
 (4\pi)^{|\bsm|}\binom{\bsnu}{\bsm}\, \big\|\partial^{\bsnu-\bsm}_\bsy B^{-1}(\cdot,\bsy)\big\|_{L^\infty(\Dref)}
 \prod_{j\in{\rm supp}(\bsm)} b_j.
\end{align*}
The assertion follows by applying Lemmata~\ref{lemma:delannoyrecursion}
and~\ref{lemma:delannoybound} with
$\Upsilon_{\bsnu}=\big\|\partial^{\bsnu}_\bsy
B^{-1}(\cdot,\bsy)\big\|_{L^\infty(\Dref)}$, $c_0 = 1/\sigma^2_{\min}$, $c
= 4\pi(1+\xi_{\bsb})/\sigma_{\min}^2$, $q = 2$, and $\beta_j=b_j$.
\quad\qed

\begin{lemma}\label{lemma:detJbound}
{Under assumptions~{\rm \ref{assumption1}}--\,{\rm \ref{assumption4}}}, {there} holds for all $\bsy\in U$ and $\bsnu\in\mathscr{F}$ that
$$
  \|\partial_{\bsy}^{\bsnu}\det J(\cdot,\bsy)\|_{L^\infty(\Dref)}
  \,\le\, d!\,(2+\xi_{\bsb} )^d\,(2\pi)^{|\bsnu|}\,
  \sum_{\bsw\leq\bsnu}|\bsw|!\,\bsb^{\bsw}\prod_{i\geq 1}S(\nu_i,w_i).
$$
\end{lemma}
\proof Let $\bsy\in U$. The proof is carried out by induction over all the
minors (subdeterminants) of the matrix
${J}(\bsx,\bsy)=:[J_{\ell,\ell'}(\bsx,\bsy)]_{\ell,\ell'=1}^d$, $\bsx\in
\Dref$. For any $q\le d$, let $J^{\bsell,\bsell'}(\bsx,\bsy) :=
[J_{\ell_t,\ell'_{t'}}(\bsx,\bsy)]_{t,t'=1}^q$ denote the $q\times q$
submatrix specified by the indices $\bsell := \{\ell_1,\ldots,\ell_q\}$
and $\bsell':=\{\ell'_1,\ldots,\ell'_q\}$ with $1\leq
\ell_1<\cdots<\ell_q\leq d$ and $1\leq \ell'_1<\cdots<\ell'_q\leq d$. We
will prove by induction on $q$ that
\begin{align} \label{eq:hyp}
 \|\partial_{\bsy}^{\bsnu} \det J^{\bsell,\bsell'}(\cdot,\bsy)\|_{L^\infty(\Dref)}
 \,\le\, q!\,(2+\xi_{\bsb} )^q\,(2\pi)^{|\bsnu|}\,\sum_{\bsw\leq\bsnu}|\bsw|!\,\bsb^{\bsw}\prod_{i\geq 1}S(\nu_i,w_i).
\end{align}
Similarly to \eqref{eq:derJ}, we have for the derivatives of matrix
elements
\begin{align}\label{eq:derqeq1}
  \partial^\bsm_\bsy J_{\ell_t,\ell'_{t'}}(\bsx,\bsy) \,=\,
  \begin{cases}
  1 + \frac{1}{\sqrt{6}}\sum_{i\ge 1} \sin(2\pi y_i)\, [\bspsi_i'(\bsx)]_{\ell_t,\ell'_{t'}} & \mbox{if } \bsm = \bszero, \\
  \frac{1}{\sqrt{6}}(2\pi)^k\sin(2\pi y_j + k \frac{\pi}{2})[\bspsi_j'(\bsx)]_{\ell_t,\ell'_{t'}}
  & \mbox{if } \bsm = k\,\bse_j,\, k\ge 1,\\
  0 & \mbox{otherwise,}
  \end{cases}
\end{align}
for $t,t'\in\{1,\ldots,q\}$. Using $\big|
[\bspsi_i'(\bsx)]_{\ell_t,\ell'_{t'}}\big| \le \|\bspsi_i'(\bsx)\|_2$ and
\eqref{eq:b-def}, we obtain for the $(1\times 1)$-minors
\begin{align*}%
  \|\partial_{\bsy}^{\bsm}\det J^{\ell_t,\ell'_{t'}}(\cdot,\bsy)\|_{L^\infty(\Dref)}
  &= \|\partial_{\bsy}^{\bsm} J_{\ell_t,\ell_{t'}}(\cdot,\bsy)\|_{L^\infty(\Dref)}\\
  &\le
  \begin{cases}
  1+\xi_{\bsb}  &\text{if}~\bsm=\bszero,\\
  (2\pi)^k\, b_j  &\text{if}~\bsm=k\,\bse_j,\; k\ge 1,\\
  0 & \text{otherwise}
\end{cases}
\end{align*}

Suppose \eqref{eq:hyp} holds for all submatrices of size up to
$(q-1)\times (q-1)$, and now we consider the case for a $q\times q$
submatrix. For arbitrary $t'\le q$, the Laplace cofactor expansion yields
$$
  \det J^{\bsell,\bsell'}(\bsx,\bsy)
  \,=\, \sum_{t=1}^q (-1)^{t+t'}\, J_{\ell_t,\ell'_{t'}}(\bsx,\bsy)\,
  \det J^{\bsell\setminus\{\ell_t\},\bsell'\setminus\{\ell'_{t'}\}}(\bsx,\bsy).
$$
The Leibniz product rule then gives
\begin{align*}
  &\partial^\bsnu_\bsy \det J^{\bsell,\bsell'}(\bsx,\bsy)\\
  &= \sum_{t=1}^q (-1)^{t+t'}\!\! \sum_{\bsm\le\bsnu} \!\! \binom{\bsnu}{\bsm}
  \Big(\partial^\bsm_\bsy J_{\ell_t,\ell'_{t'}}(\bsx,\bsy)\Big)
  \Big(\partial^{\bsnu-\bsm}_\bsy \det J^{\bsell\setminus\{\ell_t\},\bsell'\setminus\{\ell'_{t'}\}}(\bsx,\bsy)\Big).
\end{align*}
Together with \eqref{eq:derqeq1} and the induction hypothesis
\eqref{eq:hyp} we obtain
\begin{align}
  &\|\partial_{\bsy}^{\bsnu}\det J^{\bsell,\bsell'}(\cdot,\bsy)\|_{L^\infty(\Dref)} \nonumber\\
  &\,\le\, \sum_{t=1}^q \bigg(\|J_{\ell_t,\ell'_{t'}}(\cdot,\bsy)\|_{L^\infty(\Dref)}
  \|\partial_{\bsy}^{\bsnu}\det J^{\bsell\setminus\{\ell_t\},\bsell'\setminus\{\ell'_{t'}\}}(\cdot,\bsy)\|_{L^\infty(\Dref)} \nonumber\\
  &\qquad\quad
  +\sum_{j\ge 1} \sum_{k=1}^{\nu_j}\binom{\nu_j}{k}\|\partial^{k\bse_j}J_{\ell_t,\ell'_{t'}}(\cdot,\bsy)\|_{L^\infty(\Dref)}\nonumber\\
  &\qquad\qquad\qquad\qquad\qquad\quad \times
  \|\partial^{\bsnu-k\bse_j}\det J^{\bsell\setminus\{\ell_t\},\bsell'\setminus\{\ell'_{t'}\}}(\cdot,\bsy)\|_{L^\infty(\Dref)}\bigg)
  \nonumber\\
  &\,\le\, \sum_{t=1}^q \bigg( (1+\xi_{\bsb} )\,
  (q-1)!\,(2+\xi_{\bsb} )^{q-1}\,(2\pi)^{|\bsnu|}\,\sum_{\bsw\leq\bsnu}|\bsw|!\,\bsb^{\bsw}\prod_{i\geq 1}S(\nu_i,w_i) \nonumber\\
  &\qquad\quad
  +\sum_{j\ge 1} \sum_{k=1}^{\nu_j}\binom{\nu_j}{k} (2\pi)^k\, b_j\,
  (q-1)!\,(2+\xi_{\bsb} )^{q-1}(2\pi)^{|\bsnu|-k}\, \nonumber\\
  &\qquad\qquad\qquad\qquad\qquad\quad \times \sum_{\bsw\leq\bsnu-k\bse_j}\bigg(|\bsw|!\,\bsb^{\bsw}\,
  S(\nu_j-k,w_j)\prod_{\satop{i\ge 1}{i\ne j}}S(\nu_i,w_i) \bigg)\bigg). \label{eq:derJfinal}
\end{align}
To simplify the last term, we obtain in complete analogy with the
derivation presented in~\cite[Lemma~2.2]{periodicpaper} the identity
\begin{align*}
& \sum_{j\geq 1}\sum_{k=1}^{\nu_j}\binom{\nu_j}{k} b_j\sum_{\bsw\leq \bsnu-k\bse_j}|\bsw|!\,\bsb^{\bsw}\,
 S(\nu_j-k,w_j) \prod_{\satop{i\ge 1}{i\ne j}} S(\nu_i,w_i)\\
 &=\, \sum_{\bsw\leq\bsnu}|\bsw|!\,\bsb^{\bsw}\prod_{i\ge 1}S(\nu_i,w_i).
\end{align*}
Plugging this expression into~\eqref{eq:derJfinal} and collecting terms
yields \eqref{eq:hyp}, which in turn completes the proof of the
lemma.\quad\qed

\begin{lemma}\label{lemma:Ader} {Under assumptions~{\rm \ref{assumption1}}--\,{\rm \ref{assumption4}}}, {there} holds for
all $\bsy\in U$ and $\bsnu\in\mathscr{F}$ that
\begin{align*}
 & \|\partial_\bsy^{\bsnu}A(\cdot,\bsy)\|_{L^\infty(\Dref)}\\
  &\leq \frac{d!\,(2+\xi_{\bsb} )^d}{\sigma_{\min}^2}
  \bigg(\frac{4\pi(1+\xi_{\bsb})}{\sigma_{\min}^2}\bigg)^{|\bsnu|}
  \sum_{\bsm\leq \bsnu}(|\bsm|+1)!\, a_{|\bsm|}\,\bsm!\,
  \bsb^{\bsm}\prod_{i\geq 1}S(\nu_i,m_i),
\end{align*}
where the sequence $(a_k)_{k=0}^\infty$ is defined in \eqref{eq:ak}.
\end{lemma}

\proof Let $\bsnu\in\mathscr{F}$. By the Leibniz product rule, we have
from \eqref{eq:A} that
\[
 \partial_\bsy^{\bsnu}A(\bsx,\bsy)
 = \sum_{\bsm\leq\bsnu}\binom{\bsnu}{\bsm} \Big(\partial_\bsy^{\bsm} B^{-1}(\bsx,\bsy)\Big)\!
 \Big(\partial_\bsy^{\bsnu-\bsm}\det J(\bsx,\bsy)\Big),\,\bsx\in \Dref,\,\bsy\in U.
\]
Taking the matrix 2-norm followed by the $L^\infty$-norm over $\bsx$, and
then applying Lemmata~\ref{lemma:JTJinvbound} and~\ref{lemma:detJbound},
we obtain
\begin{align*}
  &\|\partial_\bsy^{\bsnu}A(\cdot,\bsy)\|_{L^\infty(\Dref)}\\
  &\le \sum_{\bsm\leq\bsnu}\binom{\bsnu}{\bsm}
  \big\|\partial_\bsy^{\bsm} B^{-1}(\cdot,\bsy) \big\|_{L^\infty(\Dref)}
  \big\|\partial_\bsy^{\bsnu-\bsm}\det J(\cdot,\bsy)\big\|_{L^\infty(\Dref)}\\
  &\leq  c_{|\bsnu|}\sum_{\bsm\leq \bsnu}\sum_{\bsw\leq \bsm}\sum_{\bsmu\leq \bsnu-\bsm}\binom{\bsnu}{\bsm}
  \bbA_\bsw\,\bbB_\bsmu\,\prod_{i\geq 1} \Big( S(m_i,w_i)\,S(\nu_i-m_i,\mu_i) \Big) \\
  &= c_{|\bsnu|} \sum_{\bsm\leq\bsnu} \bigg(\sum_{\bsw\leq\bsm} \binom{\bsm}{\bsw}\,\bbA_\bsw\, \bbB_{\bsm-\bsw}\bigg)
  \prod_{i\geq 1}S(\nu_i,m_i),
\end{align*}
where $\bbA_\bsw := |\bsw|!\,a_{|\bsw|}\,\bsw!\,\bsb^{\bsw}$ and
$\bbB_\bsmu := |\bsmu|!\, b^{\bsmu}$, and we overestimated some
multiplying factors to get
$$
 c_k \,:=\,
 \frac{d!\,(2+\xi_{\bsb} )^d}{\sigma_{\min}^2}\,
 \bigg(\frac{4\pi(1+\xi_{\bsb})}{\sigma_{\min}^2}\bigg)^{k},
 \qquad k\geq 0.
$$
The last equality is due to Lemma~\ref{lem:triple}. Using $a_{|\bsw|}\le
a_{|\bsm|}$ and $\bsw!\le \bsm!$, we have
\begin{align*}
  \sum_{\bsw\leq\bsm} \binom{\bsm}{\bsw}\,\bbA_\bsw\, \bbB_{\bsm-\bsw}
  &\le a_{|\bsm|}\,\bsm!\,\bsb^\bsm \sum_{\bsw\leq\bsm} \binom{\bsm}{\bsw}|\bsw|!\,|\bsm-\bsw|!,
\end{align*}
where the last sum over $\bsw$ can be rewritten as
\begin{align*}
 \sum_{j=0}^{|\bsm|} j!\,(|\bsm|-j)! \sum_{\substack{\bsw\leq \bsm\\ |\bsw|=j}}\binom{\bsm}{\bsw}
 = \sum_{j=0}^{|\bsm|} j!\,(|\bsm|-j)! \binom{|\bsm|}{j}
 = (|\bsm|+1)!,
\end{align*}
thus completing the proof.
\quad\qed

\subsection{Parametric regularity of the transported source term}

\begin{lemma}\label{lemma:fhat}
Let $f\in \calC^\infty(\mathcal{D})$ be analytic. Then there exist
constants $C_f>0$ and $\rho\ge 1$ such that
$\|\partial_\bsx^{\bsnu}f\|_{L^\infty(\mathcal{D})}\leq C_f\,
\bsnu!\,\rho^{|\bsnu|}$ for all $\bsnu\in\mathbb{N}_0^d$. Furthermore, {under the assumptions~{\rm \ref{assumption1}}--\,{\rm \ref{assumption4}}},  {there}
holds for all $\bsy\in U$ and $\bsnu\in\mathscr{F}$ that the derivatives
of $\widehat{f}(\bsx,\bsy)=f(\bsV(\bsx,\bsy))$ are bounded by
$$
 \big\|\partial_{\bsy}^{\bsnu}\widehat{f}(\cdot,\bsy)\big\|_{L^\infty(\Dref)}
 \le C_f(2\pi)^{|\bsnu|} \sum_{\bsm\leq \bsnu}
 \binom{|\bsm|+d-1}{d-1} |\bsm|!\rho^{|\bsm|} \bsb^{\bsm}\prod_{i\geq 1}S(\nu_i,m_i).
$$
\end{lemma}
\proof The bound $\|\partial_\bsx^{\bsnu}f\|_{L^\infty(\mathcal{D})}\leq
C_f\, \bsnu!\,\rho^{|\bsnu|}$ is a consequence of the Cauchy integral
formula for analytic functions of several variables (cf.,
e.g.,~\cite[Theorem~2.2.1]{hormander}).

Let $\bsy\in U$. Trivially, we see that $\|\widehat
f(\cdot,\bsy)\|_{L^\infty(D_{\rm ref})}\leq C_f$, so we may assume in the
following that $\bsnu\in\mathscr F\setminus\{\bszero\}$. We will make use
of the multivariate Fa\`{a} di Bruno formula~\cite[Theorem~3.1 and
Remark~3.3]{savits}
\begin{align}
 \partial_\bsy^{\bsnu}\widehat f(\bsx,\bsy)
 \,:=\, \sum_{\substack{\bslambda\in\mathbb N_0^d\\ 1\leq |\bslambda|\leq |\bsnu|}}
 (\partial_{\boldsymbol x}^{\bslambda}f)(\boldsymbol V(\bsx,\bsy))\,
 \kappa_{\bsnu,\bslambda}(\bsx,\bsy),\label{eq:fhatderiv}
\end{align}
where we define $\kappa_{\bsnu,\bslambda}(\bsx,\bsy)$ for general
$\bsnu\in\mathscr{F}$ and $\bslambda \in\bbZ^d$ in a recursive manner as follows:
$\kappa_{\bsnu,\bszero}(\bsx,\bsy) := \delta_{\bsnu,\bszero}$,
$\kappa_{\bsnu,\bslambda}(\bsx,\bsy) := 0$ if $|\bsnu|<|\bslambda|$ or
$\bslambda\not\geq \bszero$ (i.e., $\bslambda$ contains negative entries),
and otherwise
\begin{align}
 \kappa_{\bsnu+\bse_j,\bslambda}(\bsx,\bsy)
 \,:=\, \sum_{\ell=1}^d\sum_{\bszero\leq \bsm\leq \bsnu}
 \binom{\bsnu}{\bsm}\partial_\bsy^{\bsm+\bse_j}[\bsV(\bsx,\bsy)]_\ell\;
 \kappa_{\bsnu-\bsm,\bslambda-\bse_\ell}(\bsx,\bsy). \label{eq:faadibruno}
\end{align}

From the definition of {$\bsV$} {in \eqref{eq:V}} it is easy to see that
\begin{align*} %
  \partial^\bsm_\bsy [\bsV(\bsx,\bsy)]_\ell \,=\,
  \begin{cases}
  x_\ell + \frac{1}{\sqrt{6}}\sum_{i\ge1}\sin(2\pi y_i)\,[\bspsi_i(\bsx)]_\ell & \mbox{if } \bsm = \bszero, \\
  \frac{1}{\sqrt{6}}\,(2\pi)^k\,\sin(2\pi y_j + k\,\frac{\pi}{2})\,[\bspsi_j(\bsx)]_\ell & \mbox{if } \bsm = k\,\bse_j,\; k\ge 1,\\
  0 & \mbox{otherwise,}
  \end{cases}
\end{align*}
which yields
$$
 \big\|\partial_\bsy^{\bsm+\bse_j}[\bsV(\cdot,\bsy)]_\ell\big\|_{L^\infty(D_{\rm ref})}
 \leq \begin{cases}(2\pi)^{k+1}\, b_j & \mbox{if } \bsm = k\,\bse_j,\; k\ge 0,\\
0&\text{otherwise}.
\end{cases}
$$
Thus we obtain from~\eqref{eq:faadibruno} the recursion
$$
 \|\kappa_{\bsnu+\boldsymbol e_j,\bslambda}(\cdot,\bsy)\|_{L^\infty(D_{\rm ref})}
 \leq b_j\sum_{k=0}^{\nu_j}(2\pi)^{k+1}\binom{\nu_j}{k}
 \sum_{\substack{\ell=1\\ \lambda_\ell>0}}^d\|\kappa_{\bsnu-k \bse_j,\bslambda-\bse_\ell}(\cdot,\bsy)\|_{L^\infty(D_{\rm ref})}.
$$
By Lemma~\ref{lemma:faadibrunorecursion} with $c=2\pi$, we have for all
$\bsnu\in\mathscr{F}$ and $\bslambda\in \mathbb N_0^d$ that
$$
 \|\kappa_{\bsnu,\bslambda}(\cdot,\bsy)\|_{L^\infty(D_{\rm ref})}
 \leq (2\pi)^{|\bsnu|}\,\frac{|\bslambda|!}{\bslambda!}\sum_{\substack{\bsm\leq \bsnu\\ |\bsm|=|\bslambda|}}\bsb^{\bsm}
 \prod_{i\geq 1}S(\nu_i,m_i),
$$
which together with~\eqref{eq:fhatderiv} yields for
$\bsnu\in\mathscr{F}\setminus\{\bszero\}$ that
\begin{align*}
 \big\|\partial_\bsy^{\bsnu}\widehat f(\cdot,\bsy)\big\|_{L^\infty(D_{\rm ref})}
 &\leq C_f(2\pi)^{|\bsnu|}\sum_{\substack{\bslambda\in\mathbb N_0^d\\ 1\leq |\bslambda|\leq |\bsnu|}}
 |\bslambda|!\,\rho^{|\bslambda|}\sum_{\substack{\bsm\leq \bsnu\\ |\bsm|=|\bslambda|}}\bsb^{\bsm}\prod_{i\geq 1}S(\nu_i,m_i)\\
 &= C_f(2\pi)^{|\bsnu|}\sum_{\bszero\ne \bsm\leq \bsnu} |\bsm|!\rho^{|\bsm|}\, \bsb^{\bsm}
 \bigg(\prod_{i\geq 1}S(\nu_i,m_i)\bigg)\sum_{\substack{\bslambda\in\mathbb N_0^d\\ |\bslambda|=|\bsm|}}1.
\end{align*}
This yields the desired result since $\sum_{\bslambda\in\mathbb
N_0^d,\,|\bslambda|=|\bsm|}1 = \binom{|\bsm|+d-1}{d-1}$.
The $\bsm\ne\bszero$ condition for the above sum can be dropped since the
$\bsm=\bszero$ term is actually zero when $\bsnu\ne\bszero$. This allows
us to then extend the formula to cover also the case
$\bsnu=\bszero$.
\quad\qed

\begin{lemma}\label{lemma:fref}
Under the assumptions of Lemma~\ref{lemma:fhat}, {there}
holds for all $\bsy\in U$ and $\bsnu\in\mathscr{F}$ that
\begin{align*}
 &\|\partial_\bsy^\bsnu f_{\rm ref}(\cdot,\bsy)\|_{L^\infty(\Dref)}\\
 &\,\le\, C_f\,d!\,(2+\xi_{\bsb})^d\,(2\pi)^{|\bsnu|}\sum_{\bsm\leq\bsnu}
 \frac{(|\bsm|+d)!}{d!}\,\rho^{|\bsm|}\,\bsb^{\bsm}\prod_{i\geq 1}S(\nu_i,m_i).
\end{align*}
\end{lemma}

\proof The proof follows essentially the same steps as the proof of
Lemma~\ref{lemma:Ader}. By the Leibniz product rule, we have
$$
 \partial_\bsy^\bsnu f_{\rm ref}(\bsx,\bsy)
 \,=\, \sum_{\bsm\leq \bsnu}\binom{\bsnu}{\bsm}
 \Big(\partial_\bsy^{\bsm}\widehat{f}(\bsx,\bsy)\Big)
 \Big(\partial_\bsy^{\bsnu-\bsm}\det J(\bsx,\bsy)\Big),\,\bsx\in \Dref,\,\bsy\in U.
$$
Taking the $L^\infty$-norm over $\bsx$ and using
Lemmata~\ref{lemma:detJbound} and~\ref{lemma:fhat} yields
\begin{align*}
 &\|\partial_\bsy^\bsnu f_{\rm ref}(\cdot,\bsy)\|_{L^\infty(\Dref)}\\
 &\leq \sum_{\bsm\leq\bsnu}\binom{\bsnu}{\bsm}\|\partial_\bsy^{\bsm}\widehat f(\cdot,\bsy)\|_{L^\infty(D_{\rm ref}))}
 \|\partial_\bsy^{\bsnu-\bsm}\det J(\bsx,\bsy)\|_{L^\infty(D_{\rm ref})}\\
  &\le c_{|\bsnu|} \sum_{\bsm\leq \bsnu}\sum_{\bsw\leq \bsm}\sum_{\bsmu\leq \bsnu-\bsm}\binom{\bsnu}{\bsm}
  \bbA_\bsw\, \bbB_\bsmu \prod_{i\geq 1} \Big( S(m_i,w_i)\,S(\nu_i-m_i,\mu_i) \Big) \nonumber\\
  &= c_{|\bsnu|} \sum_{\bsm\leq\bsnu} \bigg(\sum_{\bsw\leq\bsm}\binom{\bsm}{\bsw}\, \bbA_\bsw\, \bbB_{\bsm-\bsw}\bigg)
  \prod_{i\geq 1}S(\nu_i,m_i),
\end{align*}%
where now $\bbA_\bsw :=
\binom{|\bsw|+d-1}{d-1}\,|\bsw|!\,\rho^{|\bsw|}\,\bsb^{\bsw}$,
$\bbB_{\bsmu} := |\bsmu|!\,\bsb^{\bsmu}$, and $c_k := C_f\,
d!\,(2+\xi_{\bsb})^d\, (2\pi)^k$. Using $\rho^{|\bsw|}\le \rho^{|\bsm|}$,
we have
\begin{align*}
  \sum_{\bsw\leq\bsm}\binom{\bsm}{\bsw}\, \bbA_\bsw\, \bbB_{\bsm-\bsw}
  &\le \rho^{|\bsm|}\,\bsb^\bsm \sum_{\bsw\leq\bsm} \binom{\bsm}{\bsw} \binom{|\bsw|+d-1}{d-1}\,|\bsw|!\,|\bsm-\bsw|!,
\end{align*}
where the last sum over $\bsw$ can be rewritten as
\begin{align*}
  &\sum_{j=0}^{|\bsm|} \frac{(j+d-1)!\,(|\bsm|-j)!}{(d-1)!}
  \sum_{\satop{\bsw\leq\bsm}{|\bsw|=j}} \binom{\bsm}{\bsw}
  = \sum_{j=0}^{|\bsm|} \frac{(j+d-1)!\,(|\bsm|-j)!}{(d-1)!} \binom{|\bsm|}{j} \\
  &= |\bsm|!\sum_{j=0}^{|\bsm|} \binom{j+d-1}{j}
  = |\bsm|!\,\binom{|\bsm|+d}{|\bsm|}
  = \frac{(|\bsm|+d)!}{d!}.
\end{align*}
This completes the proof.
\quad\qed

\subsection{Parametric regularity of the transported PDE}

We have the following {\em a priori} bound.

\begin{lemma}\label{lemma:apriori}
Under the assumptions of Lemma~\ref{lemma:fhat}, {there}
holds for all $\bsy\in U$ that
\begin{align} \label{eq:Cinit}
 \|\widehat{u}(\cdot,\bsy)\|_{H_0^1(\Dref)}
 \leq \frac{(1+\xi_\bsb)^2}{\sigma_{\min}^d}\,C_P\,|\Dref|^{1/2}\,
 C_f\,d!\,(2+\xi_\bsb)^d \,=: C_{\rm init},
\end{align}
where $C_P>0$ is the Poincar\'{e} constant associated with the domain
$\Dref$.
\end{lemma}

\proof We take $\widehat{v}(\bsx) = \widehat{u}(\bsx,\bsy)$ as test
function in \eqref{eq:weak2} to obtain
\begin{align} \label{eq:myLHS}
 &\int_{\Dref}(A(\bsx,\bsy)\nabla \widehat{u}(\bsx,\bsy))\cdot\nabla \widehat{u}(\bsx,\bsy)\,{\rm d}\bsx
 \,=\, \int_{\Dref}f_{\rm ref}(\bsx,\bsy)\,\widehat{u}(\bsx,\bsy)\,{\rm d}\bsx\\
 &\,\le\, C_P\, |\Dref|^{1/2}\,\|f_{\rm ref}(\cdot,\bsy)\|_{L^\infty(\Dref)}\|\widehat{u}(\cdot,\bsy)\|_{H^1_0(\Dref)}, \nonumber
\end{align}
where we used {the} Cauchy--Schwarz inequality and the Poincar\'{e} inequality
$\|\cdot\|_{L^2(\Dref)}\leq C_P\,\|\cdot\|_{H_0^1(\Dref)}$, with the
constant $C_P>0$ depending only on the domain $\Dref$.

From the definition of $A$ in \eqref{eq:A}, the left-hand side of
\eqref{eq:myLHS} can be written as
\begin{align*}
 &\int_{\Dref} \big\|\big(J^{-1}(\bsx,\bsy)\big)^\tr\, \nabla \widehat{u}(\bsx,\bsy)\big\|_2^2\,\det(J(\bsx,\bsy))\,{\rm d}\bsx \\
 &\ge \int_{\Dref} \big(\min \sigma \big(J^{-1}(\bsx,\bsy)^{\tr}\big) \big)^2\, \|\nabla \widehat{u}(\bsx,\bsy)\|_2^2\,
 \big(\min \sigma(J(\bsx,\bsy))\big)^d \,{\rm d}\bsx\\
 &\ge \frac{\sigma_{\min}^d}{\sigma_{\max}^2}\, \|\widehat{u}(\cdot,\bsy)\|_{H^1_0(\Dref)}^2,
\end{align*}
where we used assumptions~\ref{assumption3}
and~\eqref{eq:detJpositive}. We recall that
$\sigma_{\max}=1+\xi_{\bsb}$, see line below \eqref{eq:b-def}. %
The proof is
completed by combining the upper and lower bounds, and applying the bound
for $\|f_{\rm ref}(\cdot,\bsy)\|_{L^\infty(\Dref)}$ by taking $\bsnu =
\bszero$ in Lemma~\ref{lemma:fref}. \quad\qed

\begin{theorem} \label{thm:uhatregularity}
Under the assumptions of Lemma~\ref{lemma:fhat}, {there} holds
for all $\bsy\in U$ and $\bsnu\in\mathscr{F}\setminus\{\bszero\}$ that
\begin{align*}
 &\|\partial_\bsy^\bsnu \widehat{u}(\cdot,\bsy)\|_{H_0^1(\Dref)} \\
 &\le 2C_{\rm init}\bigg(\frac{4\pi d!(2\!+\!\xi_{\bsb})^d(1\!+\!\xi_{\bsb})^3}{\sigma_{\min}^{d+4}}\bigg)^{|\bsnu|}
  \!\!\sum_{\bszero\ne\bsm\leq\bsnu}\!
 \frac{(|\bsm|+d-1)!}{(d-1)!}\bsm!\bsbeta^{\bsm}\prod_{i\geq 1}S(\nu_i,m_i),
\end{align*}
where $\beta_j := (2+\sqrt{2})\max(1+\sqrt{3},\rho)\,b_j$.
\end{theorem}

\proof
We prove the theorem by induction on the order of $\bsnu$. The case $\bsnu
= \bszero$ holds by Lemma~\ref{lemma:apriori}. Let $\bsy\in U$ and
$\bsnu\in\mathscr{F}\setminus\{\mathbf{0}\}$, and suppose that the result
holds for all multi-indices of order less than $|\bsnu|$.

We differentiate both sides of \eqref{eq:weak2} and apply the Leibniz product
rule to obtain for all $\widehat{v}\in H_0^1(\Dref)$ that
\begin{align*}
 &\int_{\Dref} \!\!\!\bigg(
 \sum_{\bsm\le\bsnu} \binom{\bsnu}{\bsm} \big(\partial^{\bsnu-\bsm}_\bsy A(\bsx,\bsy)\big)\,
 \nabla \big(\partial^{\bsm}_\bsy \widehat{u}(\bsx,\bsy)\big)\bigg)
 \cdot \nabla \widehat{v}(\bsx)\,{\rm d}\bsx\\
 &= \int_{\Dref} \!\!\big(\partial^{\bsnu}_\bsy f_{\rm ref}(\bsx,\bsy)\big)\,\widehat{v}(\bsx)\,{\rm d}\bsx.
\end{align*}
Taking $\widehat{v}(\bsx) = \partial^\bsnu_\bsy \widehat{u}(\bsx,\bsy)$
and separating out the $\bsm=\bsnu$ term, we obtain
\begin{align*}
 &\int_{\Dref} \big(A(\bsx,\bsy)\nabla \big(\partial^\bsnu_\bsy \widehat{u}(\bsx,\bsy)\big)\big)
 \cdot \nabla \big(\partial^\bsnu_\bsy \widehat{u}(\bsx,\bsy)\big)\,{\rm d}\bsx \\
 &\,=\, - \sum_{\satop{\bsm\le\bsnu}{\bsm\ne\bsnu}} \binom{\bsnu}{\bsm}
 \int_{\Dref} \big(\big(\partial^{\bsnu-\bsm}_\bsy A(\bsx,\bsy)\big)\,
 \nabla \big(\partial^{\bsm}_\bsy \widehat{u}(\bsx,\bsy)\big)\big)
 \cdot \nabla \big(\partial^\bsnu_\bsy \widehat{u}(\bsx,\bsy)\big)\,{\rm d}\bsx \\
 &\qquad + \int_{\Dref} \big(\partial^{\bsnu}_\bsy f_{\rm ref}(\bsx,\bsy)\big)\,
 \big(\partial^\bsnu_\bsy \widehat{u}(\bsx,\bsy)\big)\,{\rm d}\bsx.
\end{align*}
Using similar steps as in the proof of Lemma~\ref{lemma:apriori}, we then
arrive at
\begin{align*}
 &\frac{\sigma_{\min}^d}{(1+\xi_\bsb)^2}\|\partial_\bsy^\bsnu\widehat{u}(\cdot,\bsy)\|_{H_0^1(\Dref)}\\
 &\le\, \sum_{\satop{\bsm\le\bsnu}{\bsm\ne\bsnu}} \binom{\bsnu}{\bsm}\,
 \|\partial^{\bsnu-\bsm}_\bsy A(\cdot,\bsy)\|_{L^\infty(\Dref)} \,
 \|\partial^{\bsm}_\bsy \widehat{u}(\cdot,\bsy)\|_{H^1_0(\Dref)} \\
 &\qquad + C_P\, |\Dref|^{1/2}\,\|\partial^{\bsnu}_\bsy f_{\rm ref}(\cdot,\bsy)\|_{L^\infty(\Dref)}.
\end{align*}
Combining this with Lemmata~\ref{lemma:Ader} and~\ref{lemma:fref}, we
obtain
\begin{align*}
 &\|\partial_\bsy^\bsnu\widehat{u}(\cdot,\bsy)\|_{H_0^1(\Dref)} \\
 &\le \frac{(1+\xi_\bsb)^2}{\sigma_{\min}^d} \sum_{\satop{\bsm\le\bsnu}{\bsm\ne\bsnu}} \binom{\bsnu}{\bsm}
 \|\partial^{\bsm}_\bsy \widehat{u}(\cdot,\bsy)\|_{H^1_0(\Dref)}\frac{d!\,(2+\xi_{\bsb} )^d}{\sigma_{\min}^2}
  \bigg(\frac{4\pi(1+\xi_{\bsb})}{\sigma_{\min}^2}\bigg)^{|\bsnu-\bsm|} \\
 &\qquad\qquad \times
 \!\!\!
  \sum_{\bsmu\le\bsnu-\bsm}(|\bsmu|+1)!\, a_{|\bsmu|}\,\bsmu!\,\bsb^{\bsmu}\prod_{i\geq 1}S(\nu_i-m_i,\mu_i) \\
 &\quad + \frac{(1+\xi_\bsb)^2}{\sigma_{\min}^d}\, C_P\, |\Dref|^{1/2}\,C_f\,d!\,(2+\xi_{\bsb})^d\,(2\pi)^{|\bsnu|}\\
&\qquad\qquad \times \sum_{\bsm\leq\bsnu} \frac{(|\bsm|+d)!}{d!}\,\rho^{|\bsm|}\,\bsb^{\bsm}\prod_{i\geq 1}S(\nu_i,m_i).
\end{align*}
With some further estimations, this can be expressed in the form of the
recursive bound in Lemma~\ref{lemma:nightmarelemma} with, using $a_k\le
(1+\sqrt{3})^k$ (which follows easily from \eqref{eq:ak}),
\begin{align*}
  \beta_j &:= \max(1+\sqrt{3},\rho)\,b_j, \\
  c_0 &:= \frac{(1+\xi_\bsb)^2}{\sigma_{\min}^d}\,C_P\, |\Dref|^{1/2}\,C_f\,d!\,(2+\xi_{\bsb})^d = C_{\rm init}, \\
  c &:= \frac{(1+\xi_\bsb)^2}{\sigma_{\min}^d}\frac{d!\,(2+\xi_{\bsb} )^d}{\sigma_{\min}^2}
  \frac{4\pi(1+\xi_{\bsb})}{\sigma_{\min}^2}.
\end{align*}
In particular, we used $|\bsnu-\bsm|\ge 1$ for the exponent in
$(\frac{4\pi(1+\xi_{\bsb})}{\sigma_{\min}^2})^{|\bsnu-\bsm|}$, which
allowed us to enlarge the constant as defined in $c$.

We then conclude from Lemmata~\ref{lemma:nightmarelemma},
\ref{lem:P-bound} and~\ref{lemma:tau-bound} that for $\bsnu\ne\bszero$
\begin{align*}
 \|\partial_\bsy^\bsnu\widehat{u}(\cdot,\bsy)\|_{H_0^1(\Dref)}\!
 \le 2c_0c^{|\bsnu|}\!\!\sum_{\bszero\ne\bsm\leq\bsnu}\!
 (2\!+\!\sqrt{2})^{|\bsm|}\frac{(|\bsm|\!+\!d\!-\!1)!}{(d\!-\!1)!}\bsm!\bsbeta^{\bsm}\!\prod_{i\geq 1}\!S(\nu_i,m_i).
\end{align*}
Redefining $\beta_j$ to incorporate the factor $2+\sqrt{2}$ completes the
proof.
\quad\qed

\section{Error analysis}\label{sec:4}
{In practice, solving~\eqref{eq:weak2} is only possible after truncating the input random field~\eqref{eq:V} to finitely many terms and discretizing the spatial domain using, e.g., finite elements. In this section, we analyze the errors incurred by \emph{dimension truncation}, \emph{finite element discretization}, and \emph{QMC cubature} for the expected value of the dimensionally-truncated, finite element solution of~\eqref{eq:weak2}. We also discuss the construction of QMC cubatures with higher order convergence rates independently of the dimension, and present a combined error estimate for the problem.}
\subsection{Dimension truncation error}\label{sec:dim}

We are interested in finding a rate for
$$
\bigg\|\int_U (\widehat u(\cdot,\bsy)-\widehat u_s(\cdot,\bsy))\,\rd\bsy\bigg\|_{L^1(\Dref)},
$$
where $\widehat u$ denotes the solution to~\eqref{eq:weak2} as before and
the dimensionally-truncated PDE solution is defined by $$\widehat
u_s(\cdot,\bsy):=\widehat
u(\cdot,(y_1,\ldots,y_s,0,0,\ldots))\quad\text{for}~\bsy\in U,~s\in\mathbb
N.$$ In order to derive a dimension-truncation error bound for this
problem, we employ the Taylor series approach used
in~\cite{GGKSS2019,guth22} together with a change of variable technique.
(The Neumann series approach used in, e.g., \cite{gantner} is not
applicable here because the transported problem cannot be recast in affine
parametric form.)

\begin{theorem}\label{thm:dimtrunc}
{Suppose that the assumptions~{\rm\ref{assumption1}}--\,{\rm\ref{assumption6}} hold}. Then the di\-mensionally-truncated solution to~\eqref{eq:weak2} satisfies the error bounds
$$
\bigg\|\int_U (\widehat u(\cdot,\bsy)-\widehat u_s(\cdot,\bsy))\,{\rm d}\bsy\bigg\|_{H_0^1(D_{\rm ref})}=\mathcal O(s^{-2/p+1})
$$
and
$$
\bigg\|\int_U (\widehat u(\cdot,\bsy)-\widehat u_s(\cdot,\bsy))\,{\rm d}\bsy\bigg\|_{L^1(D_{\rm ref})}=\mathcal O(s^{-2/p+1}),
$$
where the implied coefficients are independent of the dimension $s$.
\end{theorem}
\proof We start by defining
$$
\widetilde{\bspsi}_i(\bsx)\,:=\,\frac{1}{\sqrt 6}\boldsymbol\psi_i(\bsx),\quad i\in\mathbb N,
$$
where the fluctuations $(\|\boldsymbol\psi_i\|_{W^{1,\infty}})_{i\geq 1}\in \ell^p$ for the same $p$ as in~\ref{assumption5}. Let us consider an auxiliary PDE problem
\begin{align}\label{eq:auxpde}
\begin{cases}
-\Delta u_{\rm aff}(\bsx,\bsy)=f(\bsx),&\bsx\in \widetilde D(\bsy),\\
u_{\rm aff}(\cdot,\bsy)|_{\partial \widetilde D(\bsy)}=0,
\end{cases}\quad\text{for all}~\bsy\in[-1,1]^{\mathbb N},
\end{align}
where the domain realizations $\widetilde D(\bsy):=\widetilde{\boldsymbol{V}}(D_{\rm ref},\bsy)$ are defined instead via an \emph{affine parametric} random perturbation field
$$
\widetilde{\boldsymbol V}(\bsx,\bsy)\,:=\,\bsx+\sum_{i\geq 1}y_i \widetilde{\boldsymbol\psi}_i(\bsx)
$$
with $(\|\widetilde{\boldsymbol\psi}_i\|_{W^{1,\infty}})_{i\geq 1}\in \ell^p$ for the same $p$ as in~\ref{assumption5}. The random field $\widetilde{\boldsymbol{V}}$ satisfies the hypotheses of~\cite[Theorem~5]{nummat2016} ({note that the random variables
in~\cite{nummat2016} {are supported on} $[-1,1]$ instead of $[-1/2,1/2]$}) and the transported PDE solution corresponding to~\eqref{eq:auxpde} satisfies
$$
\|\partial_\bsy^{\bsnu}\widehat u_{\rm aff}(\cdot,\bsy)\|_{H_0^1(D_{\rm ref})}\leq C\,|\bsnu|!\,\widetilde  \bsb^{\bsnu},
$$
for some $\widetilde \bsb\in\ell^p$, {where the $p$ is the same as above}.

Let $G\in H^{-1}(D_{\rm ref})$ be arbitrary and define
$$
\widetilde F(\bsy)\,:=\,\langle G,\widehat u_{\rm aff}(\cdot,\bsy)\rangle_{H^{-1}(D_{\rm ref}),H_0^1(D_{\rm ref})}.
$$
Clearly,
$$
\partial_\bsy^{\bsnu}\widetilde F(\bsy)=\langle G,\partial_\bsy^{\bsnu}\widehat u_{\rm aff}(\cdot,\bsy)\rangle_{H^{-1}(D_{\rm ref}),H_0^1(D_{\rm ref})}.
$$
Fix $\bsy\in [-1,1]^{\mathbb N}$. Then developing the Taylor series of $\widetilde F$ about the point $(\bsy_{\leq s},\mathbf 0):=(y_1,\ldots,y_s,0,0,\ldots)$ yields
\begin{align*}
\widetilde F(\bsy)-\widetilde F(\bsy_{\leq s},\mathbf 0)&=\sum_{\ell=1}^k \sum_{\substack{|\bsnu|=\ell\\ \nu_i=0~\forall i\leq s}}\frac{\bsy^{\bsnu}}{\bsnu!}\partial_\bsy^{\bsnu}\widetilde F(\bsy_{\leq s},\mathbf 0)\\
&\quad +\sum_{\substack{|\bsnu|=k+1\\ \nu_i=0~\forall i\leq s}}\frac{k+1}{\bsnu!}\bsy^{\bsnu}\int_0^1(1-t)^k\partial_\bsy^{\bsnu}\widetilde F(\bsy_{\leq s},t\bsy_{>s})\,{\rm d}t.
\end{align*}

Let $k\geq 2$ be an undetermined integer which will be specified later. By defining $\boldsymbol\theta(\bsy):=(\sin(2\pi y_i))_{i\geq 1}$, we observe that %
$$
{\widehat{u}_{\rm per}}(\cdot,\bsy)=\widehat u_{\rm aff}(\cdot,\boldsymbol\theta(\bsy)),
$$
{where $\widehat{u}_{\rm per}(\cdot,\bsy)$ is the smooth periodic extension of the solution to~\eqref{eq:weak2} for all $\bsy\in[-1,1]^{\mathbb N}$.} Let us define $F(\bsy):=\langle G,\widehat u_{\rm per}(\cdot,\bsy)\rangle_{H^{-1}(D_{\rm ref}),H_0^1(D_{\rm ref})}$. We obtain by a simple change of variable that
\begin{align*}
F(\bsy)-F(\bsy_{\leq s},\mathbf 0)&=\sum_{\ell=1}^k \sum_{\substack{|\bsnu|=\ell\\ \nu_i=0~\forall i\leq s}}\frac{\boldsymbol\theta(\bsy)^{\bsnu}}{\bsnu!}\partial_\bsy^{\bsnu}\widetilde F(\boldsymbol\theta(\bsy_{\leq s}),\mathbf 0)\\
&\quad +\sum_{\substack{|\bsnu|=k+1\\ \nu_i=0~\forall i\leq s}}\frac{k+1}{\bsnu!}\boldsymbol\theta(\bsy)^{\bsnu}\int_0^1(1-t)^k\partial_\bsy^{\bsnu}\widetilde F(\boldsymbol\theta(\bsy_{\leq s},t\bsy_{>s}))\,{\rm d}t.
\end{align*}

Since the last equality holds pointwise for all $\bsy\in[-1,1]^{\mathbb N}$, we can integrate the equation on both sides over $U\subset [-1,1]^{\mathbb N}$. We note that
the integral $\int_{0}^1 \sin(2\pi y_i)^{\nu_i}\,\rd y_i$ takes the
value between $0$ and $1$, and gives the value $0$ when $\nu_i = 1$. Thus
the summands corresponding to $\bsnu\in\mathscr F$ with $\nu_i=1$ for
at least one $i>s$ vanish in the first term. Hence
\begin{align*}
&\int_U (F(\bsy)-F(\bsy_{\leq s},\mathbf 0))\,{\rm d}\bsy\\
&=\sum_{\ell=2}^k \sum_{\substack{|\bsnu|=\ell\\ \nu_i=0~\forall i\leq s\\ \nu_i\neq 1~\forall i>s}}\frac{1}{\bsnu!}\int_U \boldsymbol\theta(\bsy)^{\bsnu}\partial_\bsy^{\bsnu}\widetilde F(\boldsymbol\theta(\bsy_{\leq s}),\mathbf 0)\,{\rm d}\bsy\\
&\quad +\sum_{\substack{|\bsnu|=k+1\\ \nu_i=0~\forall i\leq s}}\frac{k+1}{\bsnu!}\int_U\int_0^1(1-t)^k\boldsymbol\theta(\bsy)^{\bsnu}\partial_\bsy^{\bsnu}\widetilde F(\boldsymbol\theta(\bsy_{\leq s},t\bsy_{>s}))\,{\rm d}t\,{\rm d}\bsy\\
&\leq C'\,\|G\|_{H^{-1}(D_{\rm ref})}\sum_{\ell=2}^k \sum_{\substack{|\bsnu|=\ell\\ \nu_i=0~\forall i\leq s\\ \nu_i\neq 1~\forall i>s}}\widetilde \bsb^{\bsnu}
+ C''\,\|G\|_{H^{-1}(D_{\rm ref})}\sum_{\substack{|\bsnu|=k+1\\ \nu_i=0~\forall i\leq s}}\widetilde \bsb^{\bsnu},
\end{align*}
where $C'=C\,k!$ and $C''=C\,(k+1)!$. Recalling the
definition of $F$ and taking the supremum over all $G\in H^{-1}(D_{\rm
ref})$ with $\|G\|_{H^{-1}(D_{\rm ref})}\leq 1$ yields for the left-hand
side that
\begin{align*}
&\sup_{\substack{G\in H^{-1}(D_{\rm ref})\\ \|G\|_{H^{-1}(D_{\rm ref})}\leq 1}}\int_U \langle G,\widehat u(\cdot,\bsy)-\widehat u_s(\cdot,\bsy)\rangle_{H^{-1}(D_{\rm ref}),H_0^1(D_{\rm ref})}\,{\rm d}\bsy\\
&=\sup_{\substack{G\in H^{-1}(D_{\rm ref})\\ \|G\|_{H^{-1}(D_{\rm ref})}\leq 1}}\bigg\langle G,\int_U (\widehat u(\cdot,\bsy)-\widehat u_s(\cdot,\bsy))\,{\rm d}\bsy\bigg\rangle_{H^{-1}(D_{\rm ref}),H_0^1(D_{\rm ref})}\\
&=\sup_{\substack{G\in H^{-1}(D_{\rm ref})\\ \|G\|_{H^{-1}(D_{\rm ref})}\leq 1}}\bigg\langle G,\int_U (\widehat u(\cdot,\bsy)-\widehat u_s(\cdot,\bsy))\,{\rm d}\bsy\bigg\rangle_{H^{-1}(D_{\rm ref}),H_0^1(D_{\rm ref})}\\
&=\bigg\|\int_U (\widehat u(\cdot,\bsy)-\widehat u_s(\cdot,\bsy))\,{\rm d}\bsy\bigg\|_{H_0^1(D_{\rm ref})},
\end{align*}
which yields the overall bound
\begin{align} \label{eq:two}
 \bigg\|\int_U (\widehat u(\cdot,\bsy)-\widehat u_s(\cdot,\bsy))\,{\rm d}\bsy\bigg\|_{H_0^1(D_{\rm ref})}
 \lesssim
 \underbrace{\sum_{\ell=2}^k \sum_{\substack{|\bsnu|=\ell\\ \nu_i=0~\forall i\leq s\\ \nu_i\neq 1~\forall i>s}}\widetilde \bsb^{\bsnu}}_{=:\, T_1}
 + \underbrace{\sum_{\substack{|\bsnu|=k+1\\ \nu_i=0~\forall i\leq s \\ \vphantom{\forall j}} }\widetilde \bsb^{\bsnu}}_{=:\, T_2},
\end{align}%
where the generic constant is independent of $s$. The goal will be to
balance $T_1$ and $T_2$ in \eqref{eq:two} by choosing the value of
$k$ appropriately.

The term $T_1$ in \eqref{eq:two} can be bounded similarly
to~\cite[Theorem~1]{gantner}:
\begin{align*}
 T_1
 &=\sum_{\substack{2\leq|\bsnu|\leq k\\ \nu_i=0~\forall i\leq s\\ \nu_i\neq 1~\forall i>s}}\widetilde\bsb^{\bsnu}\leq \sum_{\substack{0\neq|\bsnu|_\infty\leq k\\ \nu_i=0~\forall i\leq s\\ \nu_i\neq 1~\forall i>s}}\widetilde\bsb^{\bsnu}
 =-1+\prod_{i>s}\bigg(1+\sum_{\ell=2}^k\widetilde b_i^\ell\bigg)\\
 &=-1+\prod_{i>s}\bigg(1+\frac{1-\widetilde b_i^{k-1}}{1-\widetilde b_i}\widetilde b_i^2\bigg)
 \le -1+\prod_{i>s}\bigg(1+ \tau_k\, \widetilde b_i^2\bigg),
\end{align*}
where we abuse notation slightly by interpreting the value of the term
$\frac{1-\widetilde b_i^{k-1}}{1-\widetilde b_i}$ as $k-1$ if $\widetilde
b_i=1$, and since the sequence is $\widetilde\bsb$ is non-increasing
we defined $\tau_k := \frac{1-\widetilde{b}_1^{k-1}}{1-\widetilde{b}_1}$ if $\widetilde b_1\neq 1$ and $\tau_k:=k-1$ if $\widetilde b_1=1$.  Moreover, it is an easy consequence of~\ref{assumption6} that
$$
\sum_{i>s}\widetilde b_i^2\leq\bigg(\sum_{i\geq 1}\widetilde b_i^p\bigg)^{2/p}s^{-2/p+1},
$$
where we used the well-known Stechkin's lemma (cf., e.g.,~\cite[Lemma 3.3]{KressnerTobler}): for $0<p\leq q<\infty$ and any sequence $(a_i)_{i\geq 1}$ of real numbers ordered such that $|a_1|\geq|a_2|\geq\cdots$, there holds
\begin{align}
\bigg(\sum_{i>s} |a_i|^q\bigg)^{1/q}\leq \bigg(\sum_{i\geq 1}|a_i|^p\bigg)^{1/p}s^{-1/p+1/q}.\label{eq:stechkin}
\end{align}
Therefore
\begin{align*}
 T_1
&\leq -1+\exp\bigg(\tau_k\sum_{i>s}\widetilde b_i^2\bigg)
 = \sum_{j=1}^\infty \frac{\tau_k^j}{j!}\bigg(\sum_{i>s}\widetilde b_i^2\bigg)^j \\
&\leq\bigg(\exp\bigg(\tau_k\bigg(\sum_{i\geq 1}\widetilde b_i^p\bigg)^{2/p}\bigg)-1\bigg) s^{-2/p+1}.
\end{align*}

The term $T_2$ in \eqref{eq:two} can be estimated similarly to the
approach taken in~\cite[Theorem~4.1]{GGKSS2019}. The trivial inequality
$\frac{|\bsnu|!}{\bsnu!}\geq 1$ and the multinomial theorem imply together
that
\begin{align*}
 T_2
\leq \sum_{\substack{|\bsnu|=k+1\\ \nu_i=0~\forall i\leq s}}\frac{|\bsnu|!}{\bsnu!}\widetilde\bsb^{\bsnu}=\bigg(\sum_{i>s}\widetilde b_i\bigg)^{k+1}
\leq \bigg(\sum_{i\geq 1}\widetilde b_i^p\bigg)^{(k+1)/p}s^{(k+1)(-1/p+1)},
\end{align*}
where the final inequality is a consequence of~\eqref{eq:stechkin}.
The exponent of $s$ in the bound on $T_2$ is dominated by that of
$T_1$ by choosing $k=\lceil 1/(1-p)\rceil$.

Finally, the Poincar\'e inequality implies that
$$
\bigg\|\int_U (\widehat u(\cdot,\bsy)-\widehat u_s(\cdot,\bsy))\,{\rm d}\bsy\bigg\|_{L^2(D_{\rm ref})}=\mathcal O(s^{-2/p+1}),
$$
and the $L^1$ error can be obtained by using the Cauchy--Schwarz
inequality
\begin{align*}
&\bigg\|\int_U (\widehat u(\cdot,\bsy)-\widehat u_s(\cdot,\bsy))\,{\rm d}\bsy\bigg\|_{L^1(D_{\rm ref})}\\
&\leq |D_{\rm ref}|^{1/2}\bigg\|\int_U (\widehat u(\cdot,\bsy)-\widehat u_s(\cdot,\bsy))\,{\rm d}\bsy\bigg\|_{L^2(D_{\rm ref})}.
\end{align*}
This completes the proof.\quad\qed

\subsection{Finite element error}\label{sec:fem}
{Let the assumptions \ref{assumption1}--\ref{assumption4} and~\ref{assumption7} hold.} Let us denote by $\{V_h\}_h$ a family of finite element subspaces of $H_0^1(D_{\rm ref})$ indexed by the mesh size $h$. We assume that the finite element spaces are spanned by continuous, piecewise linear finite element basis functions and the mesh corresponding to each $V_h$ is obtained from an initial, regular triangulation of $D_{\rm ref}$ by recursive uniform bisection of simplices.

We define the dimensionally-truncated finite element solution for $\bsy\in U$ as the solution $u_{s,h}(\cdot,\bsy)\in V_h$ to
$$
\int_{D_{\rm ref}} (A_s(\bsx,\bsy)\nabla \widehat u_{s,h}(\bsx,\bsy))\cdot \nabla\widehat v_h(\bsx)\,{\rm d}\bsx=\int_{D_{\rm ref}}f_{{\rm ref},s}(\bsx,\bsy)\widehat v_h(\bsx)\,{\rm d}\bsx
$$
for all $\widehat v_h\in V_h$, where we {set} $A_s(\bsx,\bsy):=A(\bsx,(y_1,\ldots,y_s,0,0,\ldots))$ {as well as} $f_{\rm ref,s}(\bsx,\bsy):=f_{\rm ref}(\bsx,(y_1,\ldots,y_s,0,0,\ldots))$ for all $\bsy\in U$. Since the reference domain $D_{\rm ref}$ was assumed in~\ref{assumption7} to be a bounded, convex polyhedron, it follows that elements in $W^{1,\infty}(D_{\rm ref})$ have Lipschitz continuous representatives. Hence we may apply Lemma~\ref{lemma:fref} and~\cite[Theorem~3.2.1.2]{grisvard1985elliptic} to obtain
\begin{align}
\|\widehat u_s(\cdot,\bsy)\|_{H^2(D_{\rm ref})\cap H_0^1(D_{\rm ref})}\leq C,\label{eq:regularityshift}
\end{align}
where $\|\widehat v\|_{H^2(D_{\rm ref})\cap H_0^1(D_{\rm ref})}:=(\|\widehat v\|_{L^2(D_{\rm ref})}^2+\|\Delta \widehat v\|_{L^2(D_{\rm ref})}^2)^{1/2}$ and the constant $C>0$ only depends on $d$ and the diameter of $D_{\rm ref}$.

We have by C\'ea's lemma that
\begin{align}
\|\widehat u_s(\cdot,\bsy)-\widehat u_{s,h}(\cdot,\bsy)\|_{H_0^1(D_{\rm ref})}\leq C' \inf_{\widehat v_h\in V_h}\|\widehat u_s(\cdot,\bsy)-\widehat v_h\|_{H_0^1(D_{\rm ref})},\label{eq:cea}
\end{align}
for a constant $C'>0$ independent of $s$, $h$, and $\bsy$, which can be seen by letting $\widehat v_h\in V_h$ be arbitrary and deriving
\begin{align*}
&\frac{\sigma_{\min}^d}{\sigma_{\max}^2}\|\widehat u_s(\cdot,\bsy)-\widehat u_{s,h}(\cdot,\bsy)\|_{H_0^1(D_{\rm ref})}^2\\
&\leq \int_{D_{\rm ref}}(A_s(\bsx,\bsy)\nabla(\widehat u_s(\bsx,\bsy)-\widehat u_{s,h}(\bsx,\bsy)))\cdot \nabla(\widehat u_s(\bsx,\bsy)-\widehat u_{s,h}(\bsx,\bsy))\,{\rm d}\bsx\\
&= \int_{D_{\rm ref}}(A_s(\bsx,\bsy)\nabla(\widehat u_s(\bsx,\bsy)-\widehat u_{s,h}(\bsx,\bsy)))\cdot \nabla(\widehat u_s(\bsx,\bsy)-\widehat v_h(\bsx))\,{\rm d}\bsx\\
&\quad +\int_{D_{\rm ref}}(A_s(\bsx,\bsy)\nabla(\widehat u_s(\bsx,\bsy)-\widehat u_{s,h}(\bsx,\bsy)))\cdot \nabla(\widehat v_h(\bsx)-\widehat u_{s,h}(\bsx,\bsy))\,{\rm d}\bsx\\
&= \int_{D_{\rm ref}}(A_s(\bsx,\bsy)\nabla(\widehat u_s(\bsx,\bsy)-\widehat u_{s,h}(\bsx,\bsy)))\cdot \nabla(\widehat u_s(\bsx,\bsy)-\widehat v_h(\bsx))\,{\rm d}\bsx\\
&\leq \frac{\sigma_{\max}^d}{\sigma_{\min}^2}\|\widehat u_s(\cdot,\bsy)-\widehat u_{s,h}(\cdot,\bsy)\|_{H_0^1(D_{\rm ref})}\|\widehat u_{s}(\cdot,\bsy)-\widehat v_h\|_{H_0^1(D_{\rm ref})},
\end{align*}
where the penultimate step follows from Galerkin orthogonality.

The well-known approximation property (cf., e.g.,~\cite{Gilbarg,Kuo2012})
\begin{align}
\inf_{\widehat v_h\in V_h}\|\widehat v-\widehat v_h\|_{H_0^1(D_{\rm ref})}\leq C''h\|\widehat v\|_{H^2(D_{\rm ref})\cap H_0^1(D_{\rm ref})}\quad\text{as}~h\to 0\label{eq:feapprox}
\end{align}
holds for all $\widehat v\in H^2(D_{\rm ref})\cap H_0^1(D_{\rm ref})$ with a constant $C''>0$ independent of $h$. The inequalities~\eqref{eq:regularityshift}, \eqref{eq:cea}, and~\eqref{eq:feapprox} yield altogether that
\begin{align}
\|\widehat u_s(\cdot,\bsy)-\widehat u_{s,h}(\cdot,\bsy)\|_{H_0^1(D_{\rm ref})}\leq C'''h\quad\text{as}~h\to 0,\label{eq:FE_Verror}
\end{align}
where the constant $C'''>0$ is independent of $s$, $h$, and $\bsy$.

An $L^1$ error estimate can be obtained from a standard Aubin--Nitsche duality argument.
\begin{theorem}
{Under assumptions~{\rm\ref{assumption1}}--\,{\rm\ref{assumption4}} and {\rm\ref{assumption7}}}, {there} holds for the dimensionally-truncated finite element solution that
$$
\|\widehat u_s(\cdot,\bsy)-\widehat u_{s,h}(\cdot,\bsy)\|_{L^1(D_{\rm ref})}\lesssim h^2\quad\text{as}~h\to0,
$$
where the {generic constant} is independent of $s$, $h$, and $\bsy$.
\end{theorem}
\proof Since
$$
\|\widehat u_s(\cdot,\bsy)-\widehat u_{s,h}(\cdot,\bsy)\|_{L^1(D_{\rm ref})}\leq |D_{\rm ref}|^{1/2}\|\widehat u_s(\cdot,\bsy)-\widehat u_{s,h}(\cdot,\bsy)\|_{L^2(D_{\rm ref})},
$$
it suffices to prove the $L^2$ error estimate.

Let $\widehat g\in L^2(D_{\rm ref})$. For $\bsy\in U$, denote by $\widehat u_{\widehat g,s}(\cdot,\bsy)\in H_0^1(D_{\rm ref})$ the solution to
$$
\int_{D_{\rm ref}} (A_s(\bsx,\bsy)\nabla \widehat u_{\widehat g,s}(\bsx,\bsy))\cdot \nabla\widehat v(\bsx)\,{\rm d}\bsx=\langle\widehat g,\widehat v\rangle_{L^2(D_{\rm ref})}\quad\text{for all}~\widehat v\in H_0^1(D_{\rm ref}).
$$
Testing against $\widehat v=\widehat u_s(\cdot,\bsy)-\widehat u_{s,h}(\cdot,\bsy)$ and letting $\widehat v_h\in V_h$ be arbitrary, it follows from Galerkin orthogonality of the finite element solution that
\begin{align*}
&\langle g,\widehat u_s(\cdot,\bsy)-\widehat u_{s,h}(\cdot,\bsy)\rangle_{L^2(D_{\rm ref})}\\
&=\int_{D_{\rm ref}}(A_s(\bsx,\bsy)\nabla \widehat u_{\widehat g,s}(\bsx,\bsy))\cdot \nabla (\widehat u_s(\bsx,\bsy)-\widehat u_{s,h}(\bsx,\bsy))\,{\rm d}\bsx\\
&=\int_{D_{\rm ref}}(A_s(\bsx,\bsy)\nabla (\widehat u_{\widehat g,s}(\bsx,\bsy)-\widehat v_h(\bsx)))\cdot \nabla (\widehat u_s(\bsx,\bsy)-\widehat u_{s,h}(\bsx,\bsy))\,{\rm d}\bsx\\
&\leq \frac{\sigma_{\max}^d}{\sigma_{\min}^2}\|\widehat u_{\widehat g,s}(\cdot,\bsy)-\widehat v_h\|_{H_0^1(D_{\rm ref})}\|\widehat u_s(\cdot,\bsy)-\widehat u_{s,h}(\cdot,\bsy)\|_{H_0^1(D_{\rm ref})}
\end{align*}
and, in consequence,
\begin{align*}
&\langle g,\widehat u_s(\cdot,\bsy)-\widehat u_{s,h}(\cdot,\bsy)\rangle_{L^2(D_{\rm ref})}\\
&\leq \frac{\sigma_{\max}^d}{\sigma_{\min}^2}\|\widehat u_s(\cdot,\bsy)-\widehat u_{s,h}(\cdot,\bsy)\|_{H_0^1(D_{\rm ref})}\inf_{\widehat v_h\in V_h}\|\widehat u_{\widehat g,s}(\cdot,\bsy)-\widehat v_h\|_{H_0^1(D_{\rm ref})}.
\end{align*}
Then {there} holds for all $\bsy\in U$ that
\begin{align}
&\|\widehat u_s(\cdot,\bsy)-\widehat u_{s,h}(\cdot,\bsy)\|_{L^2(D_{\rm ref})}\notag\\
&=\sup_{\substack{\widehat g\in L^2(D_{\rm ref})\\ \|\widehat g\|_{L^2(D_{\rm ref})}\leq 1}}\langle \widehat g,\widehat u_s(\cdot,\bsy)-\widehat u_{s,h}(\cdot,\bsy)\rangle_{L^2(D_{\rm ref})}\notag\\
&\leq \frac{\sigma_{\max}^d}{\sigma_{\min}^2}\|\widehat u_s(\cdot,\bsy)\!-\!\widehat u_{s,h}(\cdot,\bsy)\|_{H_0^1(D_{\rm ref})}\!\sup_{\substack{\widehat g\in L^2(D_{\rm ref})\\ \|\widehat g\|_{L^2(D_{\rm ref})}\leq 1}}\!\inf_{\widehat v_h\in V_h}\|\widehat u_{\widehat g,s}(\cdot,\bsy)\!-\!\widehat v_h\|_{H_0^1(D_{\rm ref})}.\label{eq:prelimfeerror}
\end{align}
By~\eqref{eq:feapprox} and \eqref{eq:regularityshift}, we obtain
\begin{align*}
\inf_{\widehat v_h\in V_h}\|\widehat u_{\widehat g,s}(\cdot,\bsy)-\widehat v_h\|_{H_0^1(D_{\rm ref})}&\leq C''h\|\widehat u_{\widehat g,s}(\cdot,\bsy)\|_{H^2(D_{\rm ref})\cap H_0^1(D_{\rm ref})}\\
&\leq CC''h\|\widehat g\|_{L^2(D_{\rm ref})}\quad\text{as}~h\to 0.
\end{align*}
Plugging this inequality and~\eqref{eq:FE_Verror} into~\eqref{eq:prelimfeerror} yields the desired result.\quad\qed

\subsection{Quasi-Monte Carlo cubature error and a choice of weights}\label{sec:qmc}

Let us consider {QMC} cubature over the $s$-dimensional
unit cube $U_s:=[0,1]^s$ for finite $s\in\mathbb{N}$. An $n$-point QMC
rule is an equal weight cubature rule of the form
$$
\int_{U_s}F(\bsy)\,\rd\bsy\approx \frac{1}{n}\sum_{i=1}^nF(\bsy^{(i)}),
$$
where the cubature nodes $\bsy^{(1)},\ldots,\bsy^{(n)}\in [0,1]^s$ are
deterministically chosen. In the subsequent analysis, we consider the
family of \emph{rank-1 lattice rules}, where the QMC nodes in $[0,1]^s$
are defined by
\begin{align*}
\boldsymbol{y}^{(i)}:= \frac{i\boldsymbol{z}\bmod n}{n}\quad\text{for}~i\in\{1,\ldots,n\},%
\end{align*}
where $\boldsymbol{z}\in\{1,\ldots,n-1\}^s$ is called the \emph{generating
vector}.

{First we briefly recall the theory of lattice rules for $1$-periodic
integrands with absolutely convergent Fourier series: for
\[
  F(\bsy) = \sum_{\bsell\in\bbZ^s} \calF\{F\}_\bsell\,\re^{2\pi\ri\bsell\cdot\bsy},
  \quad\mbox{with}\quad
  \calF\{F\}_\bsell := \int_{U_s}F(\bsy)\,\re^{-2\pi\ri\bsell\cdot\bsy}\,\rd\bsy,
\]
the lattice rule error is given by~\cite{sloankachoyan87}
\begin{align}
 \frac{1}{n}\sum_{i=1}^n F(\bsy^{(i)})-\int_{U_s}F(\bsy)\,\rd\bsy
 = \sum_{\bsell\in\Lambda^\perp\setminus\{\bszero\}} \calF\{F\}_\bsell,\label{eq:latticeerror}
\end{align}
where the \emph{dual lattice} $\Lambda^\perp :=\{\bsell\in\mathbb Z^s:
\bsell\cdot \bsz\equiv 0~({\rm mod}~n)\}$ depends on the lattice
generating vector $\bsz$. For integrands belonging to the weighted Korobov
space endowed with the norm
\begin{align} \label{eq:r}
 \|F\|_{\alpha,\bsgamma} :=
 \sup_{\bsell\in\bbZ^s}|\calF\{F\}_\bsell|\,r_{\alpha,\bsgamma}(\bsell),
 \quad\mbox{with}\quad
 r_{\alpha,\bsgamma}(\bsell) :=
 \gamma_{{\supp}(\bsell)}^{-1}\prod_{j\in{\supp}(\bsell)}|\ell_j|^\alpha,
\end{align}
where $\alpha>1$ is the smoothness parameter,
$\bsgamma:=(\gamma_{\setu})_{\setu\subseteq\{1:s\}}$ are positive weights,
and $\supp(\bsell) := \{1\le j\le s: \ell_j\ge 0\}$, we then have the
integration error bound
\begin{align*}
 \bigg|\int_{U_s}F(\bsy)\,\rd\bsy - \frac{1}{n}\sum_{i=1}^n F(\bsy^{(i)})\bigg|
 \le \|F\|_{\alpha,\bsgamma}\,
 \sum_{\bsell\in\Lambda^{\perp}\setminus\{\mathbf 0\}}\frac{1}{r_{\alpha,\bsgamma}(\bsell)}.
\end{align*}
It is known~\cite[Theorem~5]{korobovpaper} that the
\emph{component-by-component} (\emph{CBC}) algorithm can be used to find a
generating vector satisfying the bound
\begin{align} \label{eq:cbc}
 \sum_{\bsell\in\Lambda^{\perp}\setminus\{\mathbf 0\}}\frac{1}{r_{\alpha,\bsgamma}(\bsell)}
 \leq \bigg(\frac{1}{n-1}\underset{=:\,C_{\alpha,\bsgamma}(\lambda,s)}{\underbrace{\sum_{\varnothing\neq\setu\subseteq\{1:s\}}\gamma_{\setu}^{\lambda}(2\zeta(\alpha\lambda))^{|\setu|}}}\bigg)^{1/\lambda}
 \quad\mbox{for all}\quad \lambda\in (\tfrac{1}{\alpha},1],
\end{align}
where $n$ was assumed to be prime in \cite{korobovpaper}, but the result
also generalizes to non-primes by replacing $n-1$ by the Euler totient
function.

In the following we will adapt this theory to the particular integrand
$F(\bsy) = \widehat u_{s,h}(\bsx,\bsy)$ for $\bsx\in\Dref$. We derive
\emph{smoothness-driven product-and-order dependent (SPOD) weights} and
show that the error can be bounded independent of $s$; weights of this
form have previously appeared in~\cite{spodpaper14,periodicpaper}.

\begin{theorem}
{Suppose that assumptions~{\rm\ref{assumption1}}--\,{\rm\ref{assumption5}} and {\rm\ref{assumption7}} hold.} Let $\alpha := \lfloor 1/p\rfloor+1$ with $p$ from
Assumption~\textnormal{\ref{assumption5}}, and choose the weights
$\gamma_{\varnothing}:=1$ and
\begin{align}
 \gamma_{\setu}:=\sum_{\bsm_{\setu}\in\{1:\alpha\}^{|\setu|}}
 \frac{(|\bsm_{\setu}|+d-1)!}{(d-1)!}\prod_{j\in\setu}(\widetilde C^\alpha\, m_j!\,\beta_j^{m_j}\,S(\alpha,m_j)),~\varnothing\neq\setu\subseteq\{1:s\},\label{eq:spodweights}
\end{align}
where $\widetilde
C:=\frac{2\,d!\,(2+\xi_{\bsb})^d(1+\xi_{\bsb})^3}{\sigma_{\min}^{d+4}}$
and $\beta_j:=(2+\sqrt 2)\max(1+\sqrt 3,\rho)\,b_j$. Then a lattice rule
can be constructed by the CBC algorithm such that
\begin{align*}
 {\rm error} \,:=\, \int_{\Dref}\bigg|\int_{U_s} \widehat u_{s,h}(\bsx,\bsy)\,\rd\bsy
 -\frac{1}{n}\sum_{i=1}^n \widehat u_{s,h}(\bsx,\bsy^{(i)})\bigg|\,\rd\bsx
 \,=\, \calO(n^{-1/p}),
\end{align*}
where the implied constant is {independent} of $s$.
\end{theorem}

\proof
We will begin by treating $\alpha$ and the weights $\gamma_\setu$ as
arbitrary and will specify their values later in the proof. We will assume
that $\alpha\ge 2$ is an integer, and we will express the integration
error in term of the mixed derivatives of $F(\bsy) = \widehat
u_{s,h}(\bsx,\bsy)$ of order up to~$\alpha$ in each coordinate. For
$\setu\subseteq\{1:s\}$ we denote by $\bsnu = (\bsalpha_\setu;\bszero)$
the multi-index whose component $\nu_j$ is equal to $\alpha$ if $j\in\setu$
and is $0$ otherwise. First we observe by differentiating the Fourier
series of $F$ that
\[
  \partial^{(\bsalpha_\setu;\bszero)}_\bsy F(\bsy)
  = \sum_{\bsell\in\bbZ^s} \Big(\prod_{j\in\setu} (2\pi\ri\,\ell_j)^\alpha\Big)\,
  \calF\{F\}_\bsell\,\re^{2\pi\ri\bsell\cdot\bsy},
\]
from which we conclude that the Fourier coefficients of
$\partial^{(\bsalpha_\setu;\bszero)}_\bsy F$ are given by
\begin{align} \label{eq:der-coef}
  \int_{U_s} \partial^{(\bsalpha_\setu;\bszero)}_\bsy F(\bsy)\,\re^{-2\pi\ri\bsell\cdot\bsy}\,\rd\bsy
  = \Big(\prod_{j\in\setu} (2\pi\ri\,\ell_j)^\alpha\Big)\,
  \calF\{F\}_\bsell.
\end{align}

We are ready to derive the integration error for $\widehat u_{s,h}$.
Using~\eqref{eq:latticeerror} we have
\begin{align} \label{eq:step1}
{\rm error}
 &=\int_{\Dref}\bigg|\sum_{\bsell\in\Lambda^{\perp}\setminus\{\mathbf 0\}}
 \calF\{\widehat u_{s,h}(\bsx,\cdot)\}_\bsell\bigg|\,\rd\bsx \nonumber\\
&\leq \bigg(\sup_{\bsell\in\mathbb Z^s}\int_{\Dref} r_{\alpha,\bsgamma}(\bsell)\,
|\calF\{\widehat u_{s,h}(\bsx,\cdot)\}_\bsell|\,\rd\bsx\bigg)
\sum_{\bsell\in\Lambda^{\perp}\setminus\{\mathbf 0\}}\frac{1}{r_{\alpha,\bsgamma}(\bsell)}.
\end{align}
Denoting $\setu = \supp(\bsell)$, we obtain using the definition of
$r_{\alpha,\bsgamma}(\bsell)$ in \eqref{eq:r} with \eqref{eq:der-coef}
\begin{align} \label{eq:step2}
 & \int_{\Dref} r_{\alpha,\bsgamma}(\bsell)\,
  |\calF\{\widehat u_{s,h}(\bsx,\cdot)\}_\bsell|\,\rd\bsx
 = \frac{1}{\gamma_{\setu}}
  \int_{\Dref} \Big(\prod_{j\in\setu}|\ell_j|^{\alpha}\Big)
  |\calF\{ \widehat u_{s,h}(\bsx,\cdot)\}_\bsell|\,\rd\bsx \nonumber\\
&= \frac{1}{\gamma_{\setu}(2\pi)^{\alpha|\setu|}}
  \int_{\Dref} \Big(\prod_{j\in\setu}|2\pi\ri\,\ell_j|^{\alpha}\Big)
  |\calF\{ \widehat u_{s,h}(\bsx,\cdot)\}_\bsell|\,\rd\bsx \nonumber\\
&= \frac{1}{\gamma_{\setu}(2\pi)^{\alpha|\setu|}}
  \int_{\Dref}\bigg|\int_{U_s}
  \partial^{(\bsalpha_\setu;\bszero)}_\bsy \widehat u_{s,h}(\bsx,\bsy)\,\re^{-2\pi\ri\bsell\cdot\bsy}\,\rd\bsy
  \bigg|\,\rd\bsx \nonumber\\
&\leq \frac{1}{\gamma_{\setu}(2\pi)^{\alpha|\setu|}}
 \int_{U_s}\int_{\Dref} |\partial^{(\bsalpha_\setu;\bszero)}_\bsy\widehat u_{s,h}(\bsx,\bsy)|\,\rd\bsx\,\rd\bsy \nonumber\\
&\leq \frac{|{\Dref}|^{1/2}}{\gamma_{\setu}(2\pi)^{\alpha|\setu|}}
 \int_{U_s}\|\partial^{(\bsalpha_\setu;\bszero)}_\bsy\widehat u_{s,h}(\cdot,\bsy)\|_{L^2(\Dref)}\,\rd\bsy \nonumber\\
&\leq \frac{C_P|{\Dref}|^{1/2}}{\gamma_{\setu}(2\pi)^{\alpha|\setu|}}
 \int_{U_s}\|\partial^{(\bsalpha_\setu;\bszero)}_\bsy \widehat u_{s,h}(\cdot,\bsy)\|_{H_0^1(\Dref)}\,\rd\bsy,
\end{align}
where for the inequalities we used the unimodularity of the exponential
function, changed the order of integration, and applied both
Cauchy--Schwarz and Poincar\'e inequalities.

Combining \eqref{eq:cbc}, \eqref{eq:step1} and \eqref{eq:step2} together
with Theorem~\ref{thm:uhatregularity}, it follows that the CBC algorithm
can be used to find a generating vector satisfying the bound
\begin{align*}
{\rm error} &\leq 2\,C_{\rm init}\,C_P|\Dref|^{1/2}\,
 C_{\alpha,\bsgamma}(\lambda,s)^{1/\lambda}\,\bigg(\frac{1}{n-1}\bigg)^{1/\lambda}\\
&\quad\times \max_{\setu\subseteq\{1:s\}}\frac{1}{\gamma_{\setu}}\sum_{\bsm_{\setu}\in\{1:\alpha\}^{|\setu|}}\frac{(|\bsm_{\setu}|+d-1)!}{(d-1)!}\prod_{j\in\setu}\big(\widetilde C^{\alpha}m_j!\,\beta_j^{m_j}S(\alpha,m_j)\big),
\end{align*}
where $\widetilde
C:=\frac{2\,d!\,(2+\xi_{\bsb})^d(1+\xi_{\bsb})^3}{\sigma_{\min}^{d+4}}$,
$\beta_j=(2+\sqrt 2)\max(1+\sqrt 3,\rho)b_j$, and $C_{\rm init}$ is
defined as in~\eqref{eq:Cinit}. We choose the weights by
\eqref{eq:spodweights} to ensure that the factor on the second line is
bounded by $1$.

Finally, we show that with this choice of weights there is a constant
$C_{\alpha,\bsgamma}(\lambda)$ independent of $s$ such that
$C_{\alpha,\bsgamma}(\lambda,s)\leq C_{\alpha,\bsgamma}(\lambda)<\infty$.
To this end, we plug the weights~\eqref{eq:spodweights} into the
expression for $C_{\alpha,\bsgamma}(\lambda,s)$ in \eqref{eq:cbc}, define
$S_{\max}(\alpha):=\max_{k\in\{1:\alpha\}}S(\alpha,k)$, and use Jensen's
inequality (see, e.g.,~\cite[Theorem~19]{inequalities})$$\sum_k a_k\leq
\bigg(\sum_k a_k^\lambda\bigg)^{1/\lambda}\quad\text{for}~0<\lambda\leq
1~\text{and}~a_k\geq 0,$$ to find that
\begin{align}
&C_{\alpha,\bsgamma}(\lambda,s)\notag\\
&=\!\sum_{\varnothing\neq\setu\subseteq\{1:s\}} \!\!\!
\bigg(\sum_{\bsm_{\setu}\in\{1:\alpha\}^{|\setu|}} \!\!\!\!\frac{(|\bsm_{\setu}|\!+\!d\!-\!1)!}{(d\!-\!1)!}
\prod_{j\in\setu} \Big(\widetilde C^\alpha m_j!\beta_j^{m_j}S(\alpha,m_j)\Big)\bigg)^{\lambda}
(2\zeta(\alpha\lambda))^{|\setu|} \nonumber\\
&\leq\! \sum_{\varnothing\neq\setu\subseteq\{1:s\}}\sum_{\bsm_{\setu}\in\{1:\alpha\}^{|\setu|}}\!\!
\bigg(\frac{(|\bsm_{\setu}|\!+\!d\!-\!1)!}{(d\!-\!1)!} \prod_{j\in\setu}\!
 \Big(\widetilde C^{\alpha}m_j!\beta_j^{m_j}S(\alpha,m_j)(2\zeta(\alpha\lambda))^{\frac{1}{\lambda}}\Big)\!
 \bigg)^{\lambda} \nonumber\\
&\leq \sum_{\mathbf 0\neq \bsm\in\{1:\alpha\}^s}\bigg(\frac{(|\bsm|\!+\!d\!-\!1)!}{(d\!-\!1)!}\prod_{j=1}^s\mu_j^{m_j}\bigg)^{\lambda}, \label{eq:last}
\end{align}
where we define $\mu_j:=\max(1,\widetilde C^\alpha\,
\alpha!\,S_{\max}(\alpha)\,(2\zeta(\alpha\lambda))^{1/\lambda})\,\beta_j$.
By defining the auxiliary sequence $d_j:=\beta_{\lceil j/\alpha\rceil}$,
$j\geq 1$, that is
$$
d_{j\alpha+1}=d_{j\alpha+2}=\cdots=d_{(j+1)\alpha}=\beta_{j+1},\quad j\in\mathbb N_0,
$$
the last sum in \eqref{eq:last} is bounded from above by
\begin{align*}
 \sum_{\substack{\setv\subseteq \mathbb N\\ |\setv|<\infty}}
 \bigg(\frac{(|\setv|+d-1)!}{(d-1)!}\prod_{j\in\setv}d_j\bigg)^{\lambda}
&=\sum_{\ell=0}^\infty \bigg(\frac{(\ell+d-1)!}{(d-1)!}\bigg)^\lambda
 \sum_{\substack{\setv\subseteq\mathbb N\\ |\setv|=\ell}}\prod_{j\in\setv}d_j^\lambda \\
&\leq \sum_{\ell=0}^\infty \frac{1}{\ell!}\bigg(\frac{(\ell+d-1)!}{(d-1)!}\bigg)^\lambda
\bigg(\sum_{j=1}^\infty d_j^\lambda\bigg)^\ell,
\end{align*}
where the final inequality is due to the fact that $\big(\sum_{j=1}^\infty
d_j^\lambda\big)^\ell$ includes all products of the form $d_{j_1}
d_{j_2}\ldots  d_{j_\ell}$, not only those (as in the line above) with all
$j_i, \ldots, j_\ell$ different; and since the terms with all $j_1,
\ldots, j_\ell$ different are repeated $\ell !$ times,  in the last step
we can divide by $\ell!$.

By choosing $\lambda=p$ and $\alpha=\lfloor 1/p\rfloor+1$, it is easily
seen by Assumption~\ref{assumption5} and an application of the d'Alembert ratio test
that the upper bound in the above expression is finite. This implies that
the QMC error is of order $\mathcal O(n^{-1/p})$, with the implied
constant independent of the dimension $s$ when the weights are chosen
according to~\eqref{eq:spodweights}.
\quad\qed}

\subsection{Overall error estimate}
{Let $\widehat u$ be the solution to~\eqref{eq:weak2}, $\widehat u_s$ the
corresponding dimensionally-truncated solution, and $\widehat
u_{s,h}$ the dimensionally-truncated FE solution. It is easy to see that
\begin{align}
\begin{split}
&\bigg\|\int_U \widehat u(\cdot,\bsy)\,{\rm d}\bsy-\frac{1}{n}\sum_{k=1}^n \widehat u_{s,h}(\cdot,\bsy^{(k)})\bigg\|_{L^1(D_{\rm ref})}\\
&\leq \underset{\text{dimension truncation error}}{\underbrace{\bigg\|\int_U(\widehat u(\cdot,\bsy)-\widehat u_s(\cdot,\bsy))\,{\rm d}\bsy\bigg\|_{L^1(D_{\rm ref})}}}+\underset{\text{FE error}}{\underbrace{\bigg\|\int_{U_s}(\widehat u_s(\cdot,\bsy)-\widehat u_{s,h}(\cdot,\bsy))\,{\rm d}\bsy\bigg\|_{L^1(D_{\rm ref})}}}\\
&\quad +\underset{\text{QMC error}}{\underbrace{\bigg\|\int_{U_s}\widehat u_{s,h}(\cdot,\bsy)\,{\rm d}\bsy-\frac{1}{n}\sum_{k=1}^n\widehat u_{s,h}(\cdot,\bsy^{(k)})\bigg\|_{L^1(D_{\rm ref})}}}.\end{split}\label{eq:combinederror}
\end{align}}

We combine the error bounds derived in Subsections~\ref{sec:dim}--\ref{sec:qmc}. The error decomposition~\eqref{eq:combinederror} allows us to deduce the following overall error estimate.
\begin{theorem}
Let $(\bsy_k)_{k=1}^n$ be the lattice cubature nodes in $U_s$ generated by a CBC algorithm using the construction detailed in Subsection~\ref{sec:qmc}. Suppose that for each lattice point, we solve~\eqref{eq:weak2} using one common finite element discretization in the domain $D_{\rm ref}$. Under the assumptions {\rm \ref{assumption1}}--{\rm \ref{assumption7}}, we have the combined error estimate
$$
\bigg\|\int_U \widehat u(\cdot,\bsy)\,{\rm d}\bsy-\frac{1}{n}\sum_{k=1}^n \widehat u_{s,h}(\cdot,\bsy^{(k)})\bigg\|_{L^1(D_{\rm ref})}\leq C(s^{-2/p+1}+n^{-1/p}+h^2),
$$
where $h$ denotes the mesh size of the piecewise linear finite element mesh and $C>0$ is a constant independent of $s$, $h$, and $n$.
\end{theorem}

\section{Numerical experiments}\label{sec:5}
Let us consider the domain parameterization
\begin{align*}
D(\bsy):=\{(x_1,x_2)\in \mathbb R^2:0\leq x_1\leq 1,~0\leq x_2\leq a(x_1,\bsy)\},\quad \bsy=(y_j)_{j=1}^s\in U_s,
\end{align*}
where
\begin{align}
a(x,\bsy):=1+\frac{c}{\sqrt 6}\sum_{j=1}^s\sin(2\pi y_j)\psi_j(x),\quad x\in [0,1]~{\rm and}~\bsy\in U_s\label{eq:kl}
\end{align}
with $c>0$ and $\psi_j(x):=j^{-\theta}\cos(j\pi x)$, $\theta>1$.
Thus the reference domain in this case is the unit square, and the uncertain  boundary is confined to the upper edge of the square. It is possible to write $D(\bsy)=V(D(\mathbf 0),\bsy)$ with the vector-valued expansion
$$
\boldsymbol V(\bsx,\bsy):=\bsx+\frac{1}{\sqrt 6}\sum_{j=1}^s \sin(2\pi y_j)\bspsi_j(\bsx),
$$
where
$$
\bspsi_j(\bsx):=cj^{-\theta}\begin{bmatrix}0\\ x_2\cos(j\pi x_1)\end{bmatrix}.
$$
{We note that}
\begin{align}
\|\bspsi_j'(\bsx)\|_2=\sqrt{\cos^2(j\pi x_1)+j^2\pi^2x_2^2\sin^2(j\pi x_1)}\label{eq:help1}
\end{align}
and
\begin{align}
\sup_{\bsx\in D(\mathbf 0)}\sqrt{\cos^2(j\pi x_1)+j^2\pi^2x_2^2\sin^2(j\pi x_1)}=j\pi.\label{eq:help2}
\end{align}
The identity~\eqref{eq:help1} can be easily verified from the definition
of the matrix spectral norm, and~\eqref{eq:help2} follows by observing
that $\cos^2(j \pi x_1) + j^2 \pi^2 x_2^2\sin^2(j \pi x_1)  =  1 +
\sin^2(j \pi x_1)  (j^2 \pi^2 x_2^2  - 1)$ is maximized by
$x_1=\frac{1}{2j}$ and $x_2=1$. The elements of the sequence
$(b_j)_{j=1}^\infty$ needed for the construction of the QMC weights now
simplify to
\begin{align*}
&b_j=\frac{1}{\sqrt 6}\|\bspsi_j\|_{W^{1,\infty}(\Dref)}\\
&=\frac{1}{\sqrt 6}\!\max\!\bigg\{\underset{\bsx\in D(\mathbf 0)}{\rm ess\,sup}\|\bspsi_j(\bsx)\|_2,\underset{\bsx\in D(\mathbf 0)}{\rm ess\,sup}\|\bspsi_j'(\bsx)\|_2\bigg\}\\
&=\frac{c}{\sqrt 6}j^{-{\theta}}\!\max\!\bigg\{\underset{{\bsx\in D(\mathbf 0)}}{\rm sup}|x_2\cos(j\pi x_1)|,\!\underset{\bsx\in D(\mathbf 0)}{\rm sup}\!\sqrt{\cos^2(j\pi x_1)\!+\!j^2\pi^2x_2^2\sin^2(j\pi x_1)}\bigg\}\\
&=\frac{c}{\sqrt 6}j^{-\theta}\max\{1,j\pi\}=\frac{c\pi}{\sqrt 6}j^{1-\theta}.
\end{align*}
Moreover, we can take
$$
\xi_{\bsb}=\sum_{j=1}^\infty b_j=\frac{c\pi}{\sqrt 6}\zeta(\theta-1).
$$
The expected rate of QMC convergence is $\mathcal O(n^{-\theta+1})$.

\subsection{Source problem}\label{sec:source}

\begin{figure}[H]
\centering
  \hspace*{-.18cm}\subfloat{\label{fig:Xa}
  \begin{tikzpicture}
  \node (img) {\includegraphics[trim= 0cm 0cm 2.6cm 0cm,clip,height=.34\textwidth]{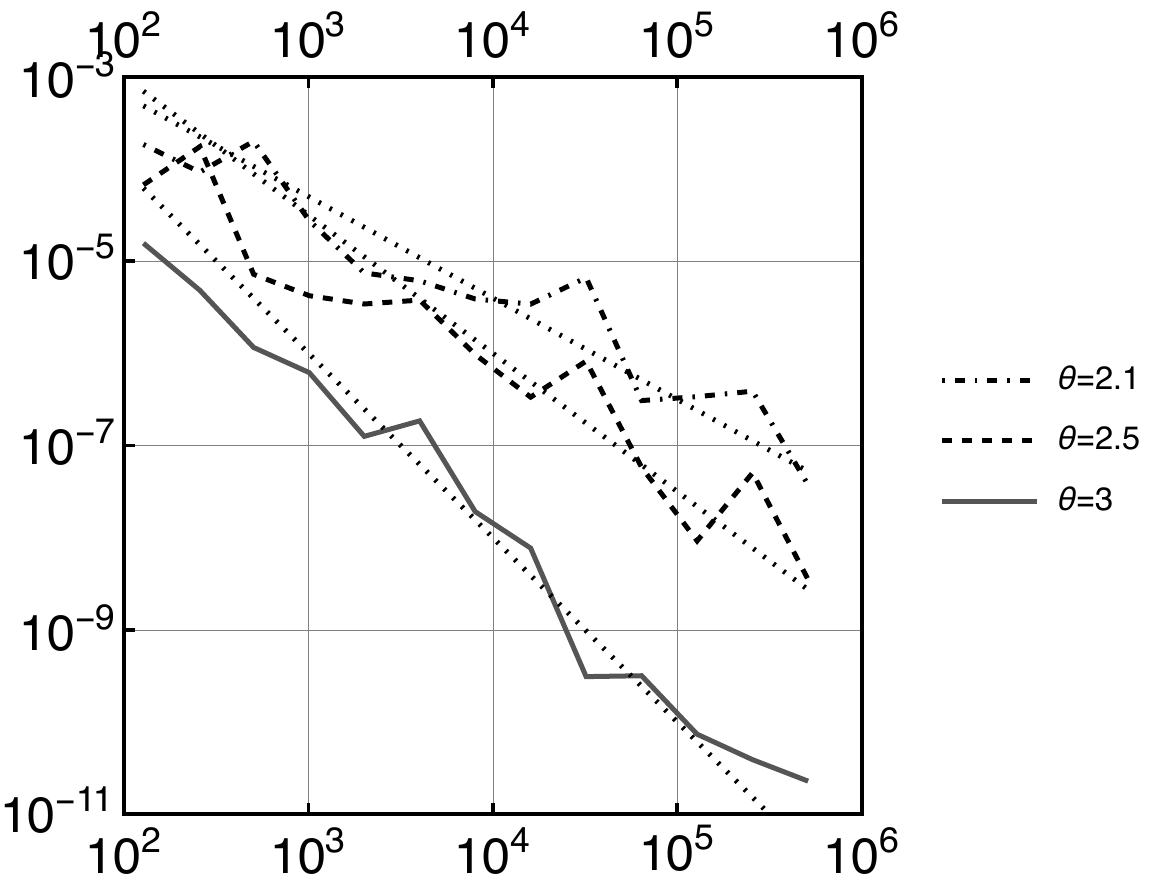}\!\!\!\!\!\!\includegraphics[trim= 9.2cm 2cm 0cm .4cm,clip,height=.4\textwidth]{Figures/source1.pdf}};
  \node[below=of img, node distance=0cm, xshift=-.62cm,yshift=1.1cm] {$n$};
  \node[left=of img, node distance=0cm, rotate=90, anchor=center,xshift=-.45cm,yshift=-0.95cm] {relative error};
    \node[above=of img, node distance=0cm, xshift=-.55cm,yshift=-1.85cm] {$\overset{\text{\footnotesize Source problem (full PDE solution):}}{\text{relative error vs.~$n$ with $\theta$}}$};
  \end{tikzpicture}
}%
  \!\!\!\!\!\!\!\!\subfloat{\label{fig:Xb}
  \begin{tikzpicture}
  \node (img) {\includegraphics[trim= 0cm 0cm 2.6cm 0cm,clip,height=.34\textwidth]{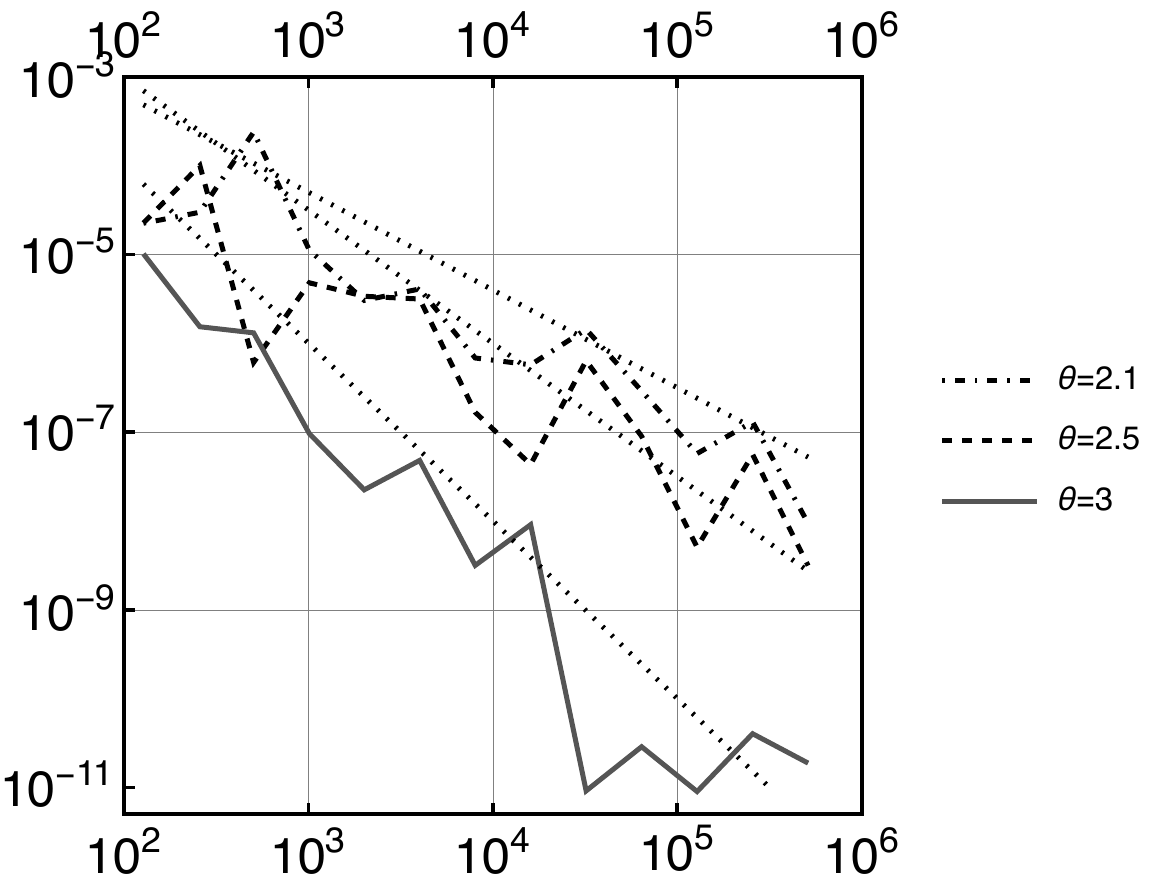}\!\!\!\!\!\!\includegraphics[trim= 9.2cm 2cm 0cm .4cm,clip,height=.4\textwidth]{Figures/source2.pdf}};
   \node[below=of img, node distance=0cm, xshift=-.62cm, yshift=1.1cm] {$n$};
  \node[left=of img, node distance=0cm, rotate=90, anchor=center,xshift=-.45cm,yshift=-0.95cm] {relative error};
      \node[above=of img, node distance=0cm, xshift=-.55cm,yshift=-1.85cm] {$\overset{\text{\footnotesize Source problem (nonlinear QoI):}}{\text{relative error vs.~$n$ with dimension $s$}}$};
  \end{tikzpicture}
}%
\caption{Left: Error criterion~\eqref{eq:numex1}. Right: Error criterion~\eqref{eq:numex2}. In both cases, the errors were computed for increasing $n$, different decay rates $\theta\in\{2.1,2.5, 3.0\}$, and fixed dimension $s=100$.}\label{fig:X}
\end{figure}

We consider the Poisson problem~\eqref{eq:model} equipped with
$f(\bsx):=x_2$ and the random domain described above for
$\theta\in\{2.1,2.5,3.0\}$ and $c=\sqrt{3/2}$.  The input random field is
truncated to $s=100$ terms. We construct a uniform regular triangulation
for the reference domain $D_{\rm ref}=D(\mathbf 0)$ with mesh size
$h=2^{-5}$, which is then mapped to each realization of the random domain
using~$\boldsymbol V$. The PDE is solved using a first order finite
element method.
For each realization of $\bsy\in U_s$, the finite element solution is
computed in $D(\bsy)$ and then mapped back onto the reference domain
$D_{\rm ref}$. This is equivalent to but simpler than computing the
solution on the reference domain, see \cite{nummat2016}.
We solve the PDE problem over QMC point sets with increasing number of
cubature nodes $n$ and compute the relative errors
\begin{align}
 \frac{\|\mathbb E[\widehat u]-Q_{s,n}(\widehat u)\|_{L^1(D_{\rm ref})}}{\|\mathbb E[\widehat u]\|_{L^2(D_{\rm ref})}},\label{eq:numex1}
\end{align}
where $Q_{s,n}$ denotes the $s$-dimensional QMC operator over $n$ cubature
nodes. As the reference solution, we use the approximate QMC solution
obtained using $n=1\,024\,207$ cubature nodes. The QMC point sets for
this experiment were constructed using the fast CBC algorithm with SPOD
weights of the form~\eqref{eq:spodweights},
where we set $\sigma_{\min}=\rho=1$ and $c=10^{-6}$ for numerical
stability. The results are displayed in the left-hand side image of
Figure~\ref{fig:X}. In this case, the observed rates $\sigma_{\theta}$ are
for $\theta = 2.1$: $\sigma_{2.1}=0.95$, $\theta = 2.5$: $\sigma_{2.5}=
1.15$, $\theta = 3$: $\sigma_{3} = 1.72$. The theoretically predicted
convergence rates are clearly observed in the numerical results: using a
higher value for the decay rate $\theta$ results in a faster rate of
convergence.

In addition, we also consider a nonlinear quantity of interest (QoI)
$$
G(\widehat u(\cdot,\bsy)):=\|\nabla \widehat u(\cdot,\bsy)\|_{L^2(D_{\rm ref})}
$$
and compute the relative errors
\begin{align}
\frac{|\mathbb E[G(\widehat u)]-Q_{s,n}(G(\widehat u))|}{|\mathbb E[G(\widehat u)]|}.\label{eq:numex2}
\end{align}
Again, the QMC solution obtained using $n=1\,024\,207$ cubature nodes is used as the reference. The results, computed using the same QMC rules as above, are displayed in the right-hand side image of Figure~\ref{fig:X}. The observed rates $\sigma_{\theta}$ are for $\theta = 2.1$: $\sigma_{2.1}=0.89$, $\theta = 2.5$: $\sigma_{2.5}= 1.09$, $\theta = 3$: $\sigma_{3} = 1.74$. Even though the use of a nonlinear QoI results in more irregular convergence behavior, the general trends are still visible, i.e., faster decay of the input random field yields a faster cubature convergence rate.

\subsection{Capacity problem}
In this section, we consider an example which is beyond the setting considered in the previous sections to demonstrate that our method also works well for a broader class of domain uncertainty quantification problems.

We will consider the Dirichlet--Neumann problem
\begin{align}
\begin{cases}
-\Delta u(\bsx,\bsy)=0\quad\text{for}~\bsx\in D(\bsy),~\bsy\in U_s,\\
\partial_{\rm n}u|_{x_1=0}=\partial_{\rm n}u|_{x_1=1}=0,\\
u|_{x_2=0}=0,~u|_{x_2=a(x_1,\bsy)}=1,
\end{cases}\label{eq:dn1}
\end{align}
and its conjugate problem
\begin{align}
\begin{cases}
-\Delta v(\bsx,\bsy)=0\quad\text{for}~\bsx\in D(\bsy),~\bsy\in U_s,\\
\partial_{\rm n}v|_{x_2=0}=\partial_{\rm n}v|_{x_2=a(x_1,\bsy)}=0,\\
v|_{x_1=0}=0,~v|_{x_1=1}=1.
\end{cases}\label{eq:dn2}
\end{align}
As the QoI, we consider the capacity
$$
{\rm cap}(D(\bsy)):=\int_{D(\bsy)}|\nabla u(\bsx,\bsy)|^2\,\rd\bsx,\quad \bsy\in U_s,
$$
where $u$ is the solution to~\eqref{eq:dn1}. We also define its conjugate as
$$
\overline{{\rm cap}}(D(\bsy)):=\int_{D(\bsy)}|\nabla v(\bsx,\bsy)|^2\,\rd\bsx,\quad \bsy\in U_s,
$$
where $v$ is the solution to~\eqref{eq:dn2}. It is well-known that (cf.,
e.g., \cite[Theorem~4-5]{ahlfors1973conformal})
\begin{align}
{\rm cap}(D(\bsy))\,\overline{{\rm cap}}(D(\bsy))=1\quad\text{for all}~\bsy\in U_s,\label{eq:reciprocal}
\end{align}
provided that $x\mapsto a(x,\bsy)$ specifies a Jordan arc such that
$a(x,\bsy)>0$ for all $x\in[0,1]$ and $\bsy\in U_s$.

The solutions to~\eqref{eq:dn1} and~\eqref{eq:dn2} are calculated
numerically using $hp$-FEM. As before, the finite element solutions
to~\eqref{eq:dn1} and~\eqref{eq:dn2} are {\em not} computed on the
reference domain $D_{\rm ref}$; for each realization of $\bsy\in U_s$,
the finite element solution is computed in $D(\bsy)$. Numerically, we can
use the reciprocal relation~\eqref{eq:reciprocal} to define an {\em a
posteriori} error estimate
\begin{align}
{\rm err}(\bsy):=|1-{\rm cap}(D(\bsy))\,\overline{{\rm cap}}(D(\bsy))|,\label{eq:aposteriori}
\end{align}
where ${\rm cap}(D(\bsy))$ and $\overline{{\rm cap}}(D(\bsy))$ are
computed using $hp$-FEM, and we tune the finite element approximation to ensure
that~\eqref{eq:aposteriori} is approximately of the same order across all
realizations of~$D(\bsy)$. One realization of $D(\bsy)$ with respective solutions $u(\bsx,\bsy)$
and $v(\bsx,\bsy)$ is shown in Figure~\ref{fig:0}.

\begin{figure}[!t]
  \centering
  \subfloat[{$u(\bsx,\bsy)$}]{\includegraphics[width=.3\textwidth]{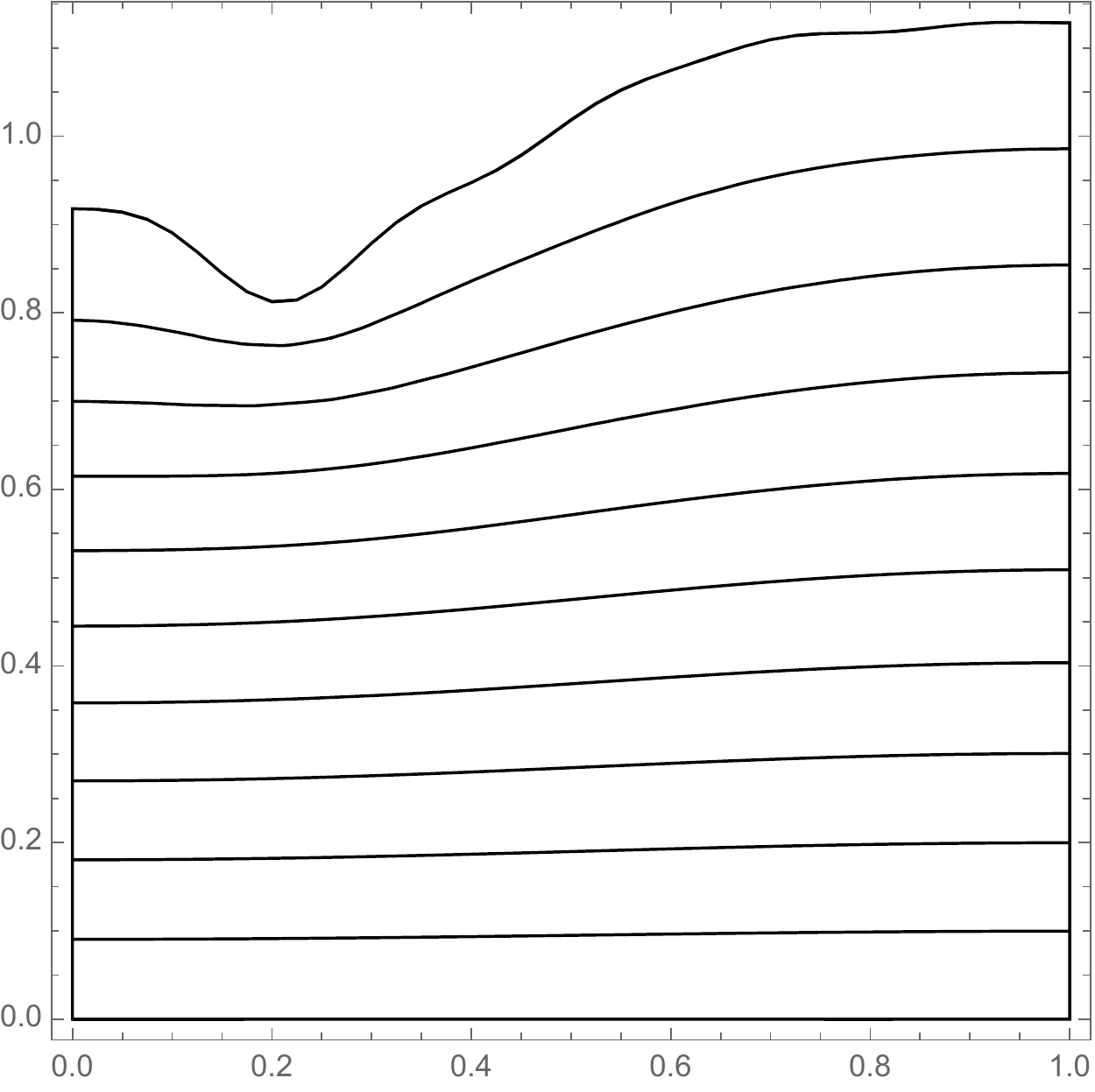}}%
  \quad
  \subfloat[{$v(\bsx,\bsy)$}]{\includegraphics[width=.3\textwidth]{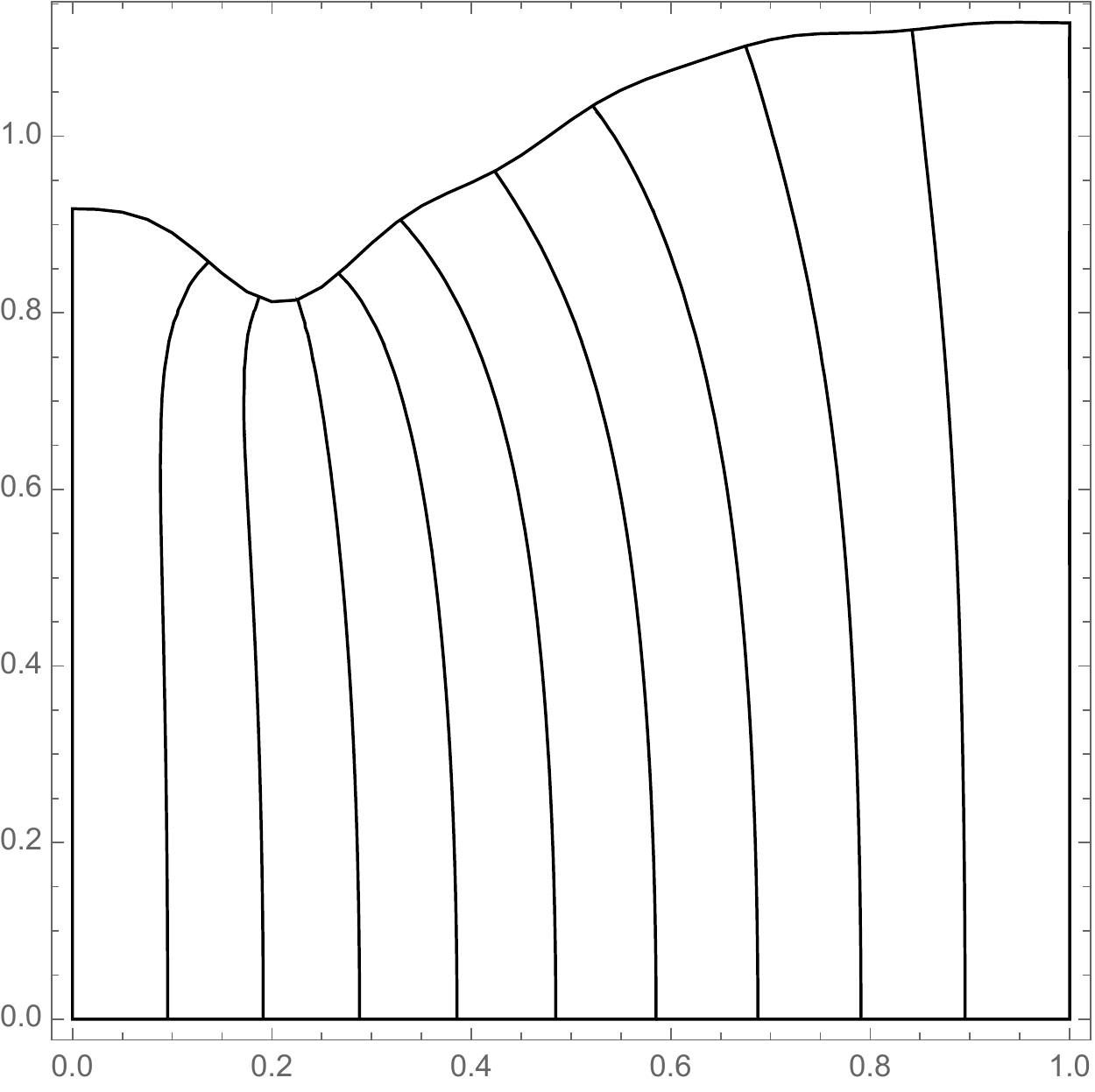}}%
  \quad
  \subfloat[Conformal map.]{\includegraphics[width=.3\textwidth]{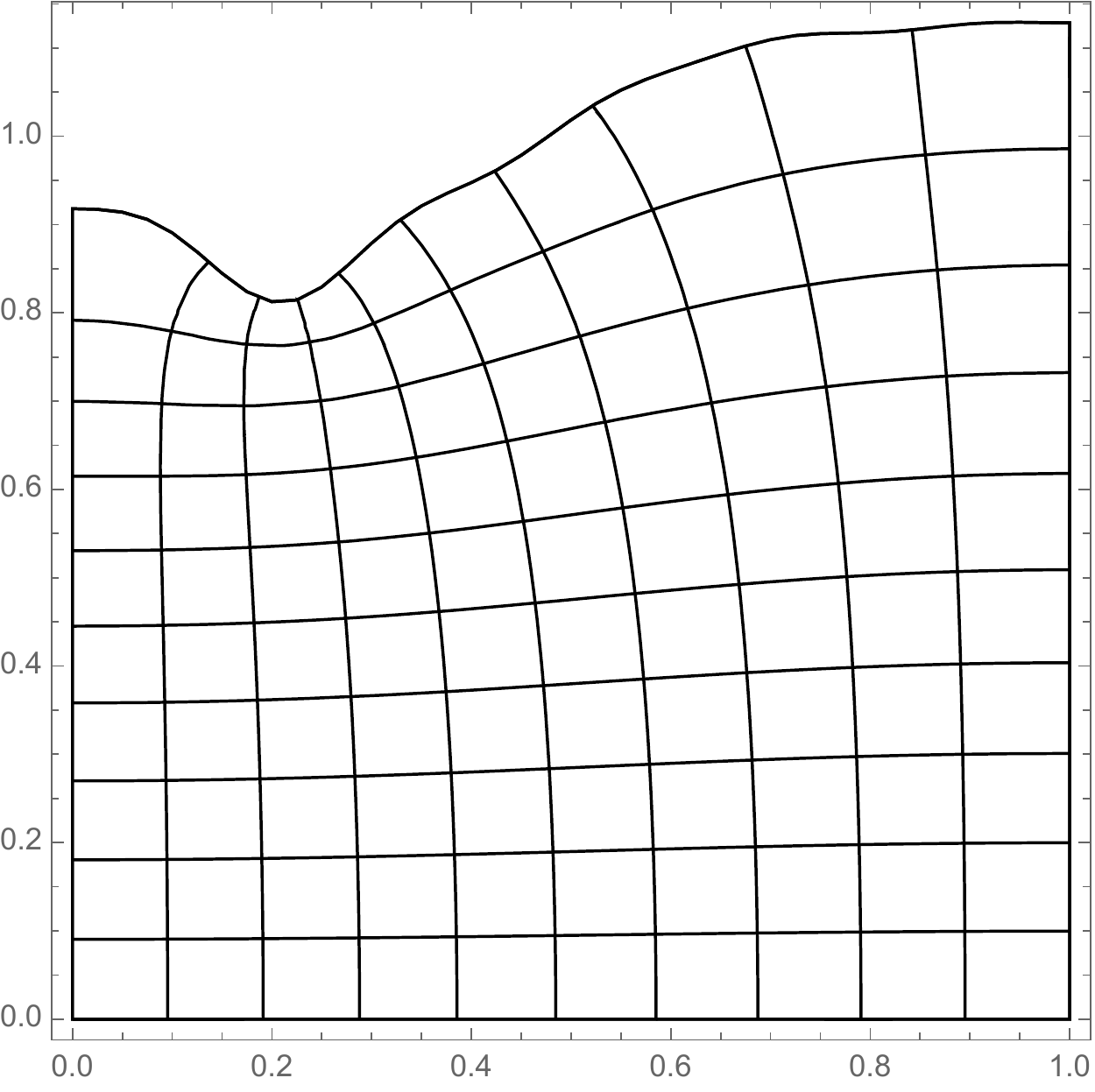}}
\caption{Sample realisation at $s=10$: Contour plots of $u(\bsx,\bsy)$ and $v(\bsx,\bsy)$.
Combination of the contour plots is the image of the conformal map of the domain.}\label{fig:0}
\end{figure}
\begin{figure}[!t]
\centering
  \hspace*{-.18cm}\subfloat[{Varying $\theta$.}]{\label{fig:1a}
  \begin{tikzpicture}
  \node (img) {\includegraphics[trim= .6cm .6cm 3.4cm .4cm,clip,height=.34\textwidth]{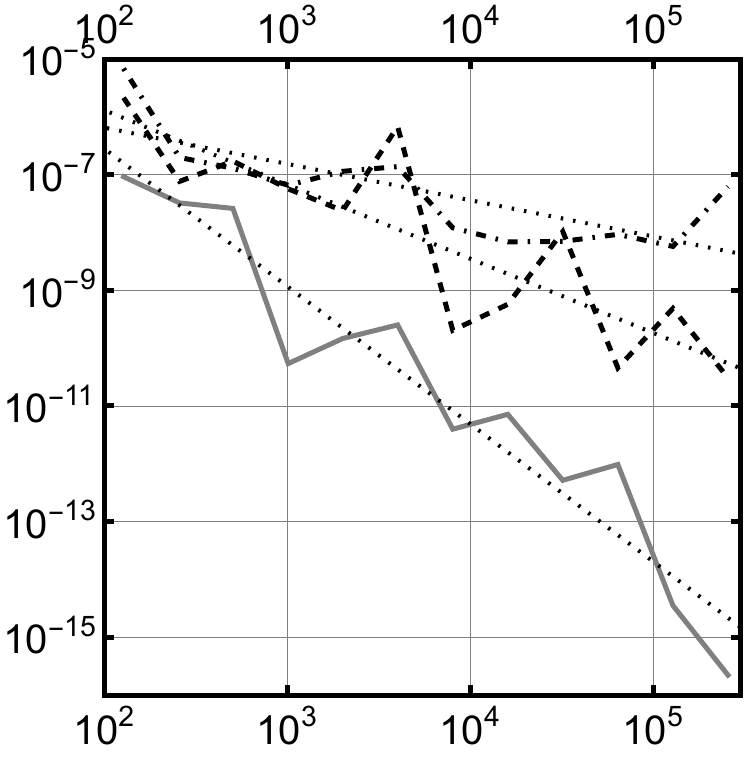}\includegraphics[trim= 9.2cm 2cm 0cm .4cm,clip,height=.4\textwidth]{Figures/fig9.pdf}};
  \node[below=of img, node distance=0cm, xshift=-.62cm,yshift=1.1cm] {$n$};
  \node[left=of img, node distance=0cm, rotate=90, anchor=center,xshift=-.45cm,yshift=-0.95cm] {relative error};
    \node[above=of img, node distance=0cm, xshift=-.65cm,yshift=-1.8cm] {$\overset{\text{\footnotesize Capacity:}}{\text{relative error vs.~$n$ with $\theta$}}$};
  \end{tikzpicture}
}%
  \!\!\!\!\!\!\!\!\subfloat[{Varying $s$.}]{\label{fig:1b}
  \begin{tikzpicture}
  \node (img) {\includegraphics[trim= .6cm .6cm 3.4cm .4cm,clip,height=.34\textwidth]{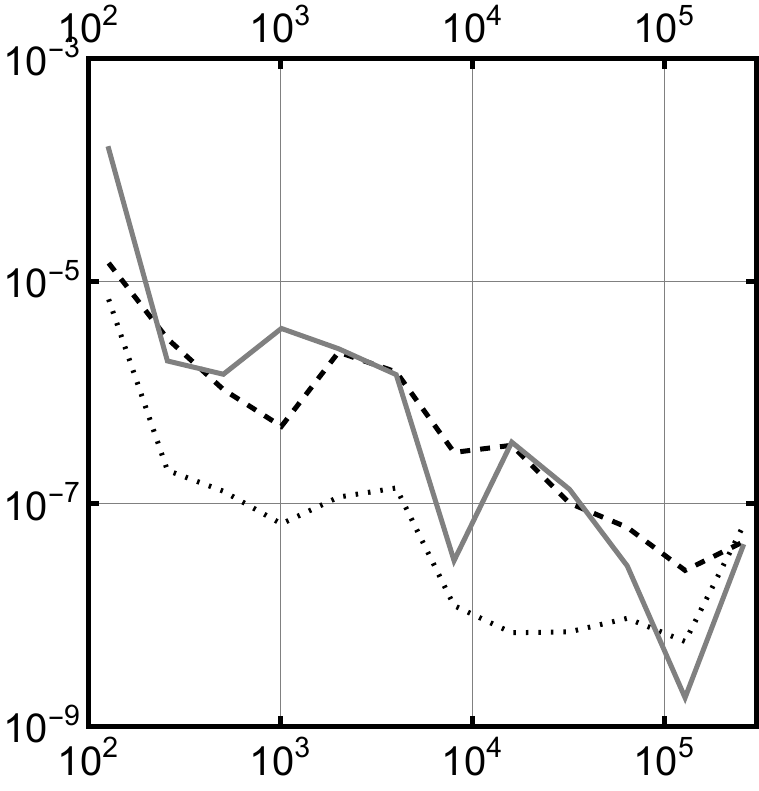}\includegraphics[trim= 9.2cm 2cm 0cm .4cm,clip,height=.4\textwidth]{Figures/fig5.pdf}};
   \node[below=of img, node distance=0cm, xshift=-.62cm, yshift=1.1cm] {$n$};
  \node[left=of img, node distance=0cm, rotate=90, anchor=center,xshift=-.45cm,yshift=-0.95cm] {relative error};
      \node[above=of img, node distance=0cm, xshift=-.65cm,yshift=-1.8cm] {$\overset{\text{\footnotesize Capacity:}}{\text{relative error vs.~$n$ with dimension $s$}}$};
  \end{tikzpicture}
}%
\\
\subfloat[{Focus on $\theta=2.5$.}]{\label{fig:1c}
\begin{tikzpicture}
\node (img) {\includegraphics[trim= .6cm .6cm 3.4cm .4cm,clip,height=.34\textwidth]{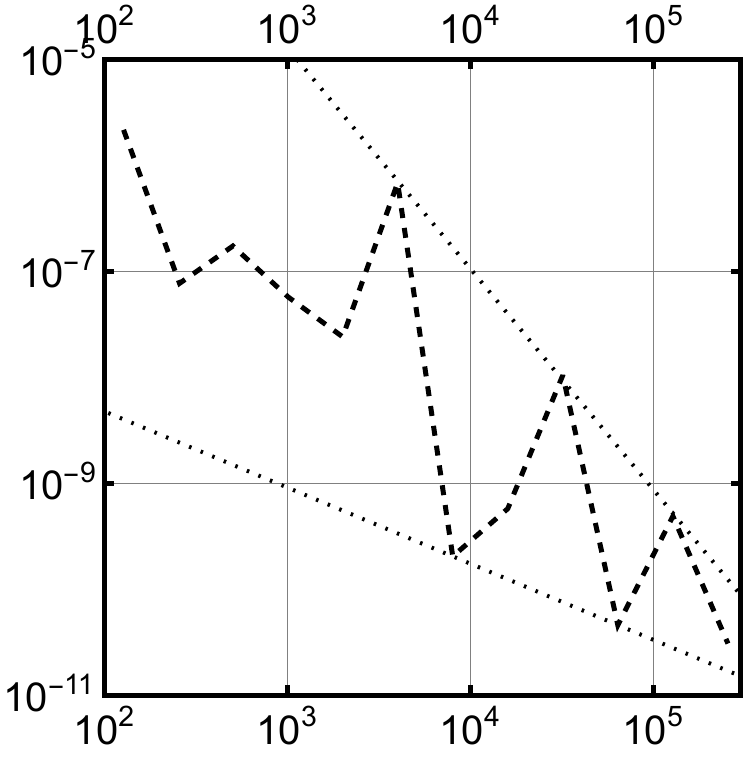}\includegraphics[trim= 9.2cm 2cm 0cm .4cm,clip,height=.4\textwidth]{Figures/fig10.pdf}};
  \node[below=of img, node distance=0cm, xshift=-.62cm, yshift=1cm] {$n$};
  \node[left=of img, node distance=0cm, rotate=90, anchor=center,xshift=-.45cm,yshift=-0.95cm] {relative error};
      \node[above=of img, node distance=0cm, xshift=-.65cm,yshift=-1.8cm] {$\overset{\text{\footnotesize Capacity:}}{\text{relative error vs.~$n$ with $\theta$}}$};
\end{tikzpicture}
}
\caption{Error criterion~\eqref{eq:relerr1} with increasing $n$. (a) Varying decay rate $\theta\in\{2.1,2.5, 3.0\}$, and fixed dimension $s=10$. (b) Varying dimensions $s\in\{10,40,100\}$, and fixed  decay rate $\theta=2.1$. (c) Detail of the case $s=10$, $\theta=2.5$.}\label{fig:1}
\end{figure}

\newpage

The numerical experiments were carried out as follows:
\begin{itemize}[align=right,leftmargin=2.4\parindent]
\item[(i)] We set $s=10$, $c=\sqrt{3/2}$, and let $\theta\in\{2.1,2.5,3.0\}$ in~\eqref{eq:kl}. We solve the problems~\eqref{eq:dn1}  and~\eqref{eq:dn2} over QMC point sets with increasing number of cubature nodes $n$ and consider the error criteria
\begin{equation}
{\rm Capacity}=\frac{|\mathbb E[{\rm cap}(D(\cdot))]-Q_{s,n}({\rm cap}(D(\cdot)))|}{|\mathbb E[{\rm cap}(D(\cdot))]|},\label{eq:relerr1}
\end{equation}
where $Q_{s,n}$ denotes the $s$-dimensional QMC operator over $n$ cubature points. As the reference solution, we use the approximate QMC solution computed using $n=1\,024\,207$.
(The reference result is computed using a tailored lattice. In two cases the general lattice reference is available and the effect on the convergence curves is negligible.) The results are displayed in Figure~\ref{fig:1a}.
\item[(ii)] We set $c=\sqrt{3/2}$, $\theta=2.1$, and let $s\in \{10,40,100\}$. We compute the relative errors~\eqref{eq:relerr1} for the problems~\eqref{eq:dn1} and~\eqref{eq:dn2} over QMC point sets with increasing number of cubature nodes $n$. Again, as the reference solution, we use the approximate QMC solution computed using $n=1\,024\,207$. The results are displayed in Figure~\ref{fig:1b}.
\end{itemize}The QMC point sets for the above experiments were computed using the fast CBC algorithm using the same specifications as in Subsection~\ref{sec:source}. {{In order to utilize the {\em a posteriori} estimate introduced above, we remark that the boundary conditions of the numerical experiments are not the homogeneous Dirichlet zero boundary conditions analyzed in the earlier theory.} %

The observed rates $\sigma_{\theta}$ are for $\theta = 2.1$: $\sigma_{2.1}=0.7$, $\theta = 2.5$: $\sigma_{2.5}= 1.3$, $\theta = 3$: $\sigma_{3} = 2.4$ (see Figure~\ref{fig:1}).
The rate does not appear to depend on the dimension $s$.
The cone of convergence for $\theta = 2.5$ shown in Figure~\ref{fig:1c}, is bounded by
the rates $\sigma_{2.5\text{A}}= 2.3$ (above) and $\sigma_{2.5\text{B}}= 0.6$ (below). The general convergence patterns are clearly observed. The irregularities in the convergence behavior may be explained by the nonlinear QoI (compare with the experiment in Section~\ref{sec:source}) as well as the more complicated boundary conditions: the gradient term in the expression for the capacity appears to be sensitive to obtuse angles at the Dirichlet--Neumann interface for certain realizations of the random domain.

\section{Conclusions}\label{sec:6}
We analyzed the Poisson problem subject to domain uncertainty, employing a domain mapping approach in which the shape uncertainty is transferred onto a fixed reference domain. As a result, the coefficient of the transported PDE problem on the reference domain becomes highly nonlinear with respect to the stochastic variables. To this end, we developed a novel parametric regularity analysis for the transported problem subject to such diffusion coefficients. The parametric regularity analysis is an important ingredient in the design of tailored QMC rules and we demonstrated that these methods exhibit higher order cubature convergence rates when the random perturbation field is parameterized using the model introduced in~\cite{periodicpaper} in which an infinite sequence of independent random variables enter the random field as periodic functions. The numerical results display that faster convergence in the stochastic fluctuations leads to faster QMC convergence in the computation of the PDE response, which is consistent with the theory.

Another novel feature of our work is the dimension truncation error analysis for the problem. Many studies in uncertainty quantification for PDEs with random coefficients exploit the parametric structure of the problem, using a Neumann series expansion to obtain dimension truncation error bounds. However, there have been several recent studies suggesting a Taylor series approach which is based on the parametric regularity of the problem~\cite{GGKSS2019,gk22,guth22}. The Taylor series approach can be used to obtain dimension truncation rates for non-affine parametric PDE problems, but a limitation of the aforementioned papers is that the parametric regularity bound needs to be of product-and-order dependent (POD) form in order for the Taylor-based approach to yield useful results. This issue is circumvented in the present paper by the use of a change of variables technique in order to transform our problem into a form in which the Taylor series approach can be used to produce improved dimension truncation error rates.

The results developed in this paper constitute important groundwork for other applications involving domain uncertainty such as domain shape recovery problems.

\section*{Acknowledgements}
F.~Y. Kuo and I.~H. Sloan acknowledge the support from the Australian Research Council (DP210100831). This research includes computations using the computational cluster Katana supported by Research Technology Services at UNSW Sydney~\cite{katana}.

\appendix \section{Technical results}\label{sec:tech}
\normalsize
\begin{lemma}\label{lemma:delannoyrecursion}
Let $\bsnu\in\mathscr F$, $c_0>0$, $c>0$, and $q\in\bbZ^+$. Let
$\bsbeta=(\beta_j)_{j\geq 1}$ and $(\Upsilon_{\bsnu})_{\bsnu\in\mathscr
F}$ be sequences of non-negative real numbers and suppose that
$$
 \Upsilon_{\bszero} \,\le\, c_0, \quad\mbox{and}\quad
 \Upsilon_{\bsnu} \,\le\, \sum_{\substack{\bszero\neq \bsm\leq \bsnu\\ |{\rm supp}(\bsm)|\leq q}}
 c^{|\bsm|}\binom{\bsnu}{\bsm}\Upsilon_{\bsnu-\bsm}
 \prod_{j\in{\rm supp}(\bsm)}\beta_j
$$
for all $\bsnu\in\mathscr F\setminus\{\bszero\}$.
Then
$$
\Upsilon_{\bsnu}\,\leq\, c_0\, c^{|\bsnu|}\sum_{\bsm\leq \bsnu}\bsm!\,D_q(\bsm)\,\bsbeta^{\bsm}\prod_{i\geq 1}S(\nu_i,m_i)
\qquad\text{for all}~\bsnu\in\mathscr F,
$$
where we define, with $\bse_\setu := \sum_{j\in\setu} \bse_j$,
\begin{align} \label{eq:delannoy}
 D_q(\bsm) \,:=\,\begin{cases}1&\text{if}~\bsm=\bszero,\\
 \displaystyle\sum_{\substack{\setu\subseteq{\rm supp}(\bsm)\\ 0 \ne |\setu|\leq q}}
 D_q(\bsm-\bse_\setu) & \text{if}~\bsm\in\mathscr F\setminus\{\bszero\}.
\end{cases}
\end{align}
The result is sharp in the sense that all inequalities can be replaced by
equalities. The result also extends to $q=\infty$, i.e., we can remove the
upper bounds on $|{\rm supp}(\bsm)|$ and $|\setu|$.
\end{lemma}

\proof We prove the claim by induction with respect to the order of the
multi-indices $\bsnu\in\mathscr F$. The claim is resolved immediately for
$\bsnu=\bszero$, so let us fix $\bsnu\in\mathscr F \setminus\{\bszero\}$
and suppose that the assertion holds for all multi-indices with order less
than $|\bsnu|$. Then
\begin{align*}
 &\Upsilon_{\bsnu}
 \leq \sum_{\substack{\bszero\neq \bsm\leq \bsnu\\ |{\rm supp}(\bsm)|\le q}} \!\!c^{|\bsm|}
 \binom{\bsnu}{\bsm}\bigg(\prod_{j\in{\rm supp}(\bsm)}\beta_j\bigg)
 \, c_0\,c^{|\bsnu|-|\bsm|}\\
 &\qquad\qquad\qquad\quad\times\sum_{\bsw\leq \bsnu-\bsm}\bsw!\, D_q(\bsw)\,\bsbeta^{\bsw}\prod_{i\geq 1}S(\nu_i-m_i,w_i)\\
 &= c_0c^{|\bsnu|}\!\!\!\sum_{\substack{\bszero\neq \bsm\leq \bsnu\\ |{\rm supp}(\bsm)|\le q}}\!
 \sum_{\bsw\leq \bsnu-\bsm}\!\binom{\bsnu}{\bsm}\!\bigg(\prod_{j\in{\rm supp}(\bsm)}\beta_j\bigg)
 \bsw! D_q(\bsw)\,\bsbeta^{\bsw}\prod_{i\geq 1}S(\nu_i-m_i,w_i)\\
 &= c_0c^{|\bsnu|}\sum_{\substack{\bsw\leq \bsnu\\ \bsw\neq \bsnu}}
 \bsw!\, D_q(\bsw)\,\bsbeta^{\bsw}
 \sum_{\substack{\bszero\neq \bsm\leq \bsnu-\bsw\\ |{\rm supp}(\bsm)|\le q}}\binom{\bsnu}{\bsm}\!
 \bigg(\prod_{j\in{\rm supp}(\bsm)}\beta_j\bigg)\!\prod_{i\geq 1}S(\nu_i-m_i,w_i),
\end{align*}
where we swapped the order of the sums over $\bsm$ and $\bsw$. The inner
sum over $\bsm$ can be rewritten as
\begin{align*}
 \sum_{\substack{\setu\subseteq{\rm supp}(\bsnu-\bsw)\\ 0\ne |\setu|\le q}}
 \prod_{j\in\setu} \bigg( \beta_j\sum_{m_j=1}^{\nu_j-w_j}\binom{\nu_j}{m_j}S(\nu_j-m_j,w_j)\bigg)\prod_{j\not\in\setu}S(\nu_j,w_j), 
\end{align*}
which\,can\,be\,further\,simplified\,using\,the\,identity\,(cf.,\,e.g.,\,\cite[formula~26.8.23]{dlmf})
$$
 \sum_{m_j=1}^{\nu_j-w_j}\binom{\nu_j}{m_j}S(\nu_j-m_j,w_j) = (w_j+1)\,S(\nu_j,w_j+1)
 \quad \text{for}~w_j<\nu_j.
$$
Hence
\begin{align*}
& \Upsilon_{\bsnu}\\
&\leq c_0 c^{|\bsnu|}\!\sum_{\substack{\bsw\leq \bsnu\\ \bsw\neq \bsnu}}\!\bsw!D_q(\bsw)\bsbeta^{\bsw}
    \!\!\sum_{\substack{\setu\subseteq{\rm supp}(\bsnu-\bsw)\\ 0\ne |\setu|\le q}}\!
     \prod_{j\in\setu}\! \bigg(\beta_j(w_j\!+\!1)S(\nu_j,w_j\!+\!1)\bigg)\!\prod_{j\not\in\setu}\!S(\nu_j,w_j) \\
&= c_0 c^{|\bsnu|}\!\sum_{\substack{\setu\subseteq{\rm supp}(\bsnu)\\ 0\ne |\setu|\le q}}\!
 \sum_{\bsw\leq \bsnu-\bse_\setu}\!\!
 \bsw! D_q(\bsw)\,\bsbeta^{\bsw}\! \prod_{j\in\setu}\! \bigg(\!\beta_j(w_j\!+\!1)S(\nu_j,w_j\!+\!1)\!\bigg)\!\prod_{j\not\in\setu}\!S(\nu_j,w_j)\\
&=c_0\, c^{|\bsnu|}\sum_{\substack{\setu\subseteq{\rm supp}(\bsnu)\\ 0\ne |\setu|\le q}}
 \sum_{\bse_\setu \leq \bsw'\leq \bsnu}
 \bsw'!\, D_q(\bsw'-\bse_\setu)\,\bsbeta^{\bsw'}\prod_{j\geq 1}S(\nu_j,w_j'),
\end{align*}
where we swapped the order of the sums over $\bsw$ and $\setu$ and then
carried out the change of variable $\bsw'= \bsw+\bse_\setu$. Changing the
order of the sums again leads to
\begin{align*}
 \Upsilon_{\bsnu}
&\leq c_0\, c^{|\bsnu|} \sum_{\bsw'\leq \bsnu} \bsw'!\,
 \bigg(\sum_{\substack{\setu\subseteq{\rm supp}(\bsw')\\ 0\ne |\setu|\le q}}
 D_q(\bsw'-\bse_\setu)\bigg) \bsbeta^{\bsw'} \prod_{j\geq 1}S(\nu_j,w_j').
\end{align*}
Relabeling $\bsw'$ to $\bsm$ and applying the definition of $D_q(\bsm)$
yields the desired result.

We remark that if the inequalities for $\Upsilon_0$ and the recursion on
$\Upsilon_\bsnu$ are replaced by equalities then all steps of the proof
are equalities and so is the final formula. We also remark that the proof
does not make any use of the upper bounds $|{\rm supp}(\bsm)|\le q$ and
$|\setu|\le q$, so the result extends trivially to the case $q=\infty$.
\quad\qed

Let $\bsm\in\mathscr F$. It is interesting to note that if $|{\rm
supp}(\bsm)| \le 2$, for example, $\bsm_{\setu}= k\bse_i + \ell \bse_j$
with $k,\ell\in\mathbb N$, then the numbers~\eqref{eq:delannoy} have the
following representation in closed form:
$$
 D_2(\bsm)=\sum_{t=0}^{\min\{k,\ell\}}2^t\binom{k}{t}\binom{\ell}{t}.
$$
These are so-called \emph{Delannoy numbers}. We prove an upper bound for
$D_2(\bsm)$ in general.

\begin{lemma}\label{lemma:delannoybound}
The numbers defined by~\eqref{eq:delannoy} with $q=2$ can be bounded by
\begin{align} \label{eq:delannoyalt}
 D_2(\bsm) \,\le\, a_{|\bsm|}' \,=\, |\bsm|!\,a_{|\bsm|} \quad\text{for all}~\bsm\in\mathscr F,
\end{align}
with the sequences
$$
 a_k'\,:=\, k!\,a_k,\qquad
 a_k \,:=\, \frac{(1+\sqrt 3)^{k+1}-(1-\sqrt 3)^{k+1}}{2^{k+1}\sqrt 3},\quad k\geq 0.
$$
\end{lemma}
\proof Let us first observe that the sequence $(a_k')_{k=0}^\infty$ can be
represented using the recurrence $a_0'=a_1'=1$ and
$a_k'=ka_{k-1}'+\frac{k^2-k}{2}a_{k-2}'$, $k\geq 2$; see, e.g., entry
A080599 in the \emph{On-line Encyclopedia of Integer Sequences
\textnormal{(}OEIS\textnormal{)}} and references therein. We also have
$a_0=a_1=1$ and $a_k = a_{k-1} + a_{k-2}/2$, $k\ge 2$.

We prove the claim \eqref{eq:delannoyalt} by induction with respect to the
order of the multi-indices $\bsm\in\mathscr F$ while utilizing the
recursive characterization of sequence $(a_k')_{k=0}^\infty$. We
immediately see that~\eqref{eq:delannoyalt} holds with equality when
$|\bsm|\leq 1$. Fix $\bsm\in\mathscr F$ such that $|\bsm|>1$ and suppose
that the claim holds for all multi-indices with order less than $|\bsm|$.
Then
\begin{align*}
 &D_2(\bsm)
 \,=\, \sum_{\substack{\setu\subseteq{\rm supp}(\bsm)\\ 0\ne |\setu|\le 2}} D_2(\bsm-\bse_\setu)\\
 &\,\le\, \sum_{\substack{\setu\subseteq{\rm supp}(\bsm)\\ 0\ne |\setu|\le 2}} a_{|\bsm|-|\setu|}'
  \,=\, a_{|\bsm|-1}' \binom{|{\rm supp}(\bsm)|}{1} + a_{|\bsm|-2}' \binom{|{\rm supp}(\bsm)|}{2}\\
 &\,\le\, a_{|\bsm|-1}'\binom{|\bsm|}{1} + a_{|\bsm|-2}' \binom{|\bsm|}{2}
 \,=\, |\bsm| a_{|\bsm|-1}' + \frac{|\bsm|^2\!-\!|\bsm|}{2} a_{|\bsm|-2}'
 \,=\, a_{|\bsm|}',
\end{align*}
as desired.\quad\qed

\begin{lemma} \label{lem:triple}
Let $\bsnu\in\mathscr{F}$ and let $(\bbA_\bsnu)_{\bsnu\in\mathscr{F}}$ and
$(\bbB_\bsnu)_{\bsnu\in\mathscr{F}}$ be arbitrary sequences of real
numbers. Then
\begin{align} \label{eq:triple}
  &\sum_{\bsm\leq \bsnu}\sum_{\bsw\leq \bsm}\sum_{\bsmu\leq \bsnu-\bsm}\binom{\bsnu}{\bsm}
  \bbA_\bsw\, \bbB_\bsmu \prod_{i\geq 1} \Big( S(m_i,w_i)\,S(\nu_i-m_i,\mu_i) \Big) \nonumber\\
  &\,=\, \sum_{\bsm\leq\bsnu} \bigg(\sum_{\bsw\leq\bsm} \binom{\bsm}{\bsw}\,\bbA_\bsw\, \bbB_{\bsm-\bsw}\bigg)
  \prod_{i\geq 1}S(\nu_i,m_i).
\end{align}
\end{lemma}

\proof
Moving the sum over $\bsm$ further inside, we can write
\begin{align*}
  &\mbox{LHS of~\eqref{eq:triple}}\\
  &= \sum_{\bsw\leq \bsnu} \sum_{\bsmu\leq \bsnu-\bsw} \bbA_\bsw\, \bbB_\bsmu
  \sum_{\bsw\leq \bsm\leq \bsnu}\binom{\bsnu}{\bsm} \prod_{i\geq 1} \Big( S(m_i,w_i)\,S(\nu_i-m_i,\mu_i) \Big)\\
  &= \sum_{\bsw\leq \bsnu} \sum_{\bsmu\leq \bsnu-\bsw} \bbA_\bsw\, \bbB_\bsmu
  \binom{\bsw+\bsmu}{\bsw}\prod_{i\geq 1}S(\nu_i,w_i+\mu_i) \\
  &= \sum_{\bsw\leq\bsnu} \sum_{\bsw\leq \bsmu'\leq \bsnu} \bbA_\bsw\, \bbB_{\bsmu'-\bsw} \binom{\bsmu'}{\bsw}\prod_{i\geq 1}S(\nu_i,\mu'_i) \\
  &= \sum_{\bsmu'\leq\bsnu} \bigg(\sum_{\bsw\leq\bsmu'} \binom{\bsmu'}{\bsw}\,\bbA_\bsw\, \bbB_{\bsmu'-\bsw}\bigg)
  \prod_{i\geq 1}S(\nu_i,\mu'_i)
  = \mbox{RHS of~\eqref{eq:triple}},
\end{align*}
as required. In the second equality we used the identity (cf.,
e.g.,~\cite[formula~26.8.23]{dlmf})
$$
\sum_{m=w}^\nu\binom{\nu}{m}S(m,w)S(\nu-m,\mu)=\binom{w+\mu}{w}S(\nu,w+\mu).
$$
In the third equality we introduced a change of variables $\bsmu'=
\bsw+\bsmu$. In the fourth equality we swapped the order of the sums.
Finally we relabel $\bsmu'$ to $\bsm$.
\quad\qed

\begin{lemma}\label{lemma:faadibrunorecursion}
Let $c>0$, $\bsbeta=(\beta_j)_{j\geq 1}$ and
$(\Upsilon_{\bsnu,\bslambda})_{\bsnu\in\mathscr F,\bslambda\in\mathbb
Z^d}$ be sequences of non-negative real numbers satisfying
\begin{align*}
 &\Upsilon_{\bsnu,\bszero} \,=\, \delta_{\bsnu,\bszero}, \qquad
  \Upsilon_{\bsnu,\bslambda} \,=\, 0 \quad\mbox{if}\quad |\bsnu|<|\bslambda| \mbox{ or } \bslambda\not\geq \boldsymbol 0,
  \qquad\mbox{and} \\
 &\Upsilon_{\bsnu+\bse_j,\bslambda}
 \,\le\, \beta_j\sum_{k=0}^{\nu_j} c^{k+1}\binom{\nu_j}{k}\sum_{\substack{\ell=1\\ \lambda_\ell>0}}^d\Upsilon_{\bsnu-k\bse_j,\bslambda-\bse_\ell}
 \quad\mbox{otherwise}.
\end{align*}
Then
$$
 \Upsilon_{\bsnu,\bslambda}\leq c^{|\bsnu|}\frac{|\bslambda|!}{\bslambda!}
 \sum_{\substack{\bsm\leq \bsnu \\ |\bsm|=|\bslambda|}}\bsbeta^{\bsm}\prod_{i\geq 1}S(\nu_i,m_i)\quad \text{for all}~\bsnu\in\mathscr F,~\bslambda\in\mathbb{N}_0^d.
$$
The result is sharp in the sense that all inequalities can be replaced by
equalities.
\end{lemma}

\proof We prove the claim by
induction with respect to the multi-indices $\bsnu\in\mathscr F$. The case
$\bsnu=\bszero$ is resolved immediately by observing that
$$
\Upsilon_{\bszero,\bszero}=1\quad \text{and}\quad \Upsilon_{\bszero,\bslambda}=0\quad \text{if}~\bslambda\neq\bszero
$$
as desired.

Fix $\bsnu\in\mathscr F$ and $\boldsymbol\lambda\in\mathbb N_0^d$, and
suppose that the claim holds for all pairs of multi-indices, with first
multi-index of order $\leq |\bsnu|$ and second multi-index $< \bslambda$.
If $\bslambda=\bszero$, then the claim is trivially true, so we may assume
in the following that $\bslambda\neq \bszero$. By letting $j\geq 1$ be
arbitrary, we find that
\begin{align*}
  &\Upsilon_{\bsnu+\bse_j,\bslambda}\\
  &\leq \beta_j\sum_{k=0}^{\nu_j}c^{k+1}\binom{\nu_j}{k}\!\sum_{\substack{\ell=1\\ \lambda_\ell>0}}^dc^{|\bsnu|-k}
  \frac{(|\bslambda|\!-\!1)!}{(\bslambda\!-\!\bse_\ell)!}\sum_{\substack{\bsm\leq \bsnu-k \bse_j \\|\bsm|=|\bslambda|\!-\!1}}
  \bsbeta^{\bsm}S(\nu_j\!-\!k,m_j)\prod_{i\neq j}S(\nu_i,m_i)\\
  &= c^{|\bsnu|+1}\frac{|\bslambda|!}{\bslambda!}\beta_j\sum_{k=0}^{\nu_j}\binom{\nu_j}{k}
  \sum_{\substack{\bsm\leq \bsnu-k \bse_j \\ |\bsm|=|\bslambda|-1}}\bsbeta^{\bsm}S(\nu_j-k,m_j)\prod_{i\neq j}S(\nu_i,m_i).
\end{align*}
We introduce $\bsm':=(m_1,\ldots,m_{j-1},m_{j+1},\ldots)$,
$\bsnu':=(\nu_1,\ldots,\nu_{j-1},\nu_{j+1},\ldots)$, and
$\boldsymbol{\beta'}:=(\beta_1,\ldots,\beta_{j-1},\beta_{j+1},\ldots)$.
Then
\begin{align*}
  &\Upsilon_{\bsnu+\bse_j,\bslambda}\\
  &\leq c^{|\bsnu|+1}\frac{|\bslambda|!}{\bslambda!}\beta_j\sum_{k=0}^{\nu_j}\binom{\nu_j}{k}
  \sum_{m_j=0}^{\nu_j-k}\beta_j^{m_j}S(\nu_j-k,m_j)\\
  &\qquad\qquad\qquad\qquad\quad\times\sum_{\substack{\bsm'\leq \bsnu'\\ |\bsm'|=|\bslambda|-m_j-1}}
  \boldsymbol{\beta'}^{\bsm'}\prod_{i\neq j}S(\nu_i,m_i')\\
  &= c^{|\bsnu|+1}\frac{|\bslambda|!}{\bslambda!}\sum_{m_j=0}^{\nu_j}\beta_j^{m_j+1}
  \bigg(\sum_{k=0}^{\nu_j-m_j}\binom{\nu_j}{k}S(\nu_j-k,m_j)\bigg)\\
  &\qquad\qquad\qquad\qquad\quad\times\sum_{\substack{\bsm'\leq \bsnu' \\|\bsm'|=|\bslambda|-m_j-1}}\boldsymbol{\beta'}^{\bsm'}\prod_{i\neq j}S(\nu_i,m_i').
\end{align*}
Making use of the fact that (see~\cite[formula 26.8.25]{dlmf})
$$
\sum_{k=0}^{\nu-m}\binom{\nu}{k}S(\nu-k,m)=S(\nu+1,m+1)\quad\text{for all}~\nu\in\mathbb{N}_0,~m\in\{0,\ldots,\nu\},
$$
we obtain
\begin{align*}
  &\Upsilon_{\bsnu+\bse_j,\bslambda}\\
  &\leq c^{|\bsnu|+1}\frac{|\bslambda|!}{\bslambda!}\sum_{m_j=0}^{\nu_j}\beta_j^{m_j+1}S(\nu_j+1,m_j+1)
  \sum_{\substack{\bsm'\leq \bsnu' \\|\bsm'|=|\bslambda|-m_j-1}} \boldsymbol{\beta'}^{\bsm'}\prod_{i\neq j}S(\nu_i,m_i')\\
  &= c^{|\bsnu|+1}\frac{|\bslambda|!}{\bslambda!}
  \sum_{\overline{m}_j=0}^{\nu_j+1}\beta_j^{\overline{m}_j}S(\nu_j+1,\overline{m}_j)
  \sum_{\substack{\bsm'\leq \bsnu' \\ |\bsm'|=|\bslambda|-\overline{m}_j}}\boldsymbol{\beta'}^{\bsm'}
  \prod_{i\neq j}S(\nu_i,m_i'),
\end{align*}
where we used the change of variable $\overline{m}_j = m_j+1$ in
conjunction with the property $S(\nu+1,0)=0$ for all $\nu\in\mathbb N_0$.
Recognizing the above expression as
$$
 c^{|\bsnu|+1}\frac{|\bslambda|!}{\bslambda!}
 \sum_{\substack{\bsm\leq \bsnu \\ |\bsm|=|\bslambda|}}\bsbeta^{\bsm}S(\nu_j+1,m_j)\prod_{i\neq j}S(\nu_i,m_i)
$$
concludes the proof. We remark that if the inequality for the recursion on
$\Upsilon_{\bsnu,\lambda}$ is replaced by equality then all steps of the
proof are equalities and so is the final formula.
\quad\qed

\begin{lemma}\label{lemma:nightmarelemma}
Let $C>0$ be a constant and suppose that $\bsbeta=(\beta_j)_{j\geq 1}$ is
a sequence of non-negative real numbers. Suppose that
$\Upsilon_{\bszero}\leq c_0$ and
\begin{align*}
 \Upsilon_{\bsnu}
 &\le \sum_{\satop{\bsm\leq \bsnu}{\bsm\ne\bsnu}}\binom{\bsnu}{\bsm}\Upsilon_{\bsm}\, {c^{|\bsnu-\bsm|}}
 \sum_{\bsmu\leq \bsnu-\bsm} {(|\bsmu|+1)!}\,{\bsmu!}\,\bsbeta^{\bsmu}\prod_{i\geq 1}S(\nu_i-m_i,\mu_i) \\
 &\qquad + {c_0\,c^{|\bsnu|}}\sum_{\boldsymbol m\leq \bsnu} {\frac{(|\bsm|+d)!}{d!}}\,\bsbeta^{\bsm}\prod_{i\geq 1}S(\nu_i,m_i)
 \qquad \text{for all } \bsnu\in \mathscr F.
\end{align*}
Then
\begin{align*}
 \Upsilon_{\bsnu}\le {c_0\,c^{|\bsnu|}} \sum_{\bsm\leq\bsnu}
 \bbP_\bsm\,{\bsm!}\,\bsbeta^{\bsm}\prod_{i\geq 1}S(\nu_i,m_i)
 \qquad \text{for all } \bsnu\in \mathscr F,
\end{align*}
where
\begin{align} \label{eq:defP}
  &\bbP_{\bszero}=1 \quad\text{and}\quad
   \bbP_{\bsm}= {\frac{(|\bsm|+d)!}{{\bsm!}\,d!}}+\sum_{\substack{\bsw\le\bsm\\ \bsw\ne\bsm}}
   \bbP_{\bsw}\,{(|\bsm|-|\bsw|+1)!}\quad \text{for}~\bsm\in\mathscr F.
\end{align}
The result is sharp in the sense that all inequalities can be replaced by
equalities.
\end{lemma}

\proof The proof is carried out by induction with respect to the
multi-indices $\bsnu\in\mathscr{F}$. The base step $\bsnu=\bszero$ is
trivially true, so let us assume that $\bsnu\in\mathscr
F\setminus\{\bszero\}$ is such that the claim holds for all multi-indices
with order $<|\bsnu|$. Then the recurrence together with the induction
hypothesis yields that
\begin{align} \label{eq:split}
 \Upsilon_{\bsnu}
 &\leq R_\bsnu + {c_0\,c^{|\bsnu|}} \sum_{\bsm\leq \bsnu}{\frac{(|\bsm|+d)!}{d!}}\,\bsbeta^{\bsm}\prod_{i\geq 1}S(\nu_i,m_i),
\end{align}
with
\begin{align*}
 R_\bsnu
&:=\sum_{\satop{\bsm\leq\bsnu}{\bsm\ne\bsnu}}\binom{\bsnu}{\bsm}
 \bigg({c_0\,c^{|\bsm|}} \sum_{\bsw\leq \bsm}%
 \bbP_{\bsw}%
 \,{\bsw!}\,\bsbeta^{\bsw}\prod_{i\geq 1}S(m_i,w_i)\bigg)\\
&\quad\quad \times {c^{|\bsnu-\bsm|}}\sum_{\bsmu\leq \bsnu-\bsm} {(|\bsmu|+1)!}\,{\bsmu!}\,\bsbeta^{\bsmu}
 \prod_{i\geq 1}S(\nu_i-m_i,\mu_i).
\end{align*}

To simplify our presentation we abbreviate temporarily $\bbA_\bsw :=
\bbP_{\bsw}\,{\bsw!}\,\bsbeta^{\bsw}$ and $\bbB_\bsmu :=
{(|\bsmu|+1)!}\,{\bsmu!}\,\bsbeta^{\bsmu}$. Then we can write
\begin{align} \label{eq:R1}
 R_\bsnu
 &= {c_0\,c^{|\bsnu|}} \sum_{\substack{\bsm\leq\bsnu \\ \bsm\neq\bsnu}}\binom{\bsnu}{\bsm}
 \bigg(\sum_{\bsw\leq\bsm} \bbA_\bsw \prod_{i\geq 1}S(m_i,w_i)\bigg)
 \sum_{\bsmu\leq \bsnu-\bsm} \bbB_\bsmu \prod_{i\geq 1}S(\nu_i-m_i,\mu_i) \nonumber\\
&= {c_0\,c^{|\bsnu|}}\sum_{\bsm\leq\bsnu}
 \sum_{\bsw\leq\bsm} \sum_{\bsmu\leq \bsnu-\bsm} \binom{\bsnu}{\bsm}
 \bbA_\bsw\,\bbB_\bsmu \prod_{i\geq 1} \Big( S(m_i,w_i)\,S(\nu_i-m_i,\mu_i) \Big) \nonumber\\
&\quad
- {c_0\,c^{|\bsnu|}}\sum_{\bsw\leq \bsnu} \bbA_\bsw\,\bbB_\bszero \prod_{i\geq 1} S(\nu_i,w_i) \nonumber \\
&= {c_0\,c^{|\bsnu|}}\sum_{\bsm\leq\bsnu}
 \bigg(\sum_{\bsw\leq\bsm} \binom{\bsm}{\bsw} \bbA_\bsw\,\bbB_{\bsm-\bsw}\bigg)
 \prod_{i\geq 1} S(\nu_i,m_i)\nonumber
 \\
 &\quad - {c_0\,c^{|\bsnu|}}\sum_{\bsm\leq \bsnu} \bbA_\bsm\,\bbB_\bszero \prod_{i\geq 1} S(\nu_i,m_i) \nonumber\\
&= {c_0\,c^{|\bsnu|}} \sum_{\bsm\leq\bsnu}
 \bigg(\sum_{\satop{\bsw\leq\bsm}{\bsw\ne\bsm}} \binom{\bsm}{\bsw} \bbA_\bsw\,\bbB_{\bsm-\bsw} \bigg)
 \prod_{i\geq 1} S(\nu_i,m_i),
\end{align}
where we used Lemma~\ref{lem:triple} in the third equality.

Reinserting the definitions of $\bbA_\bsw$ and~$\bbB_{\bsm-\bsw}$, the
inner sum over $\bsw$ in \eqref{eq:R1} becomes
\begin{align} \label{eq:R2}
 \sum_{\satop{\bsw\leq\bsm}{\bsw\ne\bsm}} \binom{\bsm}{\bsw} \bbA_\bsw\,\bbB_{\bsm-\bsw}
 &= \sum_{\satop{\bsw\leq\bsm}{\bsw\ne\bsm}} \binom{\bsm}{\bsw}
 \bbP_{\bsw}\,{\bsw!}\,\bsbeta^{\bsw}\,
 {(|\bsm-\bsw|+1)!}\,{(\bsm-\bsw)!}\,\bsbeta^{\bsm-\bsw} \nonumber\\
 &= {\bsm!}\,\bsbeta^\bsm \sum_{\satop{\bsw\leq\bsm}{\bsw\ne\bsm}} \bbP_{\bsw}\, {(|\bsm-\bsw|+1)!}.
\end{align}
Substituting \eqref{eq:R1} and \eqref{eq:R2} into \eqref{eq:split}, we
obtain
\begin{align*}
 \Upsilon_\bsnu
&\le {c_0\,c^{|\bsnu|}} \sum_{\bsm\leq\bsnu}
\bigg( {\frac{(|\bsm|+d)!}{{\bsm!}\,d!}} +
 \sum_{\satop{\bsw\leq\bsm}{\bsw\ne\bsm}} \bbP_{\bsw}\,{(|\bsm-\bsw|+1)!} \bigg)
 {\bsm!}\,\bsbeta^\bsm
 \prod_{i\geq 1} S(\nu_i,m_i).
\end{align*}
The result now follows from the definition of the sequence $\bbP_\bsm$.
\quad\qed

We need to derive an upper bound for the sequences
$(\bbP_{\bsm})_{\bsm\in\mathscr F}$. In anticipation of this, let us
define the auxiliary sequence
\begin{align} \label{eq:tau}
 \tau_0=1\quad\text{and}\quad \tau_k=\sum_{j=0}^{k-1} (k-j+1)\,\tau_j,\quad k\geq 1.
\end{align}
This sequence is listed as A003480 in \emph{The On-Line Encyclopedia of Integer Sequences
\textnormal{(}OEIS\textnormal{)}}. We observe that this sequence satisfies
the following property.

\begin{lemma} \label{lemma:renewal}
The sequence $(\tau_k)_{k\geq 0}$ defined by \eqref{eq:tau} satisfies
$$
\sum_{\satop{\bsw\le\bsm}{\bsw\ne\bsm}}\tau_{|\bsw|}\,|\bsw|!\,(|\bsm|-|\bsw|+1)!\,\binom{\bsm}{\bsw}
 = \tau_{|\bsm|}\,|\bsm|!
 \qquad\text{for }\bsm\in\mathscr F\setminus\{\mathbf 0\}.
$$
\end{lemma}

\proof Let $\bsm\in\mathscr F\setminus\{\mathbf 0\}$. Then
\begin{align*}
&\sum_{\satop{\bsw\le\bsm}{\bsw\ne\bsm}}\tau_{|\bsw|}\,|\bsw|!\,(|\bsm|-|\bsw|+1)! \binom{\bsm}{\bsw}
=\sum_{j=0}^{|\bsm|-1}\tau_j\,j!\,(|\bsm|-j+1)!
  \sum_{\substack{\bsw\leq\bsm\\ |\bsw|=j}}\binom{\bsm}{\bsw}\\
&=\sum_{j=0}^{|\bsm|-1}\tau_j\,j!\,(|\bsm|-j+1)!\,\binom{|\bsm|}{j}
=|\bsm|!\sum_{j=0}^{|\bsm|-1}\tau_j\,(|\bsm|-j+1) = \tau_{|\bsm|}\,|\bsm|!,
\end{align*}
as required.
\quad\qed

{\em Remark.} Lemma~\ref{lemma:renewal} may be viewed as an extension of the {\em renewal property} of sequence A003480; see corresponding entry in the OEIS for further details.

We are now ready to state an upper bound for the sequence
$(\bbP_\bsm)_{\bsm\in\mathscr F}$.

\begin{lemma} \label{lem:P-bound}
The sequence $(\bbP_{\bsm})_{\bsm\in\mathscr F}$ defined
by~\eqref{eq:defP} satisfies
$$
 \bbP_{\bsm}\le \frac{2\,\tau_{|\bsm|}\,(|\bsm|+d-1)!}{(d-1)!}
 \qquad \text{for all}~\bsm\in\mathscr F\setminus\{\bszero\},
$$
where the sequence $(\tau_k)_{k\ge 0}$ is defined by~\eqref{eq:tau}.
\end{lemma}

\proof The proof is carried out by induction with respect to the
multi-indices $\bsm\in\mathscr F\setminus\{\mathbf 0\}$.  The base step
$\bsm=\boldsymbol e_j$, $j\geq 1$, is resolved by observing that $\tau_1 =
2$ and
$$
 \bbP_{\boldsymbol e_j} = d+3 \,\le\, 4\,d \,=\, \frac{2\,\tau_{|\bse_j|}\,(|\bse_j|+d-1)!}{(d-1)!}\,.
$$

Let $\bsm\in\mathscr F\setminus\{\mathbf 0\}$ with $|\bsm|\geq 2$ and
suppose that the assertion holds for all multi-indices $\bsw\in\mathscr
F\setminus\{\mathbf 0\}$ with $|\bsw|<|\bsm|$. Then by separating out the
$\bsw=\bszero$ term in \eqref{eq:defP} and applying the induction
hypothesis to $\bbP_\bsw$ for $\bsw\ne\bszero$, we obtain
\begin{align*}
 \bbP_{\bsm}
 &\leq \frac{(|\bsm|+d)!}{\bsm!\,d!}+(|\bsm|+1)!+\sum_{\substack{\bszero\ne\bsw\leq \bsm\\\bsw\neq \bsm}}
       \frac{2\,\tau_{|\bsw|}\,(|\bsw|+d-1)!}{(d-1)!}\,(|\bsm|-|\bsw|+1)!\\
 &\leq \frac{(|\bsm|+d)!}{d!}+\frac{(|\bsm|+d)!}{d!} \\
 &\quad + \frac{2\,(|\bsm|+d-1)!}{|\bsm|!\,(d-1)!}
       \sum_{\substack{\bszero\ne\bsw\leq \bsm\\ \bsw\neq \bsm}}\,\tau_{|\bsw|}\,|\bsw|!\,(|\bsm|-|\bsw|+1)!\,
       \binom{\bsm}{\bsw}\\
 &= \frac{2\,(|\bsm|+d)!}{d!} + \frac{2\,(|\bsm|+d-1)!}{|\bsm|!\,(d-1)!}
       \Big(\tau_{|\bsm|}\,|\bsm|! - (|\bsm|+1)!\Big) \\
 &= \frac{2\,(|\bsm|+d-1)!}{(d-1)!}\, \Big(\frac{|\bsm|+d}{d} + \tau_{|\bsm|} - (|\bsm|+1) \Big) \\
 &\le \frac{2\,\tau_{|\bsm|}\,(|\bsm|+d-1)!}{(d-1)!},
\end{align*}
where we used Lemma~\ref{lemma:renewal} in the first equality. This
completes the proof.
\quad\qed

It remains to obtain an estimate on the $(\tau_k)_{k\geq 0}$ sequence.

\begin{lemma}\label{lemma:tau-bound} The sequence defined by~\eqref{eq:tau} satisfies
$$
  \tau_k\le (2+\sqrt 2)^k\quad \text{for all $k\ge 0$.}
$$
\end{lemma}

\proof The claim is clearly true with equality when $k=0$. If
$k\geq 1$, then the sequence admits the following characterization (see
the OEIS entry for this sequence and references therein):
$$
 \tau_k=\frac{(2+\sqrt{2})^{k+1}}{4\sqrt{2}}-\frac{(2-\sqrt{2})^{k+1}}{4\sqrt{2}},
$$
from which the claim readily follows.\quad\qed

\bibliographystyle{plain}
\bibliography{domainuq}

\begin{thebibliography}{10}

\bibitem{ahlfors1973conformal}
L.~V. Ahlfors.
\newblock {\em Conformal Invariants: Topics in Geometric Function Theory}.
\newblock Higher Mathematics Series. McGraw-Hill, 1973.

\bibitem{BC02}
I.~Babu\v{s}ka and P.~Chatzipantelidis.
\newblock On solving elliptic stochastic partial differential equations.
\newblock {\em Comput. Methods Appl. Mech. Engrg.}, 191:4093--4122, 2002.

\bibitem{CNT16}
J.~E. Castrill\'on-Cand\'as, F.~Nobile, and R.~F. Tempone.
\newblock Analytic regularity and collocation approximation for elliptic {PDEs}
  with random domain deformations.
\newblock {\em Comput. Math. Appl.}, 71(6):1173--1197, 2016.

\bibitem{Navier}
A.~Cohen, Ch. Schwab, and J.~Zech.
\newblock Shape holomorphy of the stationary {N}avier--{S}tokes equations.
\newblock {\em SIAM J. Math. Anal.}, 50(2):1720--1752, 2018.

\bibitem{spodpaper14}
J.~Dick, F.~Y. Kuo, Q.~T. {Le~Gia}, D.~Nuyens, and Ch. Schwab.
\newblock Higher order {QMC} {P}etrov--{G}alerkin discretization for affine
  parametric operator equations with random field inputs.
\newblock {\em SIAM J. Numer. Anal.}, 52(6):2676--2702, 2014.

\bibitem{korobovpaper}
J.~Dick, I.~H. Sloan, X.~Wang, and H.~Wo\'{z}niakowski.
\newblock Good lattice rules in weighted {K}orobov spaces with general weights.
\newblock {\em Numer. Math.}, 103(1):63--97, 2006.

\bibitem{gantner}
R.~N. Gantner.
\newblock Dimension truncation in {QMC} for affine-parametric operator
  equations.
\newblock In A.~B. Owen and P.~W. Glynn, editors, {\em Monte Carlo and
  Quasi-Monte Carlo Methods 2016}, pages 249--264, Stanford, CA, 2018.
  Springer.

\bibitem{Gilbarg}
D.~Gilbarg and N.~S. Trudinger.
\newblock {\em {E}lliptic {P}artial {D}ifferential {E}quations of {S}econd
  {O}rder}.
\newblock Springer-Verlag, 2nd edition, 2001.

\bibitem{GGKSS2019}
A.~D. Gilbert, I.~G. Graham, F.~Y. Kuo, R.~Scheichl, and I.~H. Sloan.
\newblock Analysis of quasi-{M}onte {C}arlo methods for elliptic eigenvalue
  problems with stochastic coefficients.
\newblock {\em Numer. Math.}, 142(4):863--915, 2019.

\bibitem{grisvard1985elliptic}
P.~Grisvard.
\newblock {\em Elliptic Problems in Nonsmooth Domains}.
\newblock Pitman Publishing Inc., 1985.

\bibitem{gk22}
P.~A. Guth and V.~Kaarnioja.
\newblock Generalized dimension truncation error analysis for high-dimensional
  numerical integration: lognormal setting and beyond, 2022.
\newblock Preprint at \url{https://arxiv.org/abs/2209.06176}.

\bibitem{guth22}
P.~A. Guth, V.~Kaarnioja, F.~Y. Kuo, C.~Schillings, and I.~H. Sloan.
\newblock Parabolic {PDE}-constrained optimal control under uncertainty with
  entropic risk measure using quasi-{M}onte {C}arlo integration, 2022.
\newblock Preprint at \url{https://arxiv.org/abs/2208.02767}.

\bibitem{H10}
H.~Harbrecht.
\newblock On output functionals of boundary value problems on stochastic
  domains.
\newblock {\em Math. Methods Appl. Sci.}, 33(1):91--102, 2010.

\bibitem{HKS23}
H.~Harbrecht, V.~Karnaev, and M.~Schmidlin.
\newblock Quantifying domain uncertainty in linear elasticity.
\newblock Technical Report 2023-06, Fachbereich Mathematik, Universit\"at
  Basel, Switzerland, 2023.

\bibitem{nummat2016}
H.~Harbrecht, M.~Peters, and M.~Siebenmorgen.
\newblock Analysis of the domain mapping method for elliptic diffusion problems
  on random domains.
\newblock {\em Numer. Math.}, 134(4):823--856, 2016.

\bibitem{HS21}
H.~Harbrecht and M.~Schmidlin.
\newblock Multilevel quadrature for elliptic problems on random domains by the
  coupling of {FEM} and {BEM}.
\newblock {\em Stoch. Partial Differ. Equ. Anal. Comput.}, 2021.
\newblock to appear.

\bibitem{HSS08}
H.~Harbrecht, R.~Schneider, and Ch. Schwab.
\newblock Sparse second moment analysis for elliptic problems in stochastic
  domains.
\newblock {\em Numer. Math.}, 109(3):385--414, 2008.

\bibitem{inequalities}
G.~H. Hardy, J.~E. Littlewood, and G.~P{\'{o}}lya.
\newblock {\em Inequalities}.
\newblock Cambridge University Press, Cambridge, UK, 1934.

\bibitem{HSSS18}
R.~Hiptmair, L.~Scarabosio, C.~Schillings, and Ch. Schwab.
\newblock Large deformation shape uncertainty quantification in acoustic
  scattering.
\newblock {\em Adv. Comput. Math.}, 44(5):1475--1518, 2018.

\bibitem{hormander}
L.~H{\"{o}}rmander.
\newblock {\em An Introduction to Complex Analysis in Several Variables}.
\newblock North-Holland, 3rd edition, 1990.

\bibitem{JSZ17}
C.~Jerez-Hanckes, Ch. Schwab, and J.~Zech.
\newblock Electromagnetic wave scattering by random surfaces: Shape holomorphy.
\newblock {\em Math. Models Methods Appl. Sci.}, 27(12):2229--2259, 2017.

\bibitem{kks24}
V.~Kaarnioja, F.~Y. Kuo, and I.~H. Sloan.
\newblock Lattice-based kernel approximation and serendipitous weights for
  parametric {PDEs} in very high dimensions.
\newblock To appear in: A.~Hinrichs, P. Kritzer, F. Pillichshammer (eds.).
  Monte Carlo and Quasi-Monte Carlo Methods 2022. Springer Verlag.

\bibitem{periodicpaper}
V.~Kaarnioja, F.~Y. Kuo, and I.~H. Sloan.
\newblock Uncertainty quantification using periodic random variables.
\newblock {\em SIAM J. Numer. Anal.}, 58(2):1068--1091, 2020.

\bibitem{katana}
Katana.
\newblock Published online, 2010.
\newblock DOI:10.26190/669X-A286.

\bibitem{KressnerTobler}
D.~Kressner and C.~Tobler.
\newblock Low-rank tensor {K}rylov subspace methods for parametrized linear
  systems.
\newblock {\em SIAM J. Matrix Anal. Appl.}, 32(4):1288--1316, 2011.

\bibitem{Kuo2012}
F.~Y. Kuo, Ch. Schwab, and I.~H. Sloan.
\newblock Quasi-{M}onte {C}arlo finite element methods for a class of elliptic
  partial differential equations with random coefficients.
\newblock {\em SIAM J. Numer. Anal.}, 50(6):3351--3374, 2012.

\bibitem{mohan}
P.~S. Mohan, P.~B. Nair, and A.~J. Keane.
\newblock Stochastic projection schemes for deterministic linear elliptic
  partial differential equations on random domains.
\newblock {\em Internat. J. Numer. Methods Engrg.}, 85(7):874--895, 2011.

\bibitem{dlmf}
F.~W.~J. Olver, A.~B. {Olde Daalhuis}, D.~W. Lozier, B.~I. Schneider, R.~F.
  Boisvert, C.~W. Clark, B.~R. Miller, B.~V. Saunders, H.~S. Cohl, M.~A.
  McClain, and eds.
\newblock {NIST Digital Library of Mathematical Functions}.
\newblock Release 1.1.6 of 2022-06-30., \url{http://dlmf.nist.gov/}.

\bibitem{savits}
T.~H. Savits.
\newblock Some statistical applications of {F}aa di {B}runo.
\newblock {\em J. Multivariate Anal.}, 97(10):2131--2140, 2006.

\bibitem{sloankachoyan87}
I.~H. Sloan and P.~J. Kachoyan.
\newblock Lattice methods for multiple integration: {T}heory, error analysis
  and examples.
\newblock {\em SIAM J. Numer. Anal.}, 24(1):116--128, 1987.

\bibitem{xiu}
D.~Xiu and D.~M. Tartakovsky.
\newblock Numerical methods for differential equations in random domains.
\newblock {\em SIAM J. Sci. Comput}, 28(3):1167--1185, 2006.

\end{thebibliography}
\end{document}